\documentclass[12pt]{amsart}
\usepackage{longtable}
\usepackage{amsfonts,amssymb,amscd,amsmath,epsfig}
\usepackage{latexsym}
\usepackage{mathrsfs}
\usepackage{tocvsec2}

\usepackage[]{youngtab}

\usepackage[dvipdfm,
hyperindex=true,pagebackref=true,bookmarks=true,
colorlinks=true,linkcolor=blue,citecolor=red]{hyperref}

\begin{document}

\renewcommand{\tilde}{\widetilde}
\renewcommand{\hat}{\widehat}

\newcommand{\BR}{{\mathbb R}}
\newcommand{\BQ}{{\mathbb Q}}
\newcommand{\BC}{{\mathbb C}}
\newcommand{\BP}{{\mathbb P}}
\newcommand{\BZ}{{\mathbb Z}}
\newcommand{\BN}{{\mathbb N}}
\newcommand{\BS}{{\mathbb S}}

\newcommand{\cH}{{\mathcal H}}
\newcommand{\cA}{{\mathcal A}}
\newcommand{\cB}{{\mathcal B}}
\newcommand{\ccF}{{\mathfrak F}}
\newcommand{\cD}{{\mathcal D}}
\newcommand{\cL}{{\mathcal L}}
\newcommand{\cF}{{\mathcal F}}
\newcommand{\cP}{{\mathcal P}}
\newcommand{\cX}{{\mathcal X}}
\newcommand{\cY}{{\mathcal Y}}
\newcommand{\cS}{{\mathcal S}}
\newcommand{\cSol}{\hbox{$\mathcal Sol$}}
\newcommand{\cT}{\hbox{$\mathcal T$}}

\newcommand{\Z}{{\mathbb Z}}
\newcommand{\Q}{{\mathbb Q}}
\newcommand{\N}{{\mathbb N}}
\newcommand{\C}{{\mathbb C}}
\newcommand{\R}{{\mathbb R}}
\newcommand{\X}{{\mathbb X}}
\newcommand{\Y}{{\mathbb Y}}

\newcommand{\CH}{{\mathcal H}}
\newcommand{\CA}{{\mathcal A}}

\def\HH{\mbox{${\mathcal H}$\kern-5.2pt${\mathcal H}$}}

\newcommand{\binomial}[2]{\genfrac{(}{)}{0pt}{}{ #1 }{ #2 }}
\newcommand{\qbinomial}[2]{\genfrac{[}{]}{0pt}{}{ #1 }{ #2 }_q }
\newcommand{\qbinom}[3]{\genfrac{[}{]}{0pt}{}{ #1 }{ #2 }_{ #3 } }


\def\der{\partial}
\def\tensor{\otimes}
\def\gam{\gamma} \def\Gam{\Gamma}
\def\del{\delta} \def\Del{\Delta}
\def\kap{\kappa}
\def\lam{\lambda} \def\Lam{\Lambda}
\def\Comp{{\mathbb C}}
\def\sM{{\mathcal M}}

\newtheorem{theorem}{Theorem}[section]
\newtheorem{maintheorem}[theorem]{Main Theorem}
\newtheorem{proposition}[theorem]{Proposition}
\newtheorem{definition}[theorem]{Definition}
\newtheorem{lemma}[theorem]{Lemma}
\newtheorem{corollary}[theorem]{Corollary}
\newtheorem{notation}[theorem]{Notation}
\newtheorem{remark}[theorem]{Remark}
\newtheorem{example}[theorem]{Example}

\newtheorem{theorem }{Theorem}[section]
\newtheorem{maintheorem }[theorem]{Main Theorem}
\newtheorem{proposition }[theorem]{Proposition}
\newtheorem{definition }[theorem]{Definition}
\newtheorem{lemma }[theorem]{Lemma}
\newtheorem{corollary }[theorem]{Corollary}
\newtheorem{notation }[theorem]{Notation}
\newtheorem{remark }[theorem]{Remark}
\newtheorem{example }[theorem]{Example}

\newtheorem{ maintheorem }[theorem]{Main Theorem}
\newtheorem{ theorem}{Theorem}[section]
\newtheorem{ proposition}[theorem]{Proposition}
\newtheorem{ definition}[theorem]{Definition}
\newtheorem{ lemma}[theorem]{Lemma}
\newtheorem{ corollary}[theorem]{Corollary}
\newtheorem{ notation}[theorem]{Notation}
\newtheorem{ remark}[theorem]{Remark}
\newtheorem{ example}[theorem]{Example}

\newtheorem{thm}{Theorem}[section]
\newtheorem{prop}[thm]{Proposition}
\newtheorem{lem}[thm]{Lemma}
\newtheorem{cor}[thm]{Corollary}
\newtheorem{conj}[thm]{Conjecture}
\newtheorem{con}[thm]{Conjecture}
\newtheorem{dfn}[thm]{Definition}
\newtheorem{df}[thm]{Definition}
 \newcommand{\rem}{{\bf Comment.\ }}
 \newcommand{\rmk}{{\bf Comment.\ }}
 \newcommand{\exmp}{{\bf Example.\ }}
 \newcommand{\ex}{{\bf Example.\ }}
 \newcommand{\prob}{{\bf Problem.\ }}

\newtheorem{note}{Note} 
\renewcommand{\thenote}{}
\newtheorem*{acka}{Acknowledgments}
\newtheorem{ack}{Acknowledgments}
\renewcommand{\theack}{}
\renewcommand{\appendixname}{\bf Appendix}
\renewcommand{\proof}{{\em Proof.\ }}

\hyphenation{
ap-pen-dix as-ymp-tot-ic at-trib-uted at-trib-ut-able
Bry-li-n-sky com-mu-ta-tion de-ge-ne-rate
de-riv-a-tive dis-trib-ute equi-vari-ant ex-tra-or-di-nary  
geo-met-ric griev-ance griev-ous grad-ed ho-lo-no-my ho-mo-thetic
in-fin-ite-ly in-fin-i-tes-i-mal Ha-rish Cha-n-dra mul-ti-plic-able 
non-euclid-ean non-iso-mor-phic non-smooth par-a-digm 
par-a-bol-ic pa-rab-o-loid pa-ram-e-trize phe-nom-e-non 
post-script pseu-do-dif-fer-en-tial pseu-do-fi-nite 
qua-drat-ics quad-ra-ture Han-kel rec-tan-gle semi-def-i-nite 
set-up wide-spread Euler-ian Feb-ru-ary Gauss-ian Grothen-dieck 
Hamil-ton-ian Her-mi-t-ian her-mi-t-ian Jan-u-ary 
Japan-ese Ka-shi-wa-ra Kor-te-weg Le-gendre No-vem-ber Rie-mann-ian 
Sep-tem-ber Za-mo-lo-d-chi-kov Kni-zh-nik quan-tum Op-dam
Mac-do-nald Ca-lo-ge-ro Su-ther-land Mo-ser 
Ol-sha-net-sky  Pe-re-lo-mov in-de-pen-dent ope-ra-tors 
cy-clo-to-mic ra-tio-nal de-gen-er-a-tion 
in-ter-est-ing de-for-ma-tions de-for-ma-tion pro-ce-dure 
fol-lows ope-ra-tors  pre-serve suf-fices ap-proach 
for-mu-las con-sider its com-ple-tion cor-re-spond-ing 
au-to-mor-phism be-cause pro-por-tional fi-nal-ly let-ting 
equi-v-a-lence ge-n-er-al-ized Mac-do-nald iden-ti-ties 
cor-re-s-pond sub-dia-grams par-ti-tion na-t-u-ral-ly 
or-dered stan-dard de-for-ma-tion ar-gu-ment com-bined 
sphe-r-i-cal rep-re-sen-ta-tions tri-go-no-me-t-ric
ge-n-er-al-ly speak-ing pri-m-it-ive ir-re-du-cible 
sum-ma-tion  rep-re-sen-ta-tives pro-por-ti-o-na-li-ty
ultra-sphe-ri-cal Ro-gers}

\def\ffor{\quad\hbox{ for }\quad}
\def\wwhen{\quad\hbox{ when }\quad}
\def\wwhere{\quad\hbox{ where }\quad}
\def\aand{\quad\hbox{ and }\quad}
\def\for{\  \hbox{ for } \ }
\def\iif{ \ \hbox{ if } \ }
\def\when{ \ \hbox{ when } \ }
\def\where{\  \hbox{ where } \ }
\def\and{\  \hbox{ and } \ }
\def\and{\  \hbox{ and } \ }
\def\oor{\  \hbox{ or } \ }
\def\proof{{\em Proof. \  }}

\def\equal{\stackrel{\,\mathbf{def}}{= \kern-3pt =}}

\def\la{\lambda}
\def\La{\Lambda}
\def\om{\omega}
\def\Om{\Omega}
\def\Th{\Theta}
\def\th{\theta}
\def\al{\alpha}
\def\be{\beta}
\def\ga{\gamma}
\def\ep{\epsilon}
\def\up{\upsilon}
\def\Up{\Upsilon}
\def\de{\delta}
\def\De{\Delta}
\def\ka{\kappa}
\def\kapp{\hbox{\bf \ae}}
\def\si{\sigma}
\def\Si{\Sigma}
\def\Ga{\Gamma}
\def\ze{\zeta}
\def\io{\iota}
\def\bio{b^\iota}
\def\aio{a^\iota}
\def\twio{\tilde{w}^\iota}
\def\hwio{\hat{w}^\iota}
\def\gio{\g^\iota}
\def\Bio{B^\iota}

\def\del{\delta}
\def\pa{\partial}
\def\vp{\varphi}
\def\ve{\varepsilon}
\def\inf{\infty}

\def\vph{\varphi}
\def\vps{\varpsi}
\def\vPh{\varPhi}
\def\vep{\varepsilon}
\def\vpi{{\varpi}}
\def\vth{{\vartheta}}
\def\vsi{{\varsigma}}
\def\vrh{{\varrho}}

\def\bph{\bar{\phi}}
\def\bsi{\bar{\si}}
\def\bvp{\bar{\varphi}}

\newcommand{\bS}{{\mathbf S}}
\newcommand{\bH}{{\mathbf H}}
\newcommand{\bF}{{\mathbf F}}
\newcommand{\bE}{{\mathbf E}}

\def\tal{\tilde{\alpha}}
\def\tbe{\tilde{\beta}}
\def\tde{\tilde{\delta}}
\def\tpi{\tilde{\pi}}
\def\txi{\tilde{\xi}}
\def\tPi{\tilde{\Pi}}
\def\tPhi{\tilde{\Phi}}
\def\tV{\tilde{V}}
\def\tJ{\tilde{J}}
\def\tla{\tilde{\lambda}}
\def\tga{\tilde{\gamma}}
\def\tGa{\tilde{\Gamma}}
\def\tvs{\tilde{{\varsigma}}}
\def\tu{\tilde{u}}
\def\tU{\tilde{U}}
\def\tw{\widetilde w}
\def\tW{\widetilde W}
\def\tB{\tilde B}
\def\tv{\tilde v}
\def\tV{\tilde V}
\def\tz{\tilde z}
\def\tb{\tilde b}
\def\ta{\tilde a}
\def\tih{\tilde h}
\def\trh{\tilde {\rho}}
\def\tx{\tilde x}
\def\tf{\tilde f}
\def\tg{\tilde g}
\def\tG{\tilde G}
\def\tk{\tilde k}
\def\tl{\tilde l}
\def\tL{\tilde L}
\def\tD{\tilde D}
\def\tR{\tilde R}
\def\tP{\tilde P}
\def\tH{\tilde H}
\def\tp{\tilde p}

\def\hH{\hat{H}}
\def\hh{\hat{h}}
\def\hR{\hat{R}}
\def\hY{\hat{Y}}
\def\hX{\hat{X}}
\def\hP{\hat{P}}
\def\hT{\hat{T}}
\def\hV{\hat{V}}
\def\hG{\hat{G}}
\def\hF{\hat{F}}
\def\hw{\widehat{w}}
\def\hW{\widehat{W}}
\def\hu{\hat{u}}
\def\hs{\hat{s}}
\def\hv{\hat{v}}
\def\hb{\hat{b}}
\def\hB{\widehat{B}}
\def\hze{\hat{\zeta}}
\def\hsi{\hat{\sigma}}
\def\hrh{\hat{\rho}}
\def\hth{\hat{\theta}}
\def\hy{\hat{y}}
\def\hx{\hat{x}}
\def\hz{\hat{z}}
\def\hg{\hat{g}}
\def\he{\hat{e}}
\def\hE{\widehat{E}}

\def\B{\mathbf{B}}
\def\I{\mathbf{I}}
\def\P{\mathbf{P}}
\def\G{\mathbf{G}}
\def\S{\mathbf{S}}
\def\F{\mathbf{F}}
\def\one{\mathbf{1}}
\def\Sn{\mathbf{S}_n}
\def\0{\mathbf{0}}
\def\H{\mathbf{H}}
\def\V{\mathbf{V}}

\def\f{\mathcal{F}}
\def\çF{\mathcal{F}}
\def\o{\mathcal{O}}
\def\t{\mathcal{T}}
\def\r{\mathcal{R}}
\def\l{\mathcal{L}}
\def\m{\mathcal{M}}
\def\k{\mathcal{K}}
\def\n{\mathcal{N}}
\def\d{\mathcal{D}}
\def\p{\mathcal{P}}
\def\cP{\mathcal{P}}
\def\a{\mathcal{A}}
\def\h{\mathcal{H}}
\def\c{\mathcal{C}}
\def\y{\mathcal{Y}}
\def\e{\mathcal{E}}
\def\v{\mathcal{V}}
\def\z{\mathcal{Z}}
\def\x{\mathcal{X}}
\def\s{\mathcal{S}}
\def\g{\mathcal{G}}
\def\u{\mathcal{U}}
\def\w{\mathcal{W}}
\def\i{\mathcal{I}}
\def\j{\mathcal{J}}
\def\b{\mathcal{B}}

\def\lan{\langle}
\def\llb{(\!(}
\def\ran{\rangle}
\def\rrb{)\!)}
 \def\dim{{\hbox{\rm dim}}_{\mathbb C}\,}
\def\lng{\hbox{\rm{\tiny lng}}}
\def\sht{\hbox{\rm{\tiny sht}}}
\def\sph{\hbox{\rm{\tiny sph}}}
\def\inv{\hbox{\rm{\tiny inv}}}

\def\br#1{\langle #1 \rangle}

\def\rank{\hbox{rank}}
\def\gl{\mathfrak{gl}_N}

\newcommand{\Aut}{\operatorname{Aut}}
\newcommand{\Hom}{\operatorname{Hom}}
\newcommand{\End}{\operatorname{End}}
\newcommand{\Ind}{\operatorname{Ind}}
\newcommand{\ad}{\operatorname{ad}}
\newcommand{\pr}{\operatorname{pr}}
\newcommand{\aweyl}{\tilde{\mathbb S}_n}
\newcommand{\hec}{{\mathcal H}^t_n}
\newcommand{\Func}{{\mathcal F}({\mathbb C}^n,{\mathcal H}^t_n)}
\newcommand{\tr}{\operatorname{tr}}
\newcommand{\Out}{\operatorname{Out}}
\newcommand{\Rad}{\operatorname{Rad}}
\newcommand{\Spec}{\operatorname{Spec}}
\newcommand{\id}{\operatorname{id}}
\newcommand{\Int}{\operatorname{Int}}
\newcommand{\ct} {\operatorname{ct}}

\newcommand{\rat}{{\mathbb Q}}
\newcommand{\real}{{\mathbb R}}
\newcommand{\cplx}{{\mathbb C}}
\newcommand{\zint}{{\mathbb Z}}

\newcommand{\sq}{\phantom{1}\hfill$\qed$}
\newcommand{\Rea}{\Re}
\newcommand{\Ima}{\Im}

\newcommand{\st}{\bowtie}
\newcommand{\modd}{\mbox{\,mod\,}}
\newcommand{\lr}{\langle}
\newcommand{\rr}{\rangle}
\newcommand{\eps}{\varepsilon}
\newcommand{\phk}{\phi^{(k)}}
\newcommand{\psk}{\psi^{(k)}}
\newcommand{\Res}{\mbox{Res}\;}
\newcommand{\sgn}{\mbox{sgn}}
\newcommand{\mn} {\left\{ \begin{array}{c}m\\
n\end{array}\right\}}

\def\sX{\mathscr{X}}
\def\sH{\mathscr{H}}
\def\sY{\mathscr{Y}}
\def\TT{\mathfrak{T}}
\def\JJ{\mathfrak{J}}
\def\HH{\mathfrak{H}}
\def\FF{\mathfrak{F}}
\def\GG{\mathfrak{G}}
\def\CC{\mathfrak{C}}
\def\LL{\mathfrak{L}}

\def\BB{\mathfrak{B}}
\def\AA{\mathfrak{A}}
\def\ZZ{\mathfrak{Z}}
\def\HH{\hbox{${\mathcal H}$\kern-5.2pt${\mathcal H}$}}
\def\HHH{\hbox{${\mathbb H}$\kern-4.2pt${\mathbb H}$}}
\def\tHH{\widetilde{\HH\ }}

\font\smm=msbm10 at 12pt 
\def\symbol#1{\hbox{\smm #1}}
\def\lsmash{{\symbol n}}
\def\rsmash{{\symbol o}}
\def\#{\sharp}

\font\tenbf=cmbx10
\font\tenrm=cmr10
\font\tenit=cmti10
\font\ninebf=cmbx9
\font\ninerm=cmr9
\font\nineit=cmti9
\font\eightbf=cmbx8
\font\eightrm=cmr8
\font\eightit=cmti8
\font\sevenrm=cmr7
\font\sevenbf=cmbx7


\par
{\centering
Dedicated with admiration to Yuri Ivanovich Manin  \\
on the occasion of his 80th birthday
\medskip
\par}

\title [Riemann Hypothesis for DAHA superpolynomials]
{Riemann Hypothesis for DAHA superpolynomials
and plane curve singularities}
\author[Ivan Cherednik]{Ivan Cherednik $^\dag$}

\begin{abstract}
Stable Khovanov-Rozansky polynomials of algebraic knots are
expected to coincide with certain generating functions,
superpolynomials, of nested Hilbert schemes and
flagged Jacobian factors of the corresponding plane curve
singularities. Also, these 3 families conjecturally match
the DAHA superpolynomials. These superpolynomials
can be considered as singular counterparts and 
generalizations of the Hasse-Weil zeta-functions. 
We conjecture that all $a$-coefficients of the
DAHA superpolynomials upon the substitution
$q\mapsto qt$ satisfy the Riemann Hypothesis for
sufficiently small $q$ 
for uncolored algebraic knots, presumably for
$q\le 1/2$ as $a=0$. This can be partially extended to
algebraic links at least for $a=0$. Colored links
are also considered, though mostly for rectangle Young diagrams. 
Connections with Kapranov's motivic zeta 
and the Galkin-St\"ohr zeta-functions are discussed.
\end{abstract}

\thanks{$^\dag$ \today.
\ \ \ Partially supported by NSF grant
DMS--1363138}

\address[I. Cherednik]{Department of Mathematics, UNC
Chapel Hill, North Carolina 27599, USA\\
chered@email.unc.edu}

 \def\sht{\raisebox{0.4ex}{\hbox{\rm{\tiny sht}}}}
 \def\bysame{{\bf --- }}
\let\oldt\~   
 \def\~{{\bf --}}
 \def\rr{{\mathsf r}}
 \def\ss{{\mathsf s}}
 \def\mm{{\mathsf m}}
 \def\pp{{\mathsf p}}
 \def\ll{{\mathsf l}}
 \def\aa{{\mathsf a}}
 \def\bb{{\mathsf b}}
 \def\NS{\hbox{\tiny\sf ns}}
 \def\ssum{\hbox{\small$\sum$}}
\newcommand{\comment}[1]{}
\renewcommand{\tilde}{\widetilde}
\renewcommand{\hat}{\widehat}
\renewcommand{\V}{\mathbb{V}}
\renewcommand{\S}{\mathbb{S}}
\renewcommand{\F}{\mathbb{F}}
\newcommand{\dagx}{\hbox{\tiny\mathversion{bold}$\dag$}}
\newcommand{\ddagx}{\hbox{\tiny\mathversion{bold}$\ddag$}}
\newtheorem{conjecture}[theorem]{Conjecture}
\newcommand*\toeq{
\raisebox{-0.15 em}{\,\ensuremath{
\xrightarrow{\raisebox{-0.3 em}{\ensuremath{\sim}}}}\,}
}
\newcommand{\unknot}{\hbox{\tiny\!\raisebox{0.2 em}
{$\bigcirc$}}\!}
\newcommand{\mmu}{\hbox{\mathversion{bold}$\mu$}}
\newcommand{\lla}{\hbox{\mathversion{bold}$\lambda$}}
\newcommand{\dde}{\hbox{\mathversion{bold}$\delta$}}

\newcommand\rightthreearrow{\hbox{\tiny
        $\mathrel{\vcenter{\mathsurround0pt
         \ialign{##\crcr
         \noalign{\nointerlineskip}$\rightarrow$\crcr
         \noalign{\nointerlineskip}$\rightarrow$\crcr
         \noalign{\nointerlineskip}$\rightarrow$\crcr
                }}}$ }}
\newcommand\rightfourarrow{\hbox{\tiny
        $\mathrel{\vcenter{\mathsurround0pt
         \ialign{##\crcr
         \noalign{\nointerlineskip}$\rightarrow$\crcr
         \noalign{\nointerlineskip}$\rightarrow$\crcr
         \noalign{\nointerlineskip}$\rightarrow$\crcr
         \noalign{\nointerlineskip}$\rightarrow$\crcr
                }}}$ }}

\newcommand\rightdotsarrow{\hbox{\small
        $\mathrel{\vcenter{\mathsurround0pt
         \ialign{##\crcr
         \noalign{\nointerlineskip}$\,\rightarrow$\crcr
         \noalign{\nointerlineskip}$\cdots$\crcr
         \noalign{\nointerlineskip}$\,\rightarrow$\crcr
                }}}$ }}
\newcommand\rightdotsarrowtiny{\hbox{\tiny
        $\mathrel{\vcenter{\mathsurround0pt
         \ialign{##\crcr
         \noalign{\nointerlineskip}$\,\rightarrow$\crcr
         \noalign{\nointerlineskip}$\cdots$\crcr
         \noalign{\nointerlineskip}$\,\rightarrow$\crcr
                }}}$ }}

\newcommand*{\vect}[1]{\overrightarrow{\mkern0mu#1}}
\newcommand{\twoone}
{\hbox{\rm
$\circ\!$\raisebox{-2.6pt}{$\rightarrow$}
\raisebox{-2.6pt}{$\!\!\circ\!\!\rightarrow$}
\kern-33pt\raisebox{+2.6pt}{$\rightarrow$}
\kern+14pt}}
\newcommand{\twoonetiny}
{\hbox{\tiny
$\circ\!$\raisebox{-2.pt}{$\rightarrow$}
\raisebox{-2.pt}{$\!\!\circ\!\!\rightarrow$}
\kern-23pt\raisebox{+2.pt}{$\rightarrow$}
\kern+10pt}}

\newcommand{\twotwo}
{\hbox{\rm $\circ\!\rightrightarrows\!
\raisebox{2.5pt}{\hbox{\small $\circ$}}
\kern-5.5pt\raisebox{-2.5pt}{\hbox{\small $\circ$}}
\!\rightrightarrows\,$}}
\newcommand{\twotwotiny}
{\hbox{\tiny $\circ\!\rightrightarrows\!
\raisebox{2.pt}{\hbox{\tiny $\circ$}}
\kern-4.2pt\raisebox{-2.pt}{\hbox{\tiny $\circ$}}
\!\rightrightarrows$}}
\newcommand{\tax}{\hbox{\sf[r,s]}}
\newcommand{\lxi}{\raisebox{0.5pt}{${}^\xi$}\!}
\vskip -0.0cm

\maketitle
\smallskip

{\tiny
\centerline{{\bf MSC} (2010): 14H50, 17B22, 17B45, 20C08,
20F36, 33D52, 30F10, 55N10, 57M25}
}
\smallskip

\vskip -0.5cm
\renewcommand{\baselinestretch}{1.2}
{\vbadness=10000\textmd
{\small \tableofcontents}
}
\renewcommand{\baselinestretch}{1.0}
\vfill\eject

\section*{List of basic notations}\label{SEC:BASIC}
{\small
\begin{align*}
& R\!=\!\{\al\}\!\subset\!\R^{n+1},\ \al_i\!=\!\ep_i\!-\!\ep_{i+1}
&\text{root system of type $A_n$, simple roots}&  \\
&W=\S_{n+1},\ \, P,Q,\ \,\Pi=P/Q
&\text{the Weyl group, weight/root lattices}& \\
&\tR\!=\!\{\tal\!=\![\al,j]\},\ \ \al_0\!=\![-\th,1]
&\text{affine root system, $\th=$maximal root}& \\
&\tW\!\!=\!\!\lan s_{\tal}\ran,\, \hW\!\!=\!W\lsmash 
P\!=\!\tW\rsmash\Pi
&\text{affine,\ \, extended Weyl groups:  }(\ref{ondthr})&\\
&\HH=\lan X_b,Y_b,T_i,\,q^{\pm\frac{1}{n+1}},
t^{\pm\frac{1}{2}}\ran
&\text{DAHA: Definition \ref{double}, $Y_{b}$
\tiny{$\!(\!b\!\in\!P\!)$}: }
(\ref{Ybx})& \\
&\tau_{\pm},\si\!=\!\tau_+\tau_-^{-1}\tau_+, 
\vph \text{\,: (\ref{starphi})}, \star
&\text{automorphisms \& involutions: } (\ref{tauplus})&\\
&\text{$P_b\!=\!$ Macdonald polynomials } 
&\text{$\v=$ polynomial module: Sect. } \ref{SEC:POLR}&\\
&\{H\}_{ev}=H(1)(t^{-\rho}),\ H\in \HH
&\text{coinvariant, $H(1)=H\!\Downarrow\ \in \v$: } 
(\ref{evfunct})&\\
&P_b^\circ\,=\, P_b(X)/(P_b(t^{-\rho}), b\in P_+
&\text{spherical normalization of $P_b$: } (\ref{Pbrho'})&\\
&J_\la= h_\la P_b,\ \,\la=\la(b),\ \, b\in P_+\,
&\text{$J$-polynomials for diagrams $\!\la$: }
 (\ref{P-arms-legs})&\\
&h_\la\!=\!\prod_{\Box\in\la}
(1-q^{arm(\Box)}t^{leg(\Box)+1})
&\text{arm and leg numbers in : Sect. } \ref{SEC:JPOL}&\\
&\vec{\rr}=\{\rr_1,\ldots \rr_\ell\}, \,
\vec{\ss}=\{\ss_1,\ldots \ss_\ell\}
&\text{Newton's pairs,\ \, $gcd(\rr_i,\ss_i)\!=1$\,: }
(\ref{iterrss})&\\
&\aa_1=\ss_1,\ \,\, \aa_{i}=\aa_{i-1}\rr_{i-1}\rr_{i}+\ss_{i}\,\
&\text{cable (topological)  parameters\,: } (\ref{Newtonpair})&\\
&\l\,=\,\l_{(\vec\rr^j,\vec\ss^j)}^{\Up,(b^j)},\ \,\ 
'\!\l\,=\, '\!\l_{('\vec\rr^j,'\vec\ss^j)}^{\,'\Up,('\!b^j)}
&\text{pairs of labeled/colored graphs: }
(\ref{links-x-y})&\\
&\hat{J\!D}^{min}_{(\vec\rr^j,\,\vec\ss^j)
,('\!\vec\rr^j,'\!\vec\ss^j)}
((b^j),('\!b^j);q,t)
&\text{hat-DAHA-Jones polynomials: }
(\ref{jones-bar})&\\
&\hat{\h}^{min}_{\l,\,'\!\l}(q,t,a)=
\sum_{i=0}^d \h^i(q,t) a^i
&\text{superpolynomials,  $d\!=\!deg_a(\hat{\h})$: }
(\ref{h-polynoms-hat})&\\
&\hat{H}^i(q,t)\!=\!\h^i(qt,t)\hat{}\,,\ \varsigma_i,\pi_i,S^i
&\text{hat-normalized $\h^i(qt,t)$: Conject. } \ref{CONJLIM}&\\
&\hat{H}^i_{sym}\!=\!(H^i(q,t)\!\!+\!\! 
H^i(q,\frac{1}{qt}))\,\hat{{}}
&\text{hat-symmetrization of $H^i(q,t)$: }
(\ref{symhpol})&\\
&\hat{\mathbb{H}}(q,t,a,u)\!=\!\!\!
\sum_{m=0}^\infty\!\! \left(\frac{u}{t}
\right)^{g(m)}\!\hat{\h}_m
&\text{family superpolynomials 
\tiny{$(\tau_-^m\!\ga_1)$}\,: }
(\ref{hsfam})&\\
&\h_{_{mot}}(q,t,a),\h_{_{mot}}^0\!\!\!=
\!\!\h_{_{mot}}\!(\hbox{\footnotesize{a\!=\!0}})
&\text{ flagged\, motivic\, s-\,polynomials: }
(\ref{motpolfull})&\\
&\z(q,t,a),\l(q,t,a),\ Z,L(q,t)  
&\text{\small
flagged Galkin-St\"ohr functions: }
(\ref{zaqt})&\\ 
&\text{given $\hat{H}^i$},\ \, \varpi_i=\min(\om') 
\text{\, s.t.}
&\text{{
 (weak) RH holds for} $\om\!=\!\frac{1}{q}\!>\!\om'$: } 
(\ref{varpi-noprime})&\\ 
&\text{\tiny
$\{\ga[3,2]\ga[2,1](P)\}\!=\!
\{\hat{\ga}_{3,2}(\hat{\ga}_{2,1}(P)\!\Downarrow)\}$}
&\text{coinvariants' abbreviations:  Sect.\,} \ref{SEC:ABREV}&\\
\end{align*}
}

\renewcommand{\natural}{\wr}

\setcounter{section}{0}
\setcounter{equation}{0}
\section{\sc Introduction}
\vbadness=3000
\hbadness=3000
The aim of the paper is to approach the Riemann Hypothesis,
$RH$, for DAHA superpolynomials of algebraic links colored 
by Young diagrams  upon the substitution 
$q\mapsto q t$.
The parameter $q$, a counterpart of the cardinality 
of $\mathbb F$ in the Weil conjectures, is assumed
sufficiently small, which is complementary to the classical theory.
Then $RH$ presumably holds for any $a$\~coefficients
of DAHA superpolynomials of {\em uncolored algebraic 
knots\,}; moreover, $q\le 1/2$ seems sufficient when $a\!=\!0$. 
For links, stable (any $q$) {\em irregular\,} (non-$RH$) 
zeros appear. For instance, the number of their pairs  
is conjectured to coincide with the number of 
components of uncolored algebraic {\em links\,} minus $1$ 
as $a\!=\!0$. We provide tools for finding
such bounds for any Young diagrams and arbitrary 
$a$\~coefficients; finding the exact $RH$\~range of $q$ is
much more subtle,  which is somewhat parallel to 
the theory of {\em spectral zeta-functions\,}. 


Let us try to put this conjecture into perspective and explain
the rationale behind it and its relations to 
the classical Weil conjectures. 

\subsection{\bf Superpolynomials}
\subsubsection{\sf Topological and geometric theories}
The superpolynomials have several reincarnations
in mathematics and physics; the origin is the theory
of {\em stable Khovanov-Rozansky polynomials\,}, which are
Poincar\'e polynomials of the HOMFLY-PT
triply-graded link homology \cite{Kh,KhR1}. They depend on
3 parameters $q,t,a$ and are actually infinite series in
the {\em unreduced\,} case. Generally, they are difficult
to calculate and there are
unsettled problems with the formulas for links,
in the presence of colors and in the reduced case, 
though the {\em categorification
theory\,} generally provides their definition for any 
colors (dominant weights). 
For uncolored unreduced algebraic links, they were conjectured
to coincide with the {\em ORS superpolynomials\,},
certain generating series
for {\em nested Hilbert schemes\,} of the corresponding
{\em plane curve singularties\,} \cite{ORS}.

These two families in the reduced setting are
conjecturally related to the {\em geometric superpolynomials\,}
introduced in \cite{ChP1,ChP2}. They were defined there
for any algebraic knots colored by columns (wedge powers of the
fundamental representation), developing \cite{ChD1,GM,Gor}.
Their construction is in terms of the {\em flagged
Jacobian factors\,}
of unibranch plane curve singularities.
Jacobian factors are
(indeed) factors of the corresponding
{\em compactified Jacobians\,};
the definition is entirely local. In turn, Jacobian
factors are almost directly
related to the affine Springer fibers in type $A$ (the
nil-elliptic case), and
therefore to the corresponding $p$\~adic orbital integrals;
see the end of \cite{ChP1} for some references and discussion.

Inspired by \cite{ORS,GORS}, \cite{ChP1,ChP2} and various prior
works, especially \cite{Kap,GSh}, the 
{\em geometric superpolynomials\,} can be considered
as ``singular" analogs of the {\em Weil polynomials\,},
the numerators of the Hasse-Weil zeta-functions of smooth curves.
To analyze this we switch to the $4${\tiny th} type
of superpolynomials, the DAHA ones from
\cite{CJ,CJJ,GoN} and further works;
the most comprehensive paper on them by now is \cite{ChD2}.

\subsubsection{\sf DAHA superpolynomials}
The DAHA superpolynomials deal with the combinatorial data of
iterated torus links and allow any colors (Young
diagrams). Importantly, they almost directly reflect the
{\em topological\,} type
of singularity, in contrast to the ORS construction and the
geometric superpolynomials from \cite{ChP1,ChP2}. We note that
the latter are related to {\em restricted\,} nested
Hilbert schemes of singularities, some subvarieties of those 
used in \cite{ORS} (geometrically simpler). See
Section \ref{SEC:HILB}.
The DAHA superpolynomials  are expected to be 
connected with physics superpolynomials
\cite{DGR,AS,GS,DMMSS,FGS}, which is not discussed 
in this work.

They are defined so far for 
{\em iterated torus links\,}. All initial {\em intrinsic\,}
conjectures from \cite{CJ} concerning torus knots
were proved (but the {\em positivity\,} discussed
below) and extended to any colored iterated torus links
(not only algebraic). The theory of DAHA-Jones polynomials
is very much uniform for any colors (dominant weights)
and root systems; the DAHA superpolynomials are in type $A$.
\smallskip

In other approaches, the limitations and practical problems
are more significant, especially with links and colors.
The ORS polynomials are generally difficult to calculate
since they are based on the {\em weight filtration\,} in
cohomology of the corresponding nested Hilbert schemes.
The {\em KhR polynomials\,} are known only for simple knots.
For torus knots, they were recently
calculated using Soergel bimodules \cite{Mel2};
see also 
Corollary 3.4 there, which proves Conjecture 2.7 $(ii)$ (uncolored)
from \cite{CJ}; not all details are  provided in  \cite{Mel2},
but the proof seems essentially a direct identification of
Gorsky's combinatorial formulas with those in \cite{CJ}.
See also \cite{EH,Hog} concerning $T(m\rr\pm 1,\rr)$ and 
links $T(m\rr,\rr)$.
The geometric superpolynomials and {\em flagged Jacobian
factors\,} from \cite{ChP1,ChP2}
are relatively simple to define, but 
this is done so far only for algebraic  links
colored by columns.

\subsection{\bf DAHA and Weil conjectures}
Let us present the main features of the theory of DAHA
superpolynomials as analogs of the Weil conjectures.
This connection is mostly heuristic and our $RH$ mostly
serves the sector of $q<1$, complementary to Weil's $RH$. The
key point of the whole work is the identification of
$\,t\,$ in DAHA superpolynomials with $T$ in the singular
counterpart of the Weil 
zeta-function.

\subsubsection{\sf Polynomiality and super-duality}
First of all, the DAHA theory and the
geometric construction from \cite{ChP1,ChP2}
{\em directly\,} provide superpolynomials,
counterparts of the Weil $P_1(T)$.
This is in contrast to the classical
theory, where $P_1(T)$ appears due to the
rationality theorem:
$$\ze(X,T)\equal
\exp(\,\sum_{n=1}^\infty\,T^n\,|X(\mathbb F_{q^n})|\,/\,n\,)=  
\frac{P_1(T)}{(1-T)(1-qT)},
$$
where $X$ is a smooth projective curve over
a finite field $\mathbb F=\mathbb F_q$ with $q$ elements and $P_1$
is of degree $2g$ for the genus $g$ of $X$ (its smoothness is 
actually needed only for $RH$). 

The coefficients of DAHA superpolynomials are presumably
all positive for rectangle Young diagrams and
algebraic knots (the positivity
conjectures from \cite{CJ,ChD1}), which is not present
in the Weil-Deligne theory \cite{Del1,Del2}.  Such a positivity 
hints at
a possible geometric interpretation and
``categorification" of these polynomials. However, they can be
positive for ``very" non-algebraic knots.
For links and non-rectangle 
Young diagrams, the positivity holds only upon the division by
some powers of $(1-t), (1-q)$. For instance,
$(1-t)^{\kappa-1}$ is presumably sufficient
for uncolored links with $\kappa$ components.

To avoid misunderstanding, let us emphasize that the
positivity of DAHA superpolynomials for knots
(and for links upon the division above) is neither
necessary nor sufficient for the validity of the
corresponding $RH$. However
when such positivity holds, one may expect
a geometric interpretation of DAHA superpolynomials as
in \cite{ChP1,ChP2}, which actually implies $RH$ 
for sufficiently small $q$.


Thus, a counterpart
of the existence of $P_1(T)$ is the
polynomiality for
{\em DAHA-Jones polynomials}\, from 
Theorem 1.2,\cite{CJJ} (for torus knots, any colors,
and root systems)
and its generalizations to iterated torus links. 
The passage to the superpolynomials,
Theorem 1.3 there, was announced in \cite{CJ}
(based on \cite{SV}); its complete proof
was provided in \cite{GoN}. See \cite{ChD2} for the most 
general version. We provide here many examples of 
algebraic links and the corresponding DAHA procedures,
significantly developing and extending those from 
\cite{ChD1,ChD2}. 
\vskip 0.1cm

Accordingly, the {\em super-duality\,} of DAHA superpolynomials 
matches {\em the Weil's
functional equation}: $\ze(X,q^{-1}T^{-1})=
q^{1-g}T^{2-2g}\ze(X,T)$.
 It was
conjectured in \cite{CJ} (let us mention prior \cite{GS}
in the context of physics superpolynomials)
and proved in \cite{GoN} on the basis of the
$q\leftrightarrow t$\~duality of the
modified Macdonald polynomials.
An alternative approach to the proof via roots of unity and the
generalized {\em level-rank
duality\,} was presented in \cite{CJJ}; it can work for classical
root systems and some other families. The proof
of the duality from Proposition 3 from \cite{ORS} is
parallel
to that of the motivic functional equation 
\cite{Kap}; see also Section 6 in \cite{Gal},
Section 3 in \cite{Sto} and formula (\ref{feqze}) below
concerning the functional equation for the
Galkin-St\"ohr functions.

We note that our parameter $a$\, and adding colors 
(numerically, we mostly consider rectangle Young
diagrams) do not have direct 
origins in the theory  of Hasse-Weil zeta-functions. 
The parameter $a$ is associated with {\em flagged\,} 
Jacobian factors in \cite{ChP1,ChP2}; we also
define flagged Galkin-St\"ohr functions by considering
{\em standardizable\,} flags of ideals.


\vskip -1.5cm
\subsubsection{\sf Riemann Hypothesis}
For large $q$, corresponding to the
cardinality $|\mathbb F|$, 
$RH$ for our superpolynomials
and the zeta-functions from \cite{Gal,Sto} generally fails. 
However the inequality $q\le 1/2$ is (surprisingly) sufficient  
for uncolored {\em algebraic\,} knots at $a\!=\!0$ 
we considered; $\mathbb F[[z^4,z^6+z^{7+2m}]]$
presumably make $1/2$ sharp here as $m\to\infty$.
For uncolored algebraic {\em links\,}, the number of
stable pairs of irregular zeros is conjectured to be  
the number of components minus $1$ as $a=0$. Adding colors 
to knots and links is more subtle, though rectangle Young diagrams
satisfy  $RH$ for sufficiently small $q$ upon the symmetrization
at least for $a=0$ in all examples we considered.
We note that $RH$ can hold for non-algebraic links too,
but algebraic links, more generally 
{\em positively iterated torus links\,}, are a major class 
for $RH$ (for small $q$).
\vskip 0.2cm

{\sf Deviations from the classical theory.} Focusing on the 
sector $0\!<\!q\!<\!1$ is a significant deviation. Even for $q$ 
close 
to $0$, {\em irregular\,} zeros generally appear. For instance, 
one pair (always real) of such zeros is conjectured to occur for 
uncolored algebraic links with $2$ components. 

The main change is of course that $q$ is {\em arbitrary\,}  real
in the DAHA approach. It is a counterpart of
$|\mathbb F|$ in
the Weil conjectures, but this is just a parameter
for us. Accordingly, we calculate minimal $\varpi$ such 
that $RH$ holds for all $\om\equal 1/q>\varpi$. Also, we add
colors and one more parameter $a$ due to the flagged Jacobian 
factors (or Hilbert schemes), which is a clear extension of
the Weil conjectures. 

{\em The field with one element\,} is of importance. 
It corresponds to $q\!=\!t$
in the DAHA parameters (which becomes $q=1$ after
the substitution $q\mapsto qt$), and describes the HOMFLY-PT 
polynomials in topology  However, the bound $\varpi$ is generally 
beyond $1$; $RH$ is not expected
to hold for $q=t$ for sufficiently general non-torus knots.
We note that $t=1$ can be also considered as the case of ``field
with one element"; see \cite{ChP1,ChP2}.
\vskip 0.2cm

\subsubsection{\sf On the structure of the paper}  
The motivic
Conjecture \ref{CONJ:HL} is the best we have
to clarify the meaning of the 
substitution $q\mapsto qt$. 
Finding $\pi_i$ and polynomials $S^i$ from
Conjecture \ref{CONJLIM} is an entirely algebraic procedure
and for any Young diagrams (topological connections are expected). 
Conjecture \ref{CONJRH}, the
``qualitative $RH$", 
 for small $q$  is actually a corollary of 
Conjecture \ref{CONJLIM}. 
Conjecture  \ref{CON:om2}, the ``quantitative $RH$", 
gives a bound for $q$ 
in the case of uncolored algebraic knots.  



If the motivic superpolynomials are known and coincide 
with the DAHA superpolynomials (which is conjectured), then
the validity of $RH$ for sufficiently small $q$ 
(sufficiently large $\om$), i.e. Conjectures 
\ref{CONJLIM},\ref{CONJRH}, can be checked by
an entirely algebraic and relatively straightforward
procedure. 
For uncolored algebraic knots, the geometric superpolynomials
were defined in \cite{ChP1} and they do satisfy Conjecture 
\ref{CONJLIM}, which almost directly follows from their
definition in terms of the corresponding flagged
Jacobian factors. 

The motivic superpolynomials at $a\!=\!0$
are quite close 
to the Hasse-Weil zeta-function
for singular curves and Galkin-St\"ohr zeta-functions.
The exact relation is Conjecture \ref{CONJ:HL}; 
Conjecture \ref{CONJ:HLF} is its generalization to any $a$.
Adding colors is more subtle; the {\em geometric  
superpolynomials\,} of algebraic links are known by now 
for (any) columns. The DAHA superpolynomials are 
well developed for any Young diagrams. 
\vskip 0.2cm

\subsection{\bf Motivic approach} Let us discuss some
arithmetic-geometric details. Connecting DAHA superpolynomials
with the numerators of the
Hasse-Weil zeta-functions seems {\em a priori\,} some stretch,
but we think that the following chain of steps provides a 
sufficiently solid link.

\subsubsection{\sf Kapranov's zeta}
The first step is the Kapranov zeta-function of a smooth algebraic
curve $C$ of genus $g$ over a field $k$. It is defined
via the classes of $[C^{[n]}]$ of
the $n$-fold symmetric products of $C$ in
the Grothendieck ring $K_0(V\!ar/k)$
of varieties over $k$. The motivic zeta-function of $C$
from \cite{Kap}
is then a formal series $\ze(C,u)\equal
\sum_{n\ge 0}\, u^{n} [C^{[n]}]$. Here one
can replace $[C^{[n]}]$ by $\mu([C^{[n]}])$ for
any {\em motivic measure\,} $\mu$. If $k=\mathbb F=\mathbb F_q$
and $\mu(X)$ is the number of $\mathbb F$\~points of $X$, then
this is the classical presentation of the
{\em Hasse-Weil zeta-function} for $u=T$. 

Theorem 1.1.9 from \cite{Kap}
establishes the first two Weil conjectures
(the rationality and the functional equation) in the
motivic setting. The justification is close to the Artin's
proof in the case of $\mathbb F$.

One can then extend the definition of motivic zeta
to reduced {\em singular\,}
curves $C$, replacing $C^{[n]}$ by the corresponding Hilbert
schemes of $n$ points on $C$, subschemes of length $n$ to be exact. 
Let us assume now (and later)
that $C$ is a {\em rational\,}
planar projective reduced curve of arithmetic
genus $\de$. Then Conjecture 17 
from \cite{GSh} states that in $K_0(V\!ar/\C)$:
\begin{align}\label{PandTGS}
\sum_{n\ge 0} u^{n+1-\de} [C^{[n]}]=
\sum_{0\le i\le \de} {\n}_{c}{(i)}\bigl(\frac{u}
{(1-u)(1-u[\mathbb A^1])}\bigr)
^{i+1-\de},
\end{align}
where
$\n_C^{(i)}\in \Z_+[\mathbb A^1]$.
To be exact, this is stated
for any reduced curve, not only
rational, the total arithmetic genus $g$
must be then used in the left-hand side instead
of $\de$ and the resulting expression must be divided 
by the same series for
the normalization of $C$ (calculated by Macdonald for
smooth $C$). In the right-hand side, we must set
 $i+1-\de\mapsto i-\de$. This substitution  is due to 
the Macdonald formula for $\mathbb P^1$.

\subsubsection{\sf Nested Hilbert schemes}
When the classes $[X]$ are replaced by their (topological)
Euler numbers $e(X)$, we arrive at
\begin{align}\label{PandT}
 \sum_{n\ge 0} u^{n+1-\de} e(C^{[n]})=
\sum_{0\le i\le \de} n_{C}{(i)}\left(\frac{u}
{(1-u)^2}\right)
^{i+1-\de}.
\end{align}
The rationality here was motivated by
Gopakumar and Vafa (via the BPS invariants) and justified in
\cite{PaT}.
The positivity of $n_C^{(i)}$
was deduced from the approach based on versal deformations:
\cite{FGV} for $i=\de$ and then (for arbitrary $i$)
in \cite{Sh}.
The OS-conjecture \cite{ObS}, (extended and)
proved in \cite{Ma}, is a geometric interpretation 
and an $a$-generalization 
of (\ref{PandT}) for rational planar curves $C$
and their {\em nested\,} Hilbert schemes $C^{[n\le n+m]}$. 
It is actually a local formula and one can switch
from a rational curve $C$ to its {\em germ} $\c$ at
the singular point under consideration.
Then $\c^{[n\le n+m]}=\{I_{n+m}\!\subset\! I_n\mid
\mathfrak{m}I_n\!\subset\! I_{n+m}\}$, where 
$I_n\!\in\! \c^{[n]}$ and $\mathfrak{m}$ is the 
maximal ideal in the (local) ring of $\c$.

\subsubsection{\sf ORS polynomials}
Let $\c$ be an arbitrary plane curve singularity
of arithmetic genus $\de$ (its Serre number);
the Hilbert schemes are defined correspondingly.
Considering the construction above for 
$K_0(V\!ar/\mathbb F)$ and
then applying the motivic integration 
from Example 1.3.2b from
\cite{Kap} is essentially what was suggested in \cite{ORS}
(for nested Hilbert schemes of punctual pairs). 
The reduced ORS polynomial is 
\begin{align}\label{overlinep}
&\mathscr{P}_{\!\hbox{\tiny alg}}=
\bigl(\frac{q_{st}}{a_{st}}\bigr)^{\mu}
\frac{1-q_{st}^2}{1+a_{st}^2 t_{st}}
\,\sum_{n,m\ge 0}q_{st}^{2n}\,a_{st}^{2m}\,
t_{st}^{m^2}\,\mathfrak{w}(\c^{[n\le n+m]}).
\end{align}
Here $\mu$ is the Milnor number ($\mu=2\de$ in the unibranch
case) and $\mathfrak{w}$ is the {\em weight filtration\,}
in the compactly supported cohomology of the corresponding scheme.
See the Overview and Section 4 in \cite{ORS}; 
Proposition 3 there
contains the functional equation. And also
see  Section 9.1 from \cite{GORS}.
The parameter
$t_{st}$ is associated with this filtration, 
we put $q_{st},t_{st},a_{st}$ here
to distinguish these parameters from the DAHA parameters (below).
They are really {\em standard\,} in quite a few
topological-geometric papers; see (\ref{qtarelidaha}) below.
For $t_{st}=1$, $u=q_{st}$ and at the minimal possible degree
of $a_{st}$, the sum in (\ref{overlinep}) essentially reduces to
the right-hand side of (\ref{PandT}).

Conjecture 2 of \cite{ORS} states that
$\mathscr{P}_{\!\hbox{\tiny alg}}$
coincides with the {\em reduced\,} stable 
Khovanov-Rozansky
polynomial of the corresponding link, the
Poincar\`e polynomials of the triply-graded HOMFLY-PT homology.
Accordingly, its unreduced version
$\overline{\mathscr{P}}_{\!\hbox{\tiny alg}}$ 
is related to the unreduced 
stable KhR polynomial. 
The problem with
the identification of the ORS and KhR polynomials is that 
the number of examples is very limited in both theories.
Though see \cite{Mel2} concerning the stable KhR polynomials
for torus knots via Soergel bimodules. 

\subsection{\bf DAHA approach}
The following two families of superpolynomials are
much more explicit and calculatable.
{\em DAHA superpolynomials\,} are the key for us;
their full definition will be provided in this paper.
They were defined in \cite{CJ} for torus knots,
triggered by \cite{AS}, and extended to any iterated torus
links in further papers.

\subsubsection{\sf Geometric superpolynomials}
The connection of the DAHA superpolynomials
to geometry of plane curve singularities
goes through the {\em geometric superpolynomials\,},
a class of superpolynomials introduced
in \cite{ChP1,ChP2}, which generalizes \cite{ChD1},
\cite{GM} (for $a=0$ and torus knots) and Gorsky's approach
from \cite{Gor} (a combinatorial theory
for torus knots and any powers of $a$).
{\em Flagged Jacobian factors\,} were used in
\cite{ChP1} instead of the nested Hilbert schemes in \cite{ORS};
though see Section \ref{SEC:HILB} where {\em restricted\,}
nested Hilbert schemes emerge (the pairs of ideals that become
trivial upon tensoring with the normalization ring). 

The geometric superpolynomials
do not require root systems
at and the {\em $a$\~stabilization\,} of  
the {\em DAHA-Jones-WRT polynomials\,}. They  
are uniformly defined for any root systems, but their 
$a$\~stabilization 
naturally requires classical series $A,B,C,D$. However  see 
\cite{ChE} for some superpolynomials for the Deligne-Gross 
exceptional family \cite{DG}. Presumably they
exist for any root systems and their construction is
geometric, in terms of the corresponding {\em spectral curve\,}. 

The relation between the DAHA superpolynomials and
geometric superpolynomials for algebraic knots 
from \cite{ChP1,ChP2} is confirmed in many
examples and is verifiable theoretically. One 
can connect the standard (monoidal) transformations of 
the plane curve singularities with the corresponding
recurrence relations
in the DAHA theory. 

\subsubsection{\sf ORS polynomials vs. geometric ones}
Our uncolored geometric superpolynomials are parallel to
$\mathscr{P}_{\!\hbox{\tiny alg}}$ from \cite{ORS}.
Upon the division by $(1-t)$ there is a conjectural connection
with those.
See Sections \ref{SEC:FLAGGED} and \ref{SEC:HILB} below.
They can be expressed is terms of the weight filtration 
too via a theorem due to N. Katz; see the end of \cite{ChP1}. 
The DAHA parameter $q$ (or $t^{-1}$ due to the super-duality)
is a counterpart of  cardinality of $\mathbb F$ there.
Accordingly, $q=1$  is the case
of the ``field with one element", which leads to the
DAHA theory at critical center charge for $q=1$.
When $t^{-1}$ is used, the case $t=1$ is the so-called  
"free theory". This is different from the usage of the 
weight filtration in \cite{ORS}, where
$|\mathbb F|$ corresponds to $q_{st}^2=q/t$.  
Generally, the DAHA parameters $a,q,t$
are connected with those in
(\ref{overlinep}) as follows:

\begin{align}\label{qtarelidaha}
&t=q^2_{st},\  q=
(q_{st}t_{st})^2,\  a=a_{st}^2 t_{st},\notag\\
&q^2_{st}=t,\  t_{st}=\sqrt{\frac{q}{t}},
\  a_{st}^2=a\sqrt{\frac{t}{q}}.
\end{align}
In the ORS polynomials, the motivic measure
and finite fields vanish at $t_{st}=1$, i.e. at $q=t$ in
the DAHA parameters, which is quite different 
from $q=1$ in geometric superpolynomials. 
Also, the normalization of singularities is heavily used for 
the {\em flagged Jacobian factors\,} and
our construction does not 
require ideals and Hilbert schemes, though the Galkin-St\"ohr 
zeta-functions are directly related to the latter.

Let us mention here that 
colors and arbitrary root systems 
(available in the DAHA approach) are a challenge
from the geometric-motivic perspective. However the
{\em OS conjecture\,} (the
case $t=q$ in DAHA parameters) was established and proved
in \cite{Ma} for {\em any\,} colors (Young diagrams) and algebraic
links using {\em non-reduced\,} singularities.

\comment{
In spite of such discrepancies and not too rigorous connections
with the Weil conjectures, we think
there are sufficiently strong {\em a priori\,} reasons to
expect $RH$ for the DAHA
superpolynomials. They must be recalculated via $q/t$ instead
of $q$; such new $q$ essentially coinciding with
$t_{st}$. At least, this argument can work for $a\!=\!0$ and
without colors. This is also supported by the fact that
the first two Weil conjecture hold in the DAHA setting,
which was discussed above.}

\subsubsection{\sf Galkin-St\"ohr zeta}\label{Gal-St}
We will consider here only the
case $a=0$. The motivic measure will be 
the count of points over a finite field $\mathbb{F}$
of cardinality $q$. 
With some simplifications, the zeta-function from
 \cite{Gal,Sto}
is defined as the sum $Z_\r(t)$ of $t^{dim_\mathbb{F}
(\r/\mathfrak{a})}$ over all ideals $\mathfrak{a}\subset \r$
for the local ring $\r$ of a curve singularity. This
ring is assumed
Gorenstein to ensure the functional equation for $(1-t)Z_\r(t)$,
which is for the substitution $t\mapsto 1/(tq)$. The function
 $Z_\r(t)$ is a (local) version of the Weil 
zeta for singularities. It is quite natural from number 
theoretical viewpoint, but the corresponding $RH$ generally
fails. However if we treat $q$ as a variable
(which is $q/t$ in the DAHA parameters),  then
$q\le 1/2$ is presumably sufficient for $RH$
in the unibranch case at $a\!=\!0$. 

The connection
with our geometric superpolynomials can be stated as
some ``combinatorial" identity, which seems not 
straightforward to check; for instance, 
it generally holds only for {\em planar\,} singularities
(not arbitrary Gorenstein rings). Also, the positivity of the 
coefficients of  $(1\!-\!t)Z_\r(t)$, which follows from this
connection, generally requires  plane curve singularities.
The formula from \cite{Sto} has many positive and negative terms  
canceling each other in a non-trivial way. We provide examples
of  {\em non-plane\,} curve singularities where 
the connection with our formula can be fixed, but
some non-trivial adjustments are needed for this. 
Importantly, the functional equation for  $(1\!-\!t)Z_\r(t)$
is not difficult to established (see \cite{Sto}) and we
do not see any {\em direct\,} proof
of super-duality for our geometric
(motivic) superpolynomials, without the passage to ideals.

\subsubsection{\sf Modular periods} 
The construction of DAHA superpolynomials is actually
parallel to that for  modular periods, a starting point for
$p$\~adic measures (Mazur, Manin, Katz, eigenvarieties),
and for Manin's zeta-polynomials \cite{OnRS}. Namely,  
the DAHA-Jones polynomials of torus knots $T(\rr,\ss)$ 
colored by Young diagrams $\la$ are essentially obtained by 
applying $\ga_{\rr,\ss}=$ {\tiny
$\begin{pmatrix} \rr & * \\ \ss & *\end{pmatrix}$}$\in PSL(2,\Z)$ 
to the Macdonald polynomials $P_\la$ followed by taking the 
DAHA coinvariant. The superpolynomials are
due to the $a$\~stabilization for the root systems $A_n$; 
this is important, since the super-duality holds only upon
the $a$\~stabilization.

This procedure is analogous
to taking the integral of a cusp form $\Phi(z)$ for
$z\in\mathbb H$ multiplied
by $z^{k}$ for certain integers $k$ over the paths 
$\ga[0,\imath\infty]$. Here $z^{k}$ can be seen 
(with some stretch) as counterparts of $P_\la$, the integration 
$\int \{\cdot\}\Phi(y) dy$ plays the role of the 
{\em coinvariant}. 

The latter obviously has nothing to
do with modular forms. However, it is the simplest
level-one coinvariant among all DAHA coinvariants 
of arbitrary levels $\ell>0$ 
from \cite{ChM}. They are in one-to-one correspondence
with elliptic functions of level $\ell$ (Looijenga functions 
for any root systems). Using them makes
these constructions closer to each other. 
Vice versa, a challenge is to find modular counterparts
of the DAHA-Jones polynomials and superpolynomials 
for iterated non-torus algebraic knots. 
It is not impossible that they can be related to
Manin's iterated Shimura integrals \cite{Man}.

This is connected with the following (heuristic) interpretation
of the Dirichlet $L$\~functions of conductor $\rr$ via
the {\em families\,} $T(\rr,\ast)$. The sums of DAHA 
superpolynomials over the knots in such families are
supposed to be considered. In a more conceptual way, one 
applies $\sum_{m=1}^\infty \chi(m)${\tiny
$\begin{pmatrix} 1 & m \\ 0 & 1\end{pmatrix}$} inside
the coinvariant for the corresponding Dirichlet
character $\chi$.  The $q$\~analogs of zeta 
and Dirichlet functions from \cite{ChZ} are of this kind,
but for a different action of $PSL(2,\Z)$. They are
the integrals of
$\sum_{m=1}^\infty \chi(m) q^{mx^2/2}$ for the Gaussian 
$q^{x^2/2}$ with respect to the Macdonald measure
for $A_1$. See also Section 3 in \cite{ChD2}.

\comment{
\subsection{\bf Some perspectives}\label{PERSPEC}
To put our $RH$\~ conjectures into perspective, we will
conclude this introduction with the following
remarks.

$(i)$ The construction of DAHA superpolynomials is actually
parallel to that for  modular periods, a starting point for
$p$\~adic measures (Mazur, Manin, Katz, eigenvarieties)
and Manin's zeta-polynomials \cite{OnRS}. Namely, torus knots
$T(\rr,\ss)$ play a role similar
to that of integration paths
$[0,\frac{\rr}{\ss}]$ and  the integration with a cusp-form is
some analog of taking the DAHA coinvariant. The latter 
has no relation to any modular forms. However, the
coinvariant we use is the simplest in the family
of {\em general DAHA coinvariants\,}, which are one-to-one
with elliptic functions \cite{ChM}. Using all of them
can make modular periods one step closer to DAHA theory.

$(ii)$ In DAHA theory, we deal with iterated torus links (and
their topological classification) entirely 
via $PSL(2,\Z)$, which alone is remarkable. The general DAHA
coinvariants can be actually obtained from the 
standard one via conjugations by 
the elements of $PSL(2,\Z)$ (upper or lower triangular 
matrices  are sufficient). This leads to the 
following (heuristic) connection. The Dirichlet
$L$\~functions of conductor $\rr$ can be associated with the
{\em families\,} $T(\rr,\ast)$ of torus knots and
the corresponding sums of the invariants over such
families.  
A challenge here is to  interpret the {\em iterated\,} non-torus 
algebraic knots in some similar way. 
It is not impossible that they can be some counterparts
of Manin's iterated Shimura integrals \cite{Man}.

$(iii)$ Section 3 in \cite{ChD2}
provides some link to the $q$\~zeta function and
$q$\~Dirichlet functions
from \cite{ChZ}. 
Integration with the powers of the Gaussian, which is actually
the $th$\~function of the corresponding root system,
cab be used as 
instead of applying the DAHA coinvariant. This 
provides {\em analytic\,} version  
of DAHA-Jones theory. The invariants (superpolynomials,
DAHA-Jones polynomials) 
are replaced by the integrals with powers of $\ga(\th)$
for $\ga\in PSL(2,\Z)$. Accordingly, their sums over the families
from $(ii)$ generalize the construction from \cite{ChZ}; the
$q$\~zeta there corresponds in this sense to some
iterations of Hopf links.
}

\setcounter{equation}{0}
\section{\sc Double Hecke algebras}
\subsection{\bf Affine root systems}
Let us adjust the standard DAHA definitions to
the case of the root systems $A_n$, which is
$R\!=\!\{\al\!=\!\vep_i\!-\!\vep_j, i\!\neq\! j\}$
for the basis $\{\vep_i, 1\le i\le n\!+\!1\}\in \R^{n+1}$,
orthonormal with respect to the usual euclidean form
$(\cdot,\cdot)$.
The Weyl group is $W=\S_{n+1}$; it is
generated by the reflections (transpositions)
$s_\al$ for the set of positive  roots
$R_+=\{\vep_i-\vep_j, i<j\}$; $R_-=-R_+$.
The simple roots are $\al_i=\vep_i\!-\!\vep_{i+1}.$
The weight lattice is
$P=\oplus^n_{i=1}\Z \om_i$,
where $\{\om_i\}$ are fundamental weights:
$ (\om_i,\al_j)=\de_{ij}$. Explicitly,
\begin{align}
&\om_i=\vep_1+\ldots+\vep_i-\frac{i}{n+1}(\vep_1+\ldots+\vep_{n+1})
\for i=1,\ldots,n,\\
&\rho=\om_1\!+\ldots+\!\om_n=\frac{1}{2}\bigl(
(n\!-\!1)\vep_1+(n\!-\!3)\vep_2+\ldots+(1\!-\!n)
\vep_n\bigr).\notag
\end{align}

The root lattice is denoted by $Q=\oplus_{i=1}^n \Z\al_i$.
Replacing $\Z$ by $\Z_{\pm}=\{m\in\Z, \pm m\ge 0\}$, we obtain
$P_\pm,Q_\pm$. See  e.g., \cite{Bo} or \cite{C101}.

The vectors $\ \tal=[\al,j]\in \R^{n+2}$
for $\al \in R, j \in \Z $ form the
{\em affine root system\,}
$\tR \supset R$, where $\al\in R$ are identified with
$[\al,0]$.
We add $\al_0 \equal [-\th,1]$ to the simple
roots for the {\em maximal root\,} $\th=\vep_{1}-\vep_{n+1}$.
The corresponding set
$\tR_+$ of positive roots is 
$R_+\cup \{[\al,j],\ \al\in R, \ j > 0\}$.


\subsubsection{\sf Affine Weyl group}
Given $\tal=[\al,j]\in \tR,  \ b \in P$, let
\begin{align}
&s_{\tal}(\tz)\ =\  \tz-(z,\al^\vee)\tal,\
\ b'(\tz)\ =\ [z,\ze-(z,b)]
\label{ondon}
\end{align}
for $\tz=[z,\ze]\in \R^{n+2}$.
The
{\em affine Weyl group\,}
$\tW=\lan s_{\tal}, \tal\in \tR_+\ran$
is the semidirect product $W\lsmash Q$ of
its subgroups $W=$ $\lan s_\al,
\al \in R_+\ran$ and $Q$, where $\al$ is identified with
\begin{align*}
& s_{\al}s_{[\al,\,1]}=\
s_{[-\al,\,1]}s_{\al}\for
\al\in R.
\end{align*}

The {\em extended Weyl group\,} $ \hW$ is $W\lsmash P$, where
the action is
\begin{align}
&(wb)([z,\ze])\ =\ [w(z),\ze-(z,b)] \for w\in W, b\in P.
\label{ondthr}
\end{align}
It is isomorphic to $\tW\rsmash \Pi$ for $\Pi\equal P/Q$.
The latter group consists of $\pi_0=$id\, and the
images $\pi_r$ of $\om_r$ in $P/Q$.
Note that
$\pi_r^{-1}$ is $\pi_{r^\iota}$,  where $\iota$ is
the standard involution of the {\em non-affine\,}
Dynkin diagram of $R$,
induced by $\al_i\mapsto \al_{n+1-i}$.
Generally, we set
$\iota(b)=-w_0(b)=b^{\,\iota}\,$, where $w_0$ is the
longest element
in $W$ sending $ \{1,2,\ldots,n+1\}$
to $\{n+1,\ldots,2,1\}$.

The group $\Pi$
is naturally identified with the subgroup of $\hW$ of the
elements of the length zero; the {\em length\,} is defined as
follows:
\begin{align*}
&l(\hw)=|\La(\hw)| \for \La(\hw)\equal\tR_+\cap \hw^{-1}(-\tR_+).
\end{align*}
One has $\om_r=\pi_r u_r$ for $1\le r\le n$,
where $u_r$ is the
element $u\in W$ of {\em minimal\,} length such that
$u(\om_r)\in P_-$,
equivalently, $w=w_0u$ is of {\em maximal\,} length such that
$w(\om_r)\in P_+$. The elements $u_r$ are very explicit.
Let $w^r_0$ be the longest element
in the subgroup $W_0^{r}\subset W$ of the elements
preserving $\om_r$; this subgroup is generated by simple
reflections. One has:
\begin{align}\label{ururstar}
u_r = w_0w^r_0 \hbox{\,\, and\,\, } (u_r)^{-1}=w^r_0 w_0=
u_{r^\iota} \for 1\le r\le n.
\end{align}
\smallskip

Setting $\hw = \pi_r\tw \in \hW$ for $\pi_r\in \Pi,\, \tw\in \tW,$
\,$l(\hw)$ coincides with the length of any reduced decomposition
of $\tw$ in terms of the simple reflections
$s_i,\, 0\le i\le n.$ Thus, indeed, $\Pi$ is a subgroup of
$\hW$ of the elements of length $0$.

\subsection{\bf Definition of DAHA}
We follow \cite{CJJ,CJ,C101}.
Let $\mm\!=\!n\!+\!1$; generally it is the least natural number
such that  $(P,P)=(1/\mm)\Z.$
The {\em double affine Hecke algebra, DAHA\,}, of
type $A_n$ depends
on the parameters
$q, t$ and will be defined
over the ring
$\Z_{q,t}\equal\Z[q^{\pm 1/\mm},t^{\pm 1/2}]$
formed by
polynomials in terms of $q^{\pm 1/\mm}$ and
$t^{\pm1/2}.$ Note that the coefficients of the
Macdonald polynomials will belong to
$
\Q(q,t).
$

It is convenient to use the following notation:
\begin{align*}
&t=q^k,\ \rho_k\equal \frac{k}{2}\,\sum_{\al>0}\al=
k\sum_{i=1}^n  \om_i.
\end{align*}
We set $m_{i\, i\!+\!1}=3=m_{n 0}$ for $0\le i\le n$
and $m_{ij}=2$ otherwise, generally, $\{2,3,4,6 \}$
when the number of links  between $\al_i,\al_j$ in
the affine Dynkin diagram is $\{0,1,2,3\}$.

For pairwise commutative $X_1,\ldots,X_n,$
\begin{align}
& X_{\tb}\ \equal\ \prod_{i=1}^nX_i^{l_i} q^{ j}
\iif \tb=[b,j],\ \hw(X_{\tb})\ =\ X_{\hw(\tb)},
\label{Xdex}\\
&\hbox{where\ } b=\sum_{i=1}^n l_i \om_i\in P,\ j \in
\frac{1}{\mm}\Z,\ \hw\in \hW.
\notag \end{align}
For instance, $X_0\equal X_{\al_0}=qX_\th^{-1}$.
\medskip

\begin{definition}
The double affine Hecke algebra $\HH\ $
is generated over $\Z_{q,t}$ by
the elements $\{ T_i,\ 0\le i\le n\}$,
pairwise commutative $\{X_b, \ b\in P\}$ satisfying
(\ref{Xdex})
and the group $\Pi,$ where the following relations are imposed:

(o)\ \  $ (T_i-t^{1/2})(T_i+t^{-1/2})\ =\
0,\ 0\ \le\ i\ \le\ n$;

(i)\ \ \ $ T_iT_jT_i...\ =\ T_jT_iT_j...,\ m_{ij}$
factors on each side;

(ii)\ \   $ \pi_rT_i\pi_r^{-1}\ =\ T_j \iif
\pi_r(\al_i)=\al_j$;

(iii)\  $T_iX_b \ =\ X_b X_{\al_i}^{-1} T_i^{-1} \iif\,
(b,\al_i)=1,\
0 \le i\le  n$;

(iv)\ $T_iX_b\ =\ X_b T_i\, $ when $\, (b,\al_i)=0
\for 0 \le i\le  n$;

(v)\ \ $\pi_rX_b \pi_r^{-1}= X_{\pi_r(b)}\, =\,
X_{ u^{-1}_r(b)}
 q^{(\om_{\iota(r)},b)},\ \, 1\le r\le n$.
\label{double}\end{definition}

Given $\tw \in \tW, 1\le r\le n,\ $ the product
\begin{align}
&T_{\pi_r\tw}\equal \pi_r T_{i_l}\cdots T_{i_1},\where
\tw=s_{i_l}\cdots s_{i_1} \for l=l(\tw),
\label{Twx}
\end{align}
does not depend on the choice of the reduced decomposition
of $\tw$.
Moreover,
\begin{align}
&T_{\hv}T_{\hw}\ =\ T_{\hv\hw}\  \hbox{ whenever\,}\
 l(\hv\hw)=l(\hv)+l(\hw) \for
\hv,\hw \in \hW. \label{TTx}
\end{align}
In particular, we arrive at the pairwise
commutative elements
\begin{align}
& Y_{b}\equal
\prod_{i=1}^nY_i^{l_i} \iif
b=\sum_{i=1}^n l_i\om_i\in P,\
Y_i\equal T_{\om_i},b\in P.
\label{Ybx}
\end{align}
When acting in the polynomial representation
$\v$ (see below), they are called {\em difference
Dunkl operators}.

\subsection{\bf The automorphisms}\label{sect:Aut}
The following maps can be (uniquely) extended to
automorphisms of $\HH\,$, where
$q^{1/(2\mm)}$ must be added to $\Z_{q,t}$
(see \cite{C101}, (3.2.10)--(3.2.15))\,:
\begin{align}\label{tauplus}
& \tau_+:\  X_b \mapsto X_b, \ T_i\mapsto T_i\, (i>0),\
\ Y_r \mapsto X_rY_r q^{-\frac{(\om_r,\om_r)}{2}}\,,
\\
& \tau_+:\ T_0\mapsto  q^{-1}\,X_\th T_0^{-1},\
\pi_r \mapsto q^{-\frac{(\om_r,\om_r)}{2}}X_r\pi_r\
(1\le r\le n),\notag\\
& \label{taumin}
\tau_-:\ Y_b \mapsto \,Y_b, \ T_i\mapsto T_i\, (i\ge 0),\
\ X_r \mapsto Y_r X_r q^\frac{(\om_r,\om_r)}{ 2},\\
&\tau_-(X_{\th})=
q T_0 X_\th^{-1} T_{s_{\th}}^{-1};\ \
\si\equal \tau_+\tau_-^{-1}\tau_+\, =\,
\tau_-^{-1}\tau_+\tau_-^{-1},\notag\\
&\si(X_b)=Y_b^{-1},\   \si(Y_b)=
T_{w_0}^{-1}X_{b^{\,\iota}}^{-1}T_{w_0},\ \si(T_i)=T_i (i>0).
\label{taux}
\end{align}
These automorphisms fix $\ t,\ q$
and their fractional powers, as well as the
following {\em anti-involution\,}:
\begin{align}
&\vph:\
X_b\mapsto Y_b^{-1},\, Y_b\mapsto X_b^{-1},\
T_i\mapsto T_i\ (1\le i\le n).\label{starphi}
\end{align}

The following anti-involution results directly from
the group nature of the DAHA relations:
\begin{align}\label{star-conj}
H^\star= H^{-1} \for H\in \{T_{\hw},X_b, Y_b, q, t\}.
\end{align}
To be exact, it is naturally extended to the fractional
powers of $q,t$:
$$
\star:\ t^{\frac{1}{2}} \mapsto t^{-\frac{1}{2}},\
q^{\frac{1}{2\mm}}\mapsto  q^{-\frac{1}{2\mm}}.
$$
This anti-involution serves the inner product in the theory
of the DAHA polynomial representation.

Let us list the matrices corresponding to the automorphisms and
anti-automorphisms above upon the natural projection
onto $SL_2(\Z)$, corresponding to
$\,t^{\frac{1}{2}}=1=q^{\frac{1}{2\mm}}$.
The matrix {\tiny
$\begin{pmatrix} \al & \be \\ \ga & \de\\ \end{pmatrix}$}
will then represent the map $X_b\mapsto X_b^\al Y_b^\ga,
Y_b\mapsto X_b^\be Y_b^\de$ for $b\in P$. One has:
\smallskip

\centerline{
$\!\!\!\tau_+\rightsquigarrow$
{\tiny
$\begin{pmatrix}1 & 1 \\0 & 1 \\ \end{pmatrix}$},\
$\tau_-\rightsquigarrow$
{\tiny
$\begin{pmatrix}1 & 0 \\1 & 1 \\ \end{pmatrix}$},\
$\si\rightsquigarrow$
{\tiny
$\begin{pmatrix}0 & 1 \\-1 & 0 \\ \end{pmatrix}$},\
$\vph\rightsquigarrow$
{\tiny
$\begin{pmatrix}0 & -1 \\-1 & 0 \\ \end{pmatrix}$}.\
}
\smallskip

The {\em projective\,} $PSL_2(\Z)$ (due to Steinberg)
is the group generated by $\tau_{\pm}$ subject to the
relation $\tau_+\tau_-^{-1}\tau_+=
\tau_-^{-1}\tau_+\tau_-^{-1}.$ The notation will
be $PSL_{\,2}^{\wedge}(\Z)$; it is isomorphic to the braid
group $B_3$.

\subsubsection{\sf The coinvariant}
The projective $PSL_2(\Z)$ and the {\em coinvariant\,},
to be defined now, are the main ingredients of our approach.

Any $H\in \HH$ can be uniquely represented in the form
$$
H=\sum_{a,w,b} c_{a,w,b}\, X_a T_{w} Y_b \for w\in W,
a,b\in P
$$
(the DAHA-PBW theorem, see \cite{C101}). Using this
presentation, the
{\em coinvariant\,} is a functional $\HH\to \C\,$
defined as follows:
\begin{align}\label{evfunct}
\{\,\}_{ev}:\ X_a \ \mapsto\  q^{-(\rho_k,a)},\
Y_b \ \mapsto\  q^{(\rho_k,b)},\
T_i \ \mapsto\  t^{1/2}.
\end{align}
The main symmetry of the coinvariant is
\begin{align}\label{evsym}
&\{\,\vph(H)\,\}_{ev}\,=\,\{\,H\,\}_{ev} \hbox{\, for\, }
H\in \HH.
\end{align}
Also,
$\,\{\,\iota(H)\,\}_{ev}=\{\,H\,\}_{ev}$, where we
extend $\iota$ to $\HH$ as follows:
\begin{align}\label{iotaXY}
\iota(X_b)\!=\!X_{\iota(b)}, \
\iota(Y_b)\!=\!Y_{\iota(b)}, \
T_i^\iota\!=\!T_{\iota(i)},\ 1\le i\le n.
\end{align}

The following interpretation of the coinvariant is
important.
For any $H\in \HH$,
one has: $\{H T_w Yb\}_{ev}=\{H\}_{ev}\, \chi (T_w Y_b)$,
where $\chi$ is the standard character (one-dimensional
representation) of the affine Hecke algebra
$\h_Y$, generated by $T_w, Y_b$ for
$w\in W, b\in P$;\, $\chi$\,  sends\,
$Y_b\mapsto q^{(\rho_k,b)}$ and $ T_i\mapsto t^{1/2}$.
Therefore $\{\ldots\}_{ev}$ acts via the projection
$H\mapsto H\!\Downarrow \equal H(1)$ of $\HH\,$
onto the {\em polynomial representation \,}$\v$, which is
the $\HH$\~module induced from $\chi$;
see \cite{C101,CJ,CJJ} and the next section.

\subsection{\bf Macdonald polynomials}
\subsubsection{\sf Polynomial representation}\label{SEC:POLR}
It is isomorphic to $\Z_{q,t}[X_b]$
as a vector space with the action of $T_i(0\le i\le n)$
given by
the {\em Demazure-Lusztig operators\,}:
\begin{align}
&T_i\  = \  t^{1/2} s_i\ +\
(t^{1/2}-t^{-1/2})(X_{\al_i}-1)^{-1}(s_i-1),
\ 0\le i\le n.
\label{Demazx}
\end{align}
The elements $X_b$ become the multiplication operators
and  $\pi_r (1\le r\le n)$ act via the general formula
$\hw(X_b)=X_{\hw(b)}$ for $\hw\in \hW$. Note that $\tau_-$
naturally acts in the polynomial representation.
See formula (1.37) from \cite{CJJ}, which is based on the
identity
\begin{align}\label{taumin-poly}
\tau_-(H\!\Downarrow)\,=\,\tau_-(H)\!\Downarrow\ =\,
\Bigl(\tau_-\bigl(H\bigr)\Bigr)(1) \for H\in \HH.
\end{align}
\smallskip

{\sf Symmetric Macdonald polynomials.}
The standard notation is  $P_b(X)$ for $b\in P_+$;
see \cite{M2,C101} (they are due to Kadell
for the classical root systems and due to
Rogers for $A_1$). The usual definition is
as follows. Let $c_+$ be such that $c_+\in W(c)\cap P_+$
(it is unique);
recall that $Q_+=\oplus_{i=1}^n \Z_+\al_i$.
 For $b\in P_+$, the following are the defining relations:
\begin{align*}
&P_b\! -\!\!\!\!\sum_{a\in W(b)}\!\!\! X_{a}
\in\, \oplus^{c_+\neq b}_{c_+\in b-Q_+}\,\Q(q,t) X_c
\hbox{\, and\, }
\lan \,P_b \,X_{c^{\,\iota}}\,\mu(X;q,t)\,\ran\!=\!0 \for
\\
&\hbox{all $c$ in $\oplus$ above;\ \ }
\mu(X;q,t)\!\equal\!\prod_{\al \in R_+}
\prod_{j=0}^\infty \frac{(1\!-\!X_\al q^{j})
(1\!-\!X_\al^{-1}q^{j+1})
}{
(1\!-\!X_\al t q^{j})
(1\!-\!X_\al^{-1}t^{}q^{j+1})}\,.
\end{align*}
Here and further $\lan f\ran$ is the {\em constant term\,}
of a Laurent series or polynomial $f(X)$;
$\mu$ is considered
a Laurent series of $X_b$ with
the coefficients expanded in terms of
positive powers of $q$. The coefficients of
$P_b$ belong to the field $\Q(q,t)$.
One has (see (3.3.23) from \cite{C101}):
\begin{align}\label{macdsym}
&P_b(X^{-1})\,=\,P_{b^{\iota}}(X)\,=\,
P_b^\star(X),\ \, P_{b}(q^{-\rho_k})=
P_{b}(q^{\rho_k})\\
\label{macdeval}
=\,&(P_{b}(q^{-\rho_k}))^\star=
q^{-(\rho_k,b)}
\prod_{\al>0}\prod_{j=0}^{(\al,b)-1}
\Bigl(
\frac{
1- q^{j}\,t\, X_\al(q^{\rho_k})
 }{
1- q^{j}X_\al(q^{\rho_k})
}
\Bigr).
\end{align}
Recall that $\iota(b)=b^{\,\iota}=-w_0(b)$ for $b\in P$.

DAHA provides an important alternative (operator)
approach to the $P$\~polynomials; namely, they satisfy
the (defining) relations
\begin{align}\label{macdopers}
L_f(P_b)=f(q^{-\rho_k-b})P_b,\  L_f\equal f(X_a\mapsto Y_a)
\end{align}
for any symmetric ($W$\~invariant) polynomial
$f\in \C[X_a,a\in P]^W$. Here $b\in P_+$ and the coefficient
of $X_b$ in $P_b$ is assumed $1$.
\smallskip

{\sf Spherical normalization.}
We call $P_b^\circ\equal P_b/P_b(q^{-\rho_k})$
{\em spherical Macdonald polynomials\,} for $b\in P_+$. One
has (the evaluation theorem):
\begin{align}
&P_{b}(q^{-\rho_k})=
q^{-(\rho_k,b)}
\prod_{\al>0}^{(\al,b)>0}\,\prod_{j=0}^{(\al\!,b)-1}
\Bigl(
\frac{
1- q^{j}t X_\al(q^{\rho_k})
 }{
1- q^{j}X_\al(q^{\rho_k})
}
\Bigr).\label{Pbrho'}
\end{align}

\subsection{\bf \texorpdfstring{{\mathversion{bold}$J$}}
{J}-polynomials}
They are necessary below for managing
algebraic links (spherical polynomials $P^\circ_b$
are sufficient for knots) and are important for the
justification of the {\em super-duality\,}.

For $\,b=\sum_{i=1}^n b_i \om_i\in P_+$, i.e.
for a {\em dominant\,} weight with
$\,b_i\ge 0$ for all $\,i$, the corresponding
{\em Young diagram\,} is  as
follows:
\begin{align}\label{omviavep}
&\la\!=\!\la(b)\!=\!\{\la_1\!=\!b_1\!+\!\ldots\!+\!b_n,
\la_2\!=\!b_2\!+\!\ldots\!+\!b_n,\ldots, \la_n\!=\!b_n\},\\
&b=\sum_{i=1}^n \la_i \vep_i-
\frac{|\la|}{n+1}\,\bigl(\vep_1+\ldots+\vep_{n+1}\bigr)
\for |\la|\equal\sum_{i=1}^n \la_i.\notag
\end{align}
One has: $(b,\vep_{i}\!-\!\vep_j)=b_i\!+\!\ldots\!+\!b_{j\!-\!1}=
\la_i\!-\!\la_j$,
$(b,\rho)\!=\!(|\la|\!-\!\la_1)/2$. Also,
$b^2\equal(b,b)=\sum_{i=1}^n \la_i^2-|\la|^2/(n\!+\!1)$.


Let us calculate the  set of all $[\al,j]$ in the
product from (\ref{Pbrho'}); it is
\begin{align*}
&\{\,[\al,j],\  \al=\vep_l\!-\!\vep_{m}\in R_+,\, j>0\, \mid\,
b_l+\ldots+b_{m-1}>j> 0\}.
\end{align*}

\subsubsection{\sf Their definition}\label{SEC:JPOL}
The {\em $J$\~polynomials\,}
are as follows:
\begin{align}\label{P-arms-legs}
J_\la\equal h_\la P_b \for \la=\la(b),\,
h_\la=\prod_{\Box\in\la}
(1-q^{arm(\Box)}t^{leg(\Box)+1});
\end{align}
they are $q,t$\~integral.
\smallskip

Here $arm(\Box)$ is the {\em arm number\,}, which is
the number of boxes in the same row as $\Box$
strictly after it; $leg(\Box)$ is the {\em leg number\,},
which is the number of boxes in the column of $\Box$
strictly below it. This is for the standard presentation
of $\la$\,:
$\la_1\!\ge \la_2\!\ge\!\ldots\!\ge\!\la_{n-1}\!\ge \la_n\,$
are the numbers of
boxes in the corresponding rows and the $i${\tiny th} row
is above the  $(i+1)${\tiny th}.
\smallskip

Equivalently:
\begin{align}\label{j-polynom}
&J_\la= t^{-(\rho,b)}
\prod_{p=1}^n \prod_{j=0}^{\la_{p^*\,}-1}
\Bigl( 1\!-\!q^{j} t^{\,p+1}\Bigr) P_b^\circ, \ \,
p^*= n\!-\!p\!+\!1,\, b\in P_+.\
\end{align}
See, for instance, Theorem 2.1 from \cite{GoN}.
Note that the {\em arms and\, legs\,} do not
appear in the latter presentation (in terms of
$P_b^\circ$). In this approach, 
counterparts of $J$\~polynomials exist for any 
root systems, but there are
some deviations. See \cite{ChD2}, Section 2.6.

\subsubsection{\sf Stabilization}
The following formula is important:
\begin{align}\label{stabevalx}
J_\la(t^{-\rho})\!=\!(a^2)^{-\frac{\la_1}{4}}\,
t^{\frac{|\la|}{2}}\,
&\prod_{p=1}^n \prod_{j=0}^{\la_{p\,}-1}
\Bigl( 1\!+\!q^{j}\, a\, t^{-p+1}\Bigr) \hbox{\,\, for \,\,}
 a\!=\!-t^{n+1}.
\end{align}

More generally, we have the following
$a$\~stabilization; see \cite{ChD2}. We note that the
geometric/motivic superpolynomials are defined without
any  $A_n$ (and further stabilization), but so far they
are known only for columns ($b=m\om_1$).   

\begin{proposition}\label{PROP:stab-values}
Given two Young diagrams $\la$ and $\mu$,\,
the values $P_\la(q^{\mu+\rho_k})$ are $a$\~stable, which means 
that there is a universal expression in terms
of $\,q,t,a\,$ such that its value at $a=-t^{n+1}$ coincides with
 $P_b(q^{c+\rho_k})$ for $\la=\la(b),\mu=\la(c)$,
$b,c\in P_+$ for $A_n$ with $n$ no smaller than the
number of rows in $\la$ and in $\mu$. Up to powers of $a^{1/2}$
and $t^{1/2}$,
they are rational function in  $q,t,a$. Also,
$P_\la^\circ(q^{\mu+\rho_k})$,
$\lan P_\la,P_\la\ran$
and $\lan P^\circ_\la,P^\circ_\la\ran$ are $a$\~stable
(in the same sense). \sq
\end{proposition}

\setcounter{equation}{0}
\section{\sc DAHA-Jones theory}
\subsection{\bf Iterated torus knots}\label{sec:ITER-KNOTS}
We will first introduce the data necessary in the
construction of DAHA-Jones polynomials
and DAHA superpolynomials of algebraic 
iterated torus {\em knots\,}.

\subsubsection{\sf Newton pairs}
The (algebraic) torus knots $T(\rr,\ss)$ are defined for any
positive integers
assuming that \,gcd$(\rr,\ss)=1$.
One has the symmetry $\,T(\rr,\ss)=T(\ss,\rr)$,
where we use ``$=$" for the ambient isotopy equivalence.
Also $\,T(\rr,\ss)=\unknot\,$\, if $\rr=1$ or $\ss= 1$
for the {\em unknot\,} $\unknot\,$. Here and below
see e.g., \cite{RJ,EN,ChD1} for details and/or Knot Atlas
for the corresponding invariants.

Following \cite{ChD1}, the {\em \tax-presentation\,}
of an
{\em iterated torus knots\,} (defined below) will be
$\t(\vec{\rr},\vec{\ss})$ for any two sequences of
positive integers:
\begin{align}\label{iterrss}
\vec{\rr}=\{\rr_1,\ldots \rr_\ell\}, \
\vec{\ss}=\{\ss_1,\ldots \ss_\ell\} \hbox{\, such that\,
gcd}(\rr_i,\ss_i)=1;
\end{align}
$\ell$ will be called the {\em length\,} of $\,\vec\rr,\vec\ss$.
The pairs $[\rr_i,\ss_i]$ can be interpreted as
{\em characteristic\,} or {\em Newton pairs\,}
in terms of plane curve singularities. 
The necessary and sufficient
condition for being algebraic is $\rr_i,\ss_i>0$,
which will be imposed in this paper unless stated otherwise.

\subsubsection{\sf Cabling parameters}
The above parameters are the ones needed
in the DAHA approach. However they are not optimal for 
establishing
the symmetries of our polynomials and the justification
that our construction depends only on the
corresponding knot/link. We need
the {\em cable presentation\,} for this.
It requires one more sequence of integers:
\begin{align}\label{Newtonpair}
\aa_1=\ss_1,\,\aa_{i}=\aa_{i-1}\rr_{i-1}\rr_{i}+\ss_{i}\,\
(i=2,\ldots,m).
\end{align}
See  e.g., \cite{EN}. In terms of
the {\em cabling\,} discussed below, the corresponding knots
are as follows. First, $T(\rr,\ss)=C\!ab(\ss,\rr)(\unknot)$
(note that we transpose $\rr,\ss$ here), and then we set:
\begin{align}\label{Knotsiter}
\t(\vec{\rr},\!\vec{\ss})\rightsquigarrow
C\!ab(\vec{\aa},\!\vec{\rr})(\unknot)=
\bigl(C\!ab(\aa_\ell,\rr_\ell)\cdots C\!ab(\aa_2,\rr_2)\bigr)
\bigl(T(\rr_1,\!\ss_1)\bigr).
\end{align}

Knots and links will be considered up to
{\em ambient isotopy\,}; we use \,``$=$" for it.
The {\em cabling\,} $C\!ab(\aa,\bb)(K)$ of any oriented
knot $K$ in (oriented) $\S^3$ is defined as follows;
see e.g., \cite{Mo,EN} and references therein.
We consider a small $2$\~dimensional torus
around $K$ and put there the torus knot $T(\aa,\bb)$
in the direction of $K$,
which is $C\!ab(\aa,\bb)(K)$ (up to ambient isotopy);
we set $C\!ab(\vec\aa,\vec\rr)\equal
C\!ab(\vec\aa,\vec\rr)(\unknot)$.

This procedure depends on
the order of $\aa,\bb$ and orientation of
$K$. We choose them in the standard way:
the parameter $\aa$ gives the number of
turns around $K$. This construction also depends
on the {\em framing\,} of the cable knots; we take the natural
one, associated with the parallel copy of the
torus where a given cable knot sits (its parallel copy
has zero linking number with this knot).


\subsection{\bf From knots to links}\label{sec:knots-links}
Switching to links, we need to add
{\em colors\,} to the cables above, which are  dominant
weights $b$. For knots, there is always one color, so it
gives an extra (external) piece of information on top
of the topological data above. Now adding colors becomes
an internal part of the definition.

\subsubsection{\sf Graphs and labels}
The \tax-presentation of a
{\em iterated torus link\,} will be
a union of $\kappa$ {\em colored knots\,}
\begin{align}\label{tau-link}
\l_{(\vec\rr^{\,j},\vec\ss^{\,j})}^{\,\Up,\, (b^j)}=
\Bigl(\{\t(\vec\rr^{j},\vec\ss^{j}),\, b^j\in P_+\},\,
j=1,\ldots,\kappa\Bigr) \hbox{\, together with}&\\
\hbox{\em the incidence} \hbox{\em\, matrix\,\,} \Up=(\up_{j,k}),
\hbox{\, where\, } 0\le \up_{j,k}=\up_{k,j}
\le \min\{\ell^j,\ell^{k}\}&
\notag\\
\hbox{implies\, that\,}\, [\rr_{i}^j,\ss_{i}^j]\!=\!
[\rr_{i}^k,\ss_{i}^k]\hbox{\, for all\, }
1\le i \le \up_{j,k} \hbox{\, and any\, } 1\le j,k \le \kappa&.
\notag
\end{align}
Here $\ell^j$ is the length of $\vec\rr^j=\{\rr_i^j\}$ and
$\vec\ss^j=\{\ss_i^j\}$ for $1\le j\le \kappa$; we
naturally set $\up_{j,j}=\ell^j$.

Thus $\Up$ determines a graph $\l$ with the vertices
$\{i,j\}$ identified as in (\ref{tau-link}).
The {\em paths\,} are sequences
of increasing consecutive $i$\~vertices with fixed $j$;
their orientation is from $i$ to $i+1$. The vertices for 
neighboring $i$ in the same path $j$ are naturally
connected by the {\em edges}. 
This graph is a disjoint
union of trees. Any  its subtree contains a unique 
{\em initial vertex\,} $i=i_0$ such that the  $i$\~indices 
are the distances in this subtree from $i_0$ one plus $i_0$.  
Every subtree has at least one {\em base path\,}, the one that 
intersects all other paths in this subtree.
 
For $i\le \ell^j$,
the pairs $[\rr_i^j,\ss_i^j]$ become labels, called
{\em \tax-labels\,} of the 
{\em vertices\,} $i,j$ of $\l$; the 
square brackets will be used for them.
\vskip 0.2cm

Additionally, we put {\em arrowheads\,}
at the end of {\em every\,} path (which is at the vertex with
$i=\ell^j$).
The paths with coinciding vertices but
different arrowheads will be treated as different paths.
So $j$ are the indices of all different
maximal paths in $\l$ including the arrowheads at their ends. 
The colors $b^j$ will be assigned to
the arrowheads. 
Topologically, the $j${\tiny th} path
corresponds to the knot $\t(\vec\rr^{j},\vec\ss^{j})$
colored by $b^j\in P_+$ (later, by the corresponding
Young diagrams $\la^j$).
If a graph contains no vertices;
then it is a collection of paths that are pure arrowheads
(a set of colors).


\subsubsection{\sf Topological parameters}
The $\aa$\~parameters above must be now calculated
along the paths exactly as we did for the knots
(i.e. starting from $i=1, \aa_1=\ss_1$); then
$\aa_i^j$ depends only on the corresponding vertex.
The pairs
$\{\aa_i^j,\rr_i^j\}$
will be called the {\em cab-labels} of the vertices.
Actually only the \tax-labels will be needed in the DAHA
constructions; we will call them simply {\em labels}
(and use $[\,,\,]$ only for them). However
the cab-labels are necessary for understanding the topological
symmetries.

The torus knot colored by $b\in P_+$ (or by the corresponding
$\la$) is denoted by $T_{\rr,\ss}^b$;
$C\!ab_{\aa,\rr}^b(\mathbf L)$,
equivalently $C\!ab_{0,1}^b Cab_{\aa,\rr}(\mathbf L)$, is
the cable $C\!ab(\aa,\rr)(\mathbf L)$
of a link $\mathbf L$ colored by $b$.
The color is attached to the last $C\!ab$ in
the sequence of cables.
In the absence of vertices, the notation is $\unknot^{\,b}$
(the unknot colored by $b\in P_+$) or $C\!ab(0,1)^b$. {\em
We will use the same notation $\l$ for the graph and the
corresponding link $\mathbf{L}$.}

\smallskip
The passage from the {\em \tax-presentation\,}
to the {\em cab-presentation\,}
is
\begin{align}\label{Knotsiterx}
\l\bigl(\vec\rr^{j},\vec\ss^{j},\, 1\le j\le \kappa\bigr)
\rightsquigarrow
\Bigl(\coprod_{j=1}^\kappa C\!ab(\vec{\aa}^j,\!\vec{\rr}^j)
\Bigr)(\unknot),
\end{align}
where the composition and coproduct of cables is with respect to
the graph structure and $C\!ab(\vec{\aa}^j,\!\vec{\rr}^j)=
\cdots C\!ab(\aa^j_2,\rr^j_2)T(\rr^j_1,\ss^j_1)$ is as in
(\ref{Knotsiter}).
In this work, the coproduct symbol $\,\hbox{\small$\coprod$}\,$,
which stands for the union of cables, will be simply replaced
by comma; we set $\bigl(C\!ab(\aa,\rr),C\!ab(\aa',\rr')\bigr)$
instead of $C\!ab(\aa,\rr)\coprod C\!ab(\aa',\rr')$.
The $\aa$\~parameters are calculated as above along the
corresponding paths.
See \cite{ChD2} for details and (many) examples.
Actually, we do not need much the topological aspects
in this work; the graphs are sufficient for the DAHA
construction.

\subsection{\bf Pairs of graphs} \label{sec:Splice}
This extension is necessary to incorporate {\em all\,}
iterated torus links (in this work, all {\em algebraic\,} links);
see \cite{ChD2} for details. Let 
$\{\mathcal{L},\,'\!\mathcal{L}\}$ be a pair of labeled
graphs  defined above.

\subsubsection{\sf Twisted union}
The cabling construction provides a canonical
embedding of the iterated torus links into
the solid torus. Their {\em twisted union\,} for the
pair $\{\mathcal{L},\,'\!\mathcal{L}\}$ is as follows:
we put the links for $\l$ and
$'\!\mathcal{L}$ into the horizontal solid torus
and the complementary vertical one. 

Since we consider only {\em algebraic\,} links, we will
always change here the natural orientation
of the second component by the opposite one.
Without this switch, the resulting link
is never algebraic. For instance,
\Yboxdim7pt
$\bigl\{\mathcal{L}=\{\circ\rightarrow \yng(1)\},\,
'\!\mathcal{L}=\{\circ\rightarrow\yng(1)\}\bigr\}$
represents uncolored Hopf
$2$\~links with the linking number $lk=+1$, which
is $\{\circ\rightarrow\yng(1),
\circ\rightarrow\yng(1)^\vee\}$
in \cite{ChD2}.

{\em Thus the pairs $\{\mathcal{L},\,'\!\mathcal{L}\}$ in this
work are actually $\{\mathcal{L},\,'\!\mathcal{L}^\vee\}$ in the
notation from \cite{ChD2}}; we consider only such pairs
(with $\vee$) here. We note that  
$\{\mathcal{L},\,'\!\mathcal{L}\}$ and 
$\{'\!\mathcal{L},\,\mathcal{L}\}$ result in isotopic links;
this corresponds to formula (4.24) from \cite{ChD2} in the 
DAHA setting.

\subsubsection{\sf Positivity conditions}
Arbitrary algebraic links can be obtained using
this construction including 
the {\em twisted union\,} for the pairs
of positive algebraic graphs
subject to the inequality
$\,'\ss_1 \ss_1>\,'\rr_1 \rr_1$ for any pairs of
the first vertices of
these two graphs. Then $\{\l,\,'\!\l\}$ is
called a {\em positive pair\,}. See e.g. \cite{New}
\cite{ChD2}. These inequalities are imposed upon the full
usage of the symmetries
of the corresponding {\em splice diagram\,}. In the absence 
of vertices or if there is only one vertex with $\rr_1=1$,
the pair $\{\rr_1,\ss_1\}$ is technically treated as  
$\{0,1\}$ for $\l$ (or for $'\!\l$, with primes), so
the inequality above holds tautologically. Recall
that $C\!ab(0,1)L=L$ for any link $L$ (it is a path along the
link).  Note that the transposition of $\rr_1$ and
$\ss_1$ (only in the first pair!) does not change
the isotopy type of the corresponding component, but this
may influence the resulting twisted union. Let us comment on it.

Let $'\!\l$ be a pure arrow (no vertices)
colored by $\la$ and $\l$ is a positive graph. Then 
the pair $\{\l,'\!\!\l\}$ is positive and it
corresponds to adding ``the meridian" colored by $\la$
to $\l$. The meridian and its linking number with $\l$
(always positive) obviously may change when $\rr_1,\ss_1$
are transposed.
In the DAHA construction, this pair will be treated as
taking the coinvariant of $P_{\la}(X_b\mapsto Y_b)$ applied to
the {\em pre-polynomial\,} 
associated to $\mathcal{L}$. If
the color here is $\la=\yng(1)$ (common
in our numerical examples), then
\Yboxdim5pt
$P_{\yng(1)}$ is the $\S_{n+1}$\~orbit sum of $X_{\om_1}$.
\Yboxdim7pt
\smallskip

See \cite{EN,New} and \cite{ChD2} for details.
The theory in \cite{EN} is without colors, as well as
that in \cite{ObS,ORS}. Attaching colors to the branches
can be incorporated topologically using framed links, but 
this requires more involved combinatorial definitions
(and some usage of {\em the skein\,}).  Colors are natural in the
DAHA construction, but colored DAHA superpolynomials
are of course more complicated than the uncolored ones.
See also \cite{Ma,ChP2}.

\subsubsection{\sf Algebraic links} 
We provide here only basic facts; see \cite{EN} for details and
references, especially Theorem 9.4 there.
Generally, one begins with a polynomial equation $f(x,y)=0$
considered in a neighborhood of an isolated singularity
$0=(x=0,y=0)$. Its intersection
with a small $3$-dimensional sphere in
$\C^2$ around $0$ is called an {\em algebraic link\,}.
Assuming that $\rr^j_i,\ss^j_i>0\,$,
any {\em labeled graph\,} $\l=\l^\Up_{(\vec\rr^j,\vec\ss^j)}$
(in the \tax-presentation, without colors) corresponds to a 
germ of {\em plane curve singularity\,}  at $0$. If these
inequalities hold, the graph is called {\em positive\,}.
Using positive pairs  $\{\l,'\!\!\l\}$ provide all of them.

The corresponding (germs of) singularities are unions of  {\em
unibranch\,} components corresponding for the paths of $\Up$
(numbered by $j$),  which are given as follows:
\begin{align}\label{yxcurve}
y = c^j_1\,x^{\ss^j_1/\rr^j_1}
(1+c^j_2\,x^{\ss^j_2/(\rr^j_1\rr^j_2)}
\bigl(1+c^j_3\,
x^{\ss^j_3/(\rr^j_1\rr^j_2\rr^j_3)}
\Bigl(\ldots\Bigr)\bigr)) \hbox{\, at\, } 0.
\end{align}
The numbers $\rr,\ss$ are obtained from the corresponding
labels; the parameters $c_i^j\in \C$ must be
sufficiently general here.
The simplest example is the equation
$y^{\rr\kappa}= x^{\ss\kappa}$ under gcd$(\rr,\ss)=1$,
which corresponds to the torus link $T(\rr\kappa,\ss\kappa)$
with $\kappa$ knot components isotopic to $T(\rr,\ss)$. The
pairwise linking numbers here are all equal to
\,$\rr\,\ss$\, in this case.

The unibranch components and their
(pairwise) linking numbers uniquely determine the corresponding
germ due to the Reeve theorem; see e.g. \cite{EN}.
All linking numbers must be strictly positive for
algebraic links. The DAHA constructions works
for any (not only positive) labels. The above discussion
is of course in the absence of colors.

\subsection{\bf DAHA-Jones polynomials}
They can be uniformly defined for any (reduced,
irreducible) root systems $R$; for its twisted
affinization $\tilde{R}$, to be exact.  We need only 
$R=A_{n}$ in this work. The $P,J$\~polynomials
and the necessary DAHA tools are from the previous sections.

\subsubsection{\sf Data and ingredients}
The combinatorial data will be the \tax-labeled graphs
$\l_{(\vec\rr^{j},\vec\ss^{j})}^{\,\Up, (b^j)}$
from (\ref{tau-link}) and their pairs. Recall that
$$
1\le j\le \kappa,\,\, \vec\rr^{\,j}=\{\rr_i^j\},\,\,
\vec\ss^{\,j}=\{\ss_i^j\},\,\, 1\le i\le \ell^j,
$$
and $\Up$ is the {\em incidence graph/matrix\,};
the {\em arrowheads\,} (at the ends of all {\em paths\,})
are colored by $b^j\in P_+$. The incidence graph is
not supposed to be connected here and
the paths can contain no vertices; see
(\ref{tau-link}). The construction below will be for
{\em two\,} arbitrary such graphs $\l$, $'\!\l$ (the second
can be empty).

In the case of algebraic {\em knots\,}, spherical polynomials
$P_\la^\circ$ are sufficient; one obviously
has $J_\la^\circ\equal J_\la/J_\la(t^\rho)=P_\la^\circ$.
Generally (for links), we need the $J$\~polynomials;
see (\ref{P-arms-legs}). For the latter reference,
let $\la^j=\la(b^j)$ for dominant $b^j$. We set:
\begin{align}\label{P-arms-legs-LCM}
&(b^1,\ldots,b^m)_{ev}^J\!=\!(\la^1,\ldots,\la^m)_{ev}^J\!=
\!LC\!M\bigl(J_{\la^1}(t^{\rho}),\ldots,
J_{\la^m}(t^{\rho})\bigr),
\end{align}
where $LC\!M$ is normalized by the condition that it is a
$q,t$\~polynomial with the constant term $1$.

One has the following combinatorially transparent formula:
\begin{align}\label{P-arms-legs-union}
&(\la^1,\ldots,\la^m)_{ev}^J\ =\
(\la^1\!\!\vee\cdots\vee\!\la^m)_{ev}^J\,, \hbox{\,\,\, where}\\
&\la^1\!\!\vee\cdots\vee\!\la^m
\hbox{\, is the union of diagrams\, }
\{\la^j\}.\notag
\end{align}
The latter union is by definition the smallest Young diagram
containing all diagrams $\la^1,\ldots,\la^m$.

Note that the $J$\~polynomials in the $A_n$\~case
are not {\em minimal\,} integral (without $q,t$\~denominators)
even for $t=q$. They are important for the $a$\~stabilization
(including HOMFLY-PT polynomials for {\em links\,})
and for the super-duality. However, the switch from $P$ to $J$
does not influence the DAHA-construction for iterated knots, 
though their role is still important even for knots; 
see \cite{GoN},\cite{ChD2}. 
\smallskip

Let us now go to the DAHA construction.
Recall that $H\!\!\Downarrow\, \equal H(1)$, where the
action of $H\in \HH$ in $\v$ is used.
We represent torus knots $T(\rr,\ss)$
by the matrices $\ga[\rr,\ss]=\ga_{\rr,\ss}\in PSL_{\,2}(\Z)$
with the
first column $(\rr,\ss)^{tr}$ ($tr$ is the transposition)
for $\,\rr,\ss\in \N$ provided \,gcd$(\rr,\ss)=1$.
Let $\hat{\ga}_{\rr,\ss}\in PSL_{\,2}^{\wedge}(\Z)$ be
{\em any\,} pullback of $\ga_{\rr,\ss}$ to the projective
$PSL_{\,2}(\Z)$.

\subsubsection{\sf Pre-polynomials}\label{JONITER}
The definition is for any pair $\{\l,'\!\!\l\}$ from
(\ref{tau-link}) (the {\em positivity\,} of $\l,'\!\!\l$
and the pair is 
actually not needed). Let
\begin{align}\label{links-x-y}
\l=\l_{(\vec\rr^j,\vec\ss^j)}^{\Up,(b^j)},\
'\!\l=\, '\!\l_{('\vec\rr^j,'\vec\ss^j)}^{\,'\Up,('\!b^j)}
\where b^j,\,'\!b^j\in P_+,
\end{align}
$$
1\le j\le \kappa,\,'\!\kappa\for \l,\, '\!\l,\
\vec\rr^j= (\rr_i^j \mid 1\le i\le \ell^j),\,
'\vec\rr^j=('\!\rr_i^j \mid 1\le i\le\, '\!\ell^j).
$$

First of all, we lift $(\rr^j_i,\ss^j_i)^{tr}$,
 $('\!\rr^j_i,'\!\ss^j_i)^{tr}$\, to\,
$\hat{\ga}_i^j,\, '\hat{\ga}_i^j\in PSL_{\,2}^{\wedge}(\Z)$
as above.
Then for any vertex of $\l$, i.e. a  pair $\{i,j\}$, we
begin the (inductive) definition with:
$$
\P_{\ell^j+1}^j\equal J_{b^j},\
\ga_{\ell^j+1}^j\equal\hbox{id} \for 1\le j\le \kappa.
$$
Recall that $\ell^j=0$ when the $j$\hbox{\tiny th}
path contains only an arrowhead and $\up_{j,k}=0$
if the corresponding paths do not intersect.
\smallskip

For a given {\em path} with the index $j$, we
define the {\em pre-polynomials\,} $\P_{i}^{j}$
by induction with respect
to $i$, starting with $i=\ell^{j}$ down to $i=0$:
\begin{align}\label{hatPij}
\P_{i}^{j}=
\prod_{1\le k\le \kappa}^ {\up(k,j)=\,i}
\bigl(\,\hat{\ga}_{i+1}^{k}(\P_{i+1}^{k})
\!\!\Downarrow\,\bigr).
\end{align}
I.e. the last product is over all paths $k$ that enter
(intersect) the path for $j$ exactly at the index
$i\ge 0$, including $k=j$\, when\, $i=\ell^j$.
Note that
$\P_{\ell^{j}}^{j}=\prod_{1\le k\le
\kappa}^ {\up(k,j)=\ell^j} J_{b^k}$ for a {\em base path\,}
$j$, where this product is over all $\kappa$
{\em arrowheads\,} from (originated at) the
vertex $\{i=\ell^j,j\}$.

The polynomial $\P_{0}^{j}$ actually depends only on the
corresponding {\em subtree\,} for
{\em any\,}
path $j$\, there. If $\l$ is a union of subtrees, then
$\P_{0}^{j}$
is the product of the corresponding polynomials
$\bigl(\,\hat{\ga}_{1}^{k}
(\P_{1}^{k})\!\!\Downarrow\,\bigr)$
over all these subtrees. 
The {\em pre-polynomial\,} is defined then as 
$\,\P_{0}^{tot}=\P_{0}^{j}$ (the end of the above 
inductive procedure).
The pre-polynomial $'\P_{0}^{tot}$ for $'\!\l$
is defined in the same way.

\subsubsection{\sf Finale}
Using the notations $\mathbf{b}=(b^j),\,\mathbf{'b}=('b^j)$,
the {\em DAHA-Jones polynomial\,} for the $J$\~polynomials
$J_{b^j}$, $J_{'\!b^j}$  and a fixed index $1\le j_o\le \kappa$
(which determines the normalization) is as follows:
\begin{align}\label{jones-hat}
J\!D_{\,(\vec\rr^j,\,\vec\ss^j)\,,
\,('\!\vec\rr^j,\,'\!\vec\ss^j)}
^{R,\,j_o,\,\Up,\,'\!\Up}((b^j),('\!b^j)\, ; \,q,t)\, =&\,
J\!D^{j_o,\,\Up,\,'\!\Up}_{\,(\vec\rr^j,\,\vec\ss^j)
\, , (\,'\!\vec\rr^j,\,'\!\vec\ss^j)}
(\mathbf{b},\mathbf{'b}\, ; \,q,t)\\
=\,J\!D^{j_o}_{\l,\,'\!\l}\equal
\Bigl\{\,'\P_0^{tot}(Y)&\,\P_0^{tot}/
J_{b^{j_o}}(q^{-\rho_k})\,\Bigr\}_{ev}.
\notag
\end{align}
Here $j_o$ can be $\emptyset$, which means that there is no
divisions by the evaluations at $q^{-\rho_k}=t^{-\rho}$;
$R$ is the root system, which is $A_n$.
\smallskip

In the case of iterated torus {\em knots\,} (when there is
only one path) and in the absence of $'\!\l$,
we arrive at formula (2.12) from \cite{ChD1}:
\begin{align}\label{jones-ditx}
& J\!D_{\vec\rr,\vec\ss}
(b;q,t)\! =\!
\Bigl\{\hat{\ga_{1}}\Bigr(
\cdots\Bigl(\hat{\ga}_{\ell-1}
\Bigl(\bigl(\hat{\ga}_\ell(P_b)/
P_b(q^{-\rho_k})\bigr)\!\Downarrow
\Bigr)\!\Downarrow\Bigr) \cdots\Bigr)\Bigr\}_{ev}.
\end{align}
It includes only one $b\in P_+$ and therefore one can
use $P_b$ instead of $J_\la$.

The simplest link is for the union 
of any number of arrowheads
colored by $b^1,\ldots, b^\kappa$. It is represented
by the graph $\rightdotsarrow$. A single $\rightarrow$ is
associated with 
$P_b/P_b(q^{-\rho_k})$; generally we arrive at the product: 
\begin{align}\label{simp-link-P}
\P_{\ell}^1=
J_{\ell+1}^1\cdots J_{\ell+1}^\kappa
/J_{b^{j_o}}(q^{-\rho_k})=
J_{b^1}\cdots J_{b^\kappa}
/J_{b^{j_o}}(q^{-\rho_k}),
\end{align}
where $j_\circ$ can be $\emptyset$. It will be later allowed
to divide by the $LC\!M$ of the evaluations for all $J_{b^j}$
here, for the {\em minimal normalization\,} from
(\ref{jones-bar}). 


\subsection{\bf Polynomiality etc}
The following theorem and other statements in this
and the next sections are
from \cite{ChD1} and \cite{ChD2}.

\begin{theorem}\label{THM-integr-Jones}
For any choice of the normalization
index $1\le j_o\le \kappa$ (it can be $\emptyset$),
the DAHA-Jones polynomial \,
$J\!D^{j_o}_{\l,\,'\!\l}\,$
is indeed a polynomial in terms of $q,t$
up to a factor $q^\bullet t^\bullet$,
where the powers $\bullet$ can be rational.
Modulo such factors,
it does not depend
on the particular choice of the lifts \,$\ga^j_i\in PSL_2(\Z)$
and $\hat{\ga}^j_i\in PSL_{\,2}^{\wedge}(\Z)$\, for $1\le  i
\le\ell^j$. 

Up to the $q^\bullet t^\bullet$\~equivalence
and $\Q$\~proportionality
we pick the {\sf hat-normalization}, denoted
$\hat{J\!D}^{j_o}_{\l,\,'\!\l}\,$,
as follows.
It is a $q,t$\~polynomial not divisible by
$\,q\,$ and by $\,t\,$, with integral coefficients of all
$q,t$\~monomials such that their $GCD$ is $1$ and, finally,
the coefficient of the minimal pure
power of $t$ is assumed positive. \sq
\end{theorem}

\subsubsection{\sf Minimal normalization}
The $q,t$\~integrality and other claims
from Theorem \ref{THM-integr-Jones} hold for the
following modifications of DAHA-Jones polynomials
(which is sharper and does not require picking $j_o$).
Make $j_o=\emptyset$ (no normalization in $\hat{\P}$),
we set:
\begin{align}\label{jones-bar}
\hat{J\!D}^{min}_{\,(\vec\rr^j,\,\vec\ss^j)
\, , (\,'\!\vec\rr^j,\,'\!\vec\ss^j)}
(\mathbf{b},\mathbf{'b}\, ; \,q,t)
=\,\hat{J\!D}^{min}_{\l,\,'\!\l}\equal
\Bigl\{\frac{\, \vph\hbox{\tiny $\circ$\,}
\iota\,(\,'\hat{\P}_0^{tot})\,
\hat{\P}_0^{tot}}{
(\mathbf{b},\mathbf{'b})_{ev}^J}\Bigr\}_{ev},
\end{align}
in the notation from (\ref{P-arms-legs-union}),
(\ref{iotaXY}).
We put $j_o=min$ if the minimal normalization
(division by the corresponding $LCM$) 
is taken.
\smallskip

Some modification is needed for $t=1$ to ensure the
connection with the HOMFLY-PT polynomials. We take
$P^{(\!k=1\!)}_b$ instead of $J_\la$
for $\la=\la(b)$ and $t\!=\!q$,
which do not depend on $q$ at all; they simply coincide with
the corresponding Schur functions. The $LCM$
$(b^1,\ldots,b^m)_{ev}$ then must be understood correspondingly.
Recall that the polynomials $J_\la$ generally have extra factors
vs. the Schur functions in this case.

\subsubsection{\sf Topological symmetries}
The polynomial $\hat{J\!D}^{j_o}_{\l,\,'\!\l}$ defined in
Theorem \ref{THM-integr-Jones}\,
and $\hat{J\!D}^{min}_{\l,\,'\!\l}$ introduced above
actually depend
only on the topological link corresponding to the pair of graphs
$\{\l,\,'\!\l\}$.
For instance, the reduction of the
vertices with $\rr=1$ can be applied in
$\l$ or in $'\!\l$. Also, the transposition
$[\rr_1^j,\ss_1^j]\mapsto $ $[\ss_1^j,\rr_1^j]$  (only
for $i=1$) does not influence
$\hat{J\!D}^{j_\circ}_{\l}$ or $\hat{J\!D}^{min}_{\l}$
provided that $\,'\!\l=\emptyset$, and the pairs
$\{\l,'\!\l\}$ and $\{'\!\l,\l\}$ result in coinciding 
polynomials.
\vskip 0.2cm

The justification of this and other symmetries is
essentially parallel to Theorem 1.2 from
\cite{CJJ}. Let us discuss torus knots.
Essentially, one needs to check here that 
$T(m\rr+\ss,\rr)$ results in the same
DAHA-Jones polynomial as the ``2-cable" corresponding to
the tree $[m,1]\rightarrow [\rr,\ss]$. Topologically,
$T(m\rr+\ss,\rr)$ is isotopic to
$C\!ab(m\rr+\ss,\rr)T(m,1)$, since
$T(m,1)$ is unknot. The corresponding
relation for the $J\!D$\~polynomials readily
follows from the commutativity $\tau_-^m$
with $\Downarrow$, which simply means that $\tau_-$ acts
in $\v$.


\subsubsection{\sf Specialization $q=1$}
We now make $q=1$, assuming that $t$ is generic and
using the notation $(b^1,\ldots,b^m)_{ev}^J$ from
(\ref{P-arms-legs-LCM}). Then
\begin{align}\label{q-1-prod}
&\frac{(b^1,\ldots,b^\kappa,\,'b^1,\ldots,
'\!b^{\,'\kappa})^J_{ev}}
{(b^1)^J_{ev}\cdots (b^\kappa)^J_{ev}
\,('b^1)^J_{ev}\cdots
('b^{\,'\kappa})^J_{ev}}(q\!=\!1)\,\,
\hat{J\!D}^{min}_{\l,\,'\!\l}(q\!=\!1)\\
=\ \,&
\prod_{j=1}^{\kappa} J\!D_{\,\vec\rr^j,\,\vec\ss^j}
\,\bigl(\,b^j\,;\,q\!=\!1,\,t)\
\prod_{j=1}^{'\kappa} J\!D_{\,'\vec\rr^j,\,'\vec\ss^j}
\,\bigl(\,'\!b^j\,;\,q\!=\!1,\,t),\notag\\
&\hbox{where\,\,}
J\!D_{\,\vec\rr,\,\vec\ss}
\,\bigl(b;\,q\!=\!1,t
\bigr)\!=\!
\hbox{\small$\prod$}_{p=1}^n J\!D_{\,\vec\rr,\,\vec\ss}\,
(\om_p;\,q\!=\!1,t)^{b_p},\notag
\end{align}
for $b=\hbox{\small$\sum$}_{p=1}^n b_p \om_p\in P_+$,
where the $J\!D$\~polynomials from
(\ref{jones-ditx}) are used (for knots). See formula (2.18)
in \cite{ChD1}.

Notice the factor in the left-hand side. It would be $1$
if the polynomials $P_b^\circ$ were taken in this construction
(without any further division) instead of $J_\la$ under
the minimal normalization. This can result in the 
$q,t$\~denominators of the resulting coinvariants, but
there will be no correction factor in the left-hand side.

\section{\sc DAHA-superpolynomials}
\subsection{\bf Existence and duality}
\subsubsection{\sf Stabilization}
Following \cite{CJ,GoN,CJJ,ChD1,ChD2}, the construction from
Theorem \ref{THM-integr-Jones} and other statements above can
be extended to the {\em DAHA- superpolynomials\,},
the result of the stabilization  of
$\hat{J\!D}^{A_n,j_o}_{\l,\,'\!\l}$
(including $j_o=min$, the minimal normalization).

The $a$\~stabilization for torus knots
was announced in \cite{CJ}; its proof was published in
\cite{GoN}.
Both approaches use \cite{SV}.
The super-duality conjecture was proposed in \cite{CJ}
(let us also mention \cite{GS})
and proven in \cite{GoN} for torus knots; see also
\cite{CJJ} for an alternative approach based on the
generalized level-rank duality. The justifications
of the $a$\~stabilization and
the super-duality was extended to
arbitrary iterated torus knots in \cite{ChD1}.

The main change for links vs. knots is that the
polynomiality of the superpolynomials for links
requires the usage of  $J_\la$\~polynomials.
Actually $\{J_\la\}$ were already employed in \cite{GoN}
for the stabilization and duality, but the construction
of (reduced) $J\!D$\~polynomials and superpolynomials for
{\em knots\,} requires only spherical $P^\circ_\la$. For
{\em links\,} (vs. knots), the role of polynomials
$\{J_\la\}$ is the key; without the usage
of $J$\~polynomials
the superpolynomials have non-trivial $t$\~denominators.

\smallskip

The sequences $\vec\rr^j,\,\vec\ss^j$ of length $\ell^j$
for the graph $\l$ and $'\!\vec\rr^j,'\!\vec\ss^j$ of
length $'\!\ell^j$ for the graph $'\!\l$ will be from the
previous sections. We always
use the DAHA-Jones polynomials under the hat-normalization,
i.e. $\hat{J\!D}^{j_o},\, \hat{J\!D}^{min}$.
For $t\!=\!q$, the Schur functions $P_{\la}^{(k\!=\!1)}$
will be employed when discussing the
connection with the HOMFLY-PT polynomials.
Recall that $\la=\la(b)$ is the Young diagram
representing $b\in P_+$.
\smallskip

We consider now $P_+\ni b=$ $\sum_{i=1}^n b_i \om_i$
for $A_n $ as
a (dominant) weight for {\em any\,} $A_m$ (for $sl_{m+1}$)
for $m\ge n-1$; we set $\om_{n}=0$ upon its restriction
to $A_{n-1}$.
See \cite{CJ,GoN,CJJ,ChD1} concerning the versions
of the following theorem for torus knots and
iterated torus {\em knots}.


\begin{theorem}\label{STABILIZ}\cite[Theorem 2.3]{ChD2}
Given a pair $\{\l,\, '\!\l\}$
colored by $\mathbf{b}=(b^j),\, '\mathbf{b}=('\!b^j)$ and
the normalization index $1\le j_o\le \kappa$ or
$j_o\!=\!min$, there exists a
unique polynomial from $\Z[q,t^{\pm 1},a]$
\begin{align}\label{h-polynoms-hat}
&\hat{\h}^{j_o}_{\l,\,'\!\l}\ =\
\hat{\h}^{\,\Up,j_o}_{(\vec\rr^j,\,\vec\ss^j),
('\vec\rr^j,\,'\vec\ss^j)}(\mathbf{b},\,
'\mathbf{b};\,q,t,a)
\end{align}
such that for any  $m\!\ge\! n\!-\!1$ and proper powers
of $q,t$ (possibly rational)\,:
\begin{align}\label{jones-sup-hat}
&\hat{\h}^{j_o}_{\l,\,'\!\l}(q,t,a\!=\!-t^{m+1})\,=\,
\pm\, q^\bullet t^\bullet\,
\hat{J\!D}^{A_m,j_o}_{\l,\,'\!\l}(q,t).
\end{align}
Let us pick $\hat{\h}$ such that  $\hat{\h}(a\!=\!0)$
is under the hat-normalization from Theorem
\ref{THM-integr-Jones}.
Then  relations (\ref{jones-sup-hat}) will automatically
hold for sufficiently large $\,m$
without any correction factors $\pm q^\bullet t^\bullet$
(and one sufficiently large
$\,m$\, is actually sufficient to fix $\,\hat{\h}\,$
uniquely). \sq
\end{theorem}

\subsubsection{\sf Symmetries}
The polynomials $\hat{\h}$ depend only on the isotopy
class of the corresponding iterated torus links.
All symmetries from the previous section
hold for $\hat{\h}=\hat{\h}^{j_o}$, including $j_o=min$.
For instance, this includes
the {\em specialization\, }
relation from (\ref{q-1-prod}) at $q=1$
for $\hat{\h}$. The exact product
formula there holds when spherical polynomials $P_b^\circ=
P_b/P_b(t^{\rho})$ are used in the formulas for $\hat{J\!D},
\hat{\h}$ instead of $J$\~polynomials.
We note that {\em mirroring\,} of torus iterated links results in  
changing $\aa_i\mapsto -\aa_i$ for all $i$ and 
in $q\mapsto q^{-1}, t\mapsto t^{-1}, a\mapsto a^{-1}$
in the superpolynomials (followed by the hat-normalization).
Also, algebraically $Y$ must be replaced by $Y^{-1}$ in the 
case of two graphs (the change or their relative orientation);
 see \cite{ChD2}.
\vskip 0.2cm

The key new feature of the stable theory
is  {\em super-duality\,}.
We switch from $\mathbf{b},\,'\mathbf{b}$ to the
corresponding sets of Young diagrams
$\{\la^j\}$,\, $\{'\!\la^j\}$; their transpositions
will be denoted by $\{\cdot\}^{tr}$.
Up to powers of $q$ and $t$
denoted here and below by $q^\bullet t^\bullet$, one has:
\begin{align}\label{iter-duality}
\hat{\h}_{\l,'\!\l}(\{\la^j\},\{'\!\la^j\}; q,t,a)\!=\!
q^{\bullet}t^{\bullet}
\hat{\h}_{\l,'\!\l}(\{\la^j\}^{tr},
\{'\!\la^j\}^{tr};t^{-1},q^{-1},a).
\end{align}

Let us discuss the $a$\~degree of $\hat{\h}$\~polynomials.
\comment{
We will assume
that $\rr_i^j,\,'\!\rr_i^j\neq 0$ for $i\!>\! 1$.
Then \,
deg${}_a \hat{\h}^{min}_{\l,\,'\!\l}
(\,\lla,\,'\!\lla\,;\,q,t,a)\,$  is no greater than
\vskip -0.3cm
\begin{align}\label{deg-a-j}
\sum_{j=1}^\kappa
\max\{1,|\,\rr^j_1|\}|\rr^j_2\cdots\rr^j_{\ell^j}|\,|\la^j|+
\sum_{j=1}^{'\!\kappa}
\max\{1,|\,'\rr^j_1|\}|\,'\rr^j_2\cdots\,'\rr^j_{\,'\!\ell^j}|
\,|'\!\la^j| -\De,&\notag\\
\for \De\!=\!|\la^1\!\vee\!\ldots\!\vee\!\la^\kappa\!\vee\,
'\!\la^1\!\vee\!\ldots\!\vee\, '\!\la^{'\!\kappa}|
 \hbox{\, for \,} \hat{\h}^{min},\
\De\!=\!|\la^{j_o}| \hbox{\, for \,} \hat{\h}^{j_o},&
\end{align}
where $|\la|$ is the number of boxes in $\la$.
\smallskip
}
We conjectured in \cite{ChD2} (for algebraic links only) that
\begin{align}\label{deg-a-jj}
\hbox{deg}_a\hat{\h}^{min}_{\l,\,'\!\l}
=&\,\hbox{\small $\sum$}_{j=1}^\kappa\
\min\{\,\rr^j_1,\ss^j_1\}\, \rr^j_2\cdots
\rr^j_{\ell^j}\,|\la^j|\!\\
&+
\hbox{\small $\sum$}_{j=1}^{\,'\!\kappa}
\ \min\{\,'\rr^j_1,\,'\ss^j_1\}\,
\,'\rr^j_2\cdots
\,'\rr^j_{\,'\!\ell^j} \,|\,'\!\la^j|-\De\notag\\
\hbox{for\ } \De&=\!|\la^1\!\vee\!
\ldots\!\vee\!\la^\kappa\!\vee\,
'\!\la^1\!\vee\!\ldots\!\vee\, '\!\la^{'\!\kappa}|
\hbox{\, from \,(\ref{P-arms-legs-union})},
\notag
\end{align}
where $|\la|$ is the number of boxes in $\la$.
A somewhat weaker statement can be justified. This is
connected with the product formula at $q=1$, which gives that
the $a$\~degree is no smaller than that in (\ref{deg-a-jj}).
The right-hand side of this
formula is actually the multiplicity
of the corresponding singularity generalized to the colored
case, e.g. $|\la|(\ss-1)$ for
$T(\rr,\ss)$ colored by $\la$, where $\ss<\rr$. 

From now on we will always impose the
{\em minimal normalization\,} of $\hat{\h}$
unless stated otherwise.

\subsubsection{\sf Family superpolynomials}\label{SEC:FAM}
There are important reasons to consider superpolynomials
in the ``families" of links. Namely,
$(a)$ the corresponding generating
functions  have better algebraic structure than
individual superpolynomials, $(b)$ DAHA provides uniform tools
for calculations within families (see below), 
$(c)$ the families are natural to match
the classical theory of zeta and \cite{ChZ}, 
$(d)$ considering families is 
similar to Iwasawa theory in number theory 
(see \cite{Mor} and below).

Generally, the definition is as follows. One 
replaces one (or several) $\ga_i$ in the formula for 
DAHA superpolynomials by  $\sum v^m\tau_{\pm}^m\ga_i$ for
$m\in \Z_+$ and a new variable/parameter $v$. 
When $i=1$ in the case
of uncolored torus iterated {\em knots\,} described by  formula 
(\ref{jones-ditx}), we set
\begin{align}\label{famiter}
&\mathbb{JD}_{\vec\rr,\vec\ss}(b; q,t,u)=
\sum_{m=0}^\infty \left(\frac{u}{t}\right)^{g(m)}\, J\!D_m\\
=\sum_{m=0}^\infty
\left(\frac{u}{t}\right)^{g(m)}& 
\Bigl\{\hat{\tau_-^m \ga_{1}}\Bigr(
\cdots\Bigl(\hat{\ga}_{\ell-1}
\Bigl(\bigl(\hat{\ga}_\ell(P_b)/
P_b(q^{-\rho_k})\bigr)\!\Downarrow
\Bigr)\!\Downarrow\Bigr) \cdots\Bigr)\Bigr\}_{ev}, \notag
\end{align}
where $g(m)$ is the arithmetic genus of the 
corresponding singularity; this is $\de$
from Section \ref{SEC:MOTZ}. Combinatorially, the pair
$\rr_1,\ss_1$ is replaced by $\rr_1,m\rr_1+\ss_1$.
The corresponding rings (see below) are naturally embedded
for $i=1$. 
The division by $t^{g(m)}$ provides the exact
super-duality $q\leftrightarrow 1/t$ (not only up to
proportionality). The corresponding {\em family
superpolynomial\,}
is defined accordingly. For $\hat{\h}_m$ corresponding
to $J\!D_m$, 
\begin{align}\label{hsfam}
\hat{\mathbb{H}}(q,t,a,u)\equal
\sum_{m=0}^\infty \left(\frac{u}{t}
\right)^{g(m)} \hat{\h}_m(q,t,a).
\end{align}

The simplest example is for the family $T(2,2m+1)$, i.e.
for the summation over $\tau_-^m\tau_+\tau_-$ with
$m\in \Z_+$. The formula is as follows
\begin{align}\label{2z+1}
\hat{\mathbb{H}}_{2,1+2\Z_+}(q,t,a,u)=
\frac{1+auq/t}{(1-u/t)(1-uq)},
\end{align}
which is simple to establish; the superpolynomials
for $T(2,2m\!+\!1)$ are well known. From the DAHA perspective,
this is the simplest application of the formula for the
{\em pre-polynomial\,} $\P_{\!2,1}=(\tau_+\tau_- (P_{\om_1}))(1)$.
Importantly the polynomials $P_b$ in (\ref{famiter}) can be
replaced by the corresponding
{\em nonsymmetric\,} polynomials $E_b$ for $b\in P_+$.
They are eigenfunctions of $Y$\~operators:
$Y_a(E_b)=q^{-(a,b)}t^{-(\rho,b)}E_b$ (for dominant $b$);
the normalization
is $E_b=X_b$ modulo lower terms. We use that
$\frac{P_b}{P_b(t^{-\rho})}=
\mathscr{P}\left(\frac{E_b}{E_b(t^{-\rho})}\right)$ 
for $\mathscr{P}=
\frac{\sum_{w\in W}t^{l(w)/2}T_w}{\sum_{w\in W} t^{l(w)}};\,$
$\mathscr{P}$ commutes with taking the coinvariant.
See Section \ref{SEC:POLR} and \cite{C101,CJ,CJJ}.
Also we can switch here to $GL_n$ instead of $A_n$.
Thus, we need
$\mathbf{E}_{2,1}=(\tau_+\tau_- (E_{\om_1}))(1)$, which
is proportional to $X_{\om_1}=X_1$. Using the standard $GL_n$ 
variables $X_i$ for $1\le i\le n$, 
$X_{\om_i}=X_1 X_2\cdots X_i$ and
$$
E_{2\om_1}\!=\!X_1^2\frac{1\!-\!q}{1\!-\!qt} +
\frac{(1\!-\!t)q}{1\!-\!qt}
(X_1\!+\!\ldots\!+\!X_n)X_1,\ 
E_{2\om_1}(t^{-\rho})\!=\!\frac{1\!-\!qt^n}{t^{n-1}(1\!-\!qt)}. 
$$ 
Applying $\mathscr{P}$ here, we arrive at some DAHA-Rosso-Jones
relation.
\Yboxdim5pt
\begin{proposition}\label{FAM1-1}
Let $\ga_1=\ga_{\rr_1,\ss_1}$ be a matrix from $PSL(2,\Z)$, 
$\ga_2=\tau_+\tau_-=\ga_{2,1}$, 
$\hat{\h}^{\yng(1)}_{\ga_1,\ga_2}$ 
the corresponding (hat-normalized)
uncolored superpolynomial, corresponding to the 
cable $C\!ab(2\rr_1\ss_1\!+\!1,2)
T(\rr_1,\ss_1)$.
Then
\begin{align}\label{fga1-1}
\hat{\h}^{\yng(1)}_{\ga_1,\ga_2}=\frac{1+aq}{1-q}
\hat{\h}_{\rr_1,\ss_1}^{\yng(2)}-
\frac{q}{1-q}\hat{\h}_{2\rr_1,2\ss_1}^{\yng(1)} \for
\rr,\ss\ge 0,
\end{align}
where $\hat{\h}_{\rr,\ss}^{\yng(2)}$ is the
superpolynomial for $\ga=\ga_{\rr,\ss}$ colored by $2\om_1$,
$\hat{\h}_{2\rr,2\ss}^{\yng(1)}$ is the
uncolored 
superpolynomial for the link $2T(\rr,\ss)$. \sq
\end{proposition}
When $\ga_1\!=\!\tau_-^m\!=\!\ga_{1,m}$\, for\, $m\!\in\!\Z_+$,
we use (\ref{taumEX}) below and obtain:
\begin{align}\label{f1-1}
&\hat{\h}^{\yng(1)}_{2,2m\!+\!1}\!=\!\frac{1+aq}{1-q}\!-\!
\frac{q}{1-q}\hat{\h}_{2,2m}^{\yng(1)};
\hbox{\,\ e.g. for\, } m\!=\!1:\\
&1+qt+aq\,=\frac{1+aq}{1-q}-
\frac{q}{1-q}\bigl(1 + aq - t(1-q)\bigr).\notag
\end{align}
Generally, this approach gives that
$\mathbb{H}(q,t,a,u)$ are rational functions in terms of 
$q,t,u,a$ (the details will be published elsewhere). 

Corollary \ref{FAM1-1} is a typical recurrence
relation in the DAHA theory of superpolynomials; 
topological justification of such DAHA formulas
within the Khovanov -Rozansky
theory is an obvious challenge. 
The calculation
of $\hat{\h}_{2\rr_2,2\ss_2}^{\yng(1)}$ can be generally
performed using the Pieri rules from (1.4.50),\cite{C101}.
See also Proposition \ref{FAM1-m} below.

Let us provide the formula for $\mathbb{H}(q,t,a,u)$
in the case of the family $T(3,3n+1)$, where $g(n)=3n$ and
the superpolynomials are sufficiently well known;
see e.g. Conjecture 7 from \cite{GORS} and \cite{CJ}. One has:

\begin{align}\label{3z+1}
\hat{\mathbb{H}}_{3,1+3\Z_+}&(q,t,a,u)=
\frac{1}{(1-(u/t)^3)(1-(uq)^3)(1-u^3q^2/t^2)}\notag\\
\times  &\Bigl( 1+u^3(\frac{q}{t^2}+\frac{q^2}{t})+
a^2u^3\,\frac{q^3}{t^3}\left(1+u^3(\frac{q}{t^2}+\frac{q^2}{t})
\right)\notag\\
&+ au^3\,\frac{q}{t}\left(\frac{1}{t^2}+\frac{q}{t}+
q^2+\frac{q}{t^2}+\frac{q^2}{t}+u^3\frac{q^3}{t^3}\right)
\Bigr).
\end{align}

\comment{
$\mathscr{P}=
\sum_wt^{l(w)/2}T_w/\sum_w t^{l(w)}$
$$
\mathscr{P}(X_1^2)=
\mathscr{P}(\e_{2\om_1})\frac{1-q t^n}{t^{n-1}(1-q)}
-\frac{(1-t)q}{(1-q)t^{n-1}}
(X_1+\ldots+X_n)\mathscr{P}(X_1).
$$ 
}

{\sf Family  C\!ab(13+2m,2)T(3,2).}
These cables can be obtained for  $\ga_1=\tau_+\tau_-^2$
and $\ga_2=\tau_-^m\tau_+\tau_-$, i.e. 
the procedure $\ga\mapsto \tau_-^m\ga$ is
used her for $\ga_2$. 
We use the above formula for $E_{2\om_1}$
rewriting it as follows:
$$
X_1^2=E_{2\om_1}\!-\!\frac{(1\!-\!t)q}{1\!-\!qt}
(X_2\!+\!\ldots\!+\!X_n)X_1.\ 
$$ 
Formula (\ref{taumin-poly}) defines the action
of $\tau_-$ in $\v$. We used in \cite{ChD2}
the notation $\dot{\tau}_-$ to avoid confusion
with the action of $\tau_-$ in $\HH$. Formula
(1.37) there states that 
$\tau_-(E_b)=q^{-b_+^2/2}t^{-(\rho,b_+)}$,
where $b_+$ is the unique weight in $P_+\cap W(b)$
for $b\in P$.  We need here the relations
\begin{align}\label{taumEX}
\tau_-(E_{2\om_1})\!=\!q^{-2}t^{1\!-\!n}E_{2\om_1}, \  
\tau_-(X_i X_j)\!=\!q^{-1/2}t^{(1\!-\!n)/2},  
1\!\le\! i\!<\!j\le n.
\end{align}
One has: $\tau_-^m(X_1^2)=$
$$\frac{(q^{-2}t^{1-n})^m}{(1\!-\!qt)}
\left(
\frac{E_{2\om_1}}{E_{2\om_1}(t^{-\rho})}
\frac{(1\!-\!qt^n)}{t^{n-1}} - (1\!-\!t)q(qt)
\,X_1(X_2\!+\!\ldots\!+\!X_n)\right).
$$

\begin{proposition}\label{FAM1-m}
For any $m\in\Z_+$, $\rr_1,\ss_1\ge 0$,
let $\ga_1\!=\!\ga_{\rr_1,\ss_1}\!\in\! PSL(2,\Z)$, 
$\ga_2(m)\!=\!\tau_-^m\tau_+\tau_-=\ga_{2,2m+1}$, 
$\hat{\h}^{\yng(1)}_{\ga_1,\ga_2(m)}$ be
the corresponding superpolynomial for the cable
$C\!ab(2\rr_1\ss_1\!+\!2m\!+\!1,2)
T(\rr_1,\ss_1)$.
Then
\begin{align}\label{fga1-m}
\hat{\h}^{\yng(1)}_{\ga_1,\ga_2(m)}
=&\frac{1\!+\!aq}{1\!-\!qt}
\left(1\!+\!(qt)^m\frac{q(1\!-\!t)}{1-q}\right)
\hat{\h}_{\rr_1,\ss_1}^{\yng(2)}\!-\!
(qt)^m\frac{q}{1\!-\!q}\hat{\h}_{2\rr_1,2\ss_1}^{\yng(1)}\notag\\
&=\!\frac{1\!-\!(qt)^m}{1\!-\!qt}(1\!+\!aq)
\hat{\h}_{\rr_1,\ss_1}^{\yng(2)}+
(qt)^m\hat{\h}_{\ga_1,\ga_2(0)}^{\yng(1)},
\end{align}
where $\hat{\h}_{\rr,\ss}^{\yng(2)}$,
$\hat{\h}_{2\rr,2\ss}^{\yng(1)}$ are as
in the previous proposition (which is used to obtain
the second formula from the first). The formula
for $\ga_2(m)$ in terms of $\ga_2(0)$  
actually works for any $m\in \Z$. \sq
\end{proposition}
In the case of $\ga_1=\ga_{3,2}$ these formulas are 
an important addition to the (mostly numerical)
analysis of this family in \cite{ChD1}. For instance, for $m=-5$
we obtain the superpolynomial $(\mathbf c)$ from
Section 4.1 there:
$$
({\bf c})\ \vec\rr=\{3,2\},\,
\vec\ss\!=\!\{2,-9\}\for 
C\!ab(3,2)T(3,2);\, 
\h_{\vec\rr,\vec\ss}(\square;q,t,a)=
$$
\renewcommand{\baselinestretch}{0.5} 
{\small
\(
1-q^2+q t+q^2 t
-q^3 t+q^2 t^2+q^3 t^3+a^3 \left(-\frac{q^4}{t^2}
-\frac{q^5}{t}\right)+
a^2 \bigl(q^3-q^4-q^5-\frac{q^3}{t^2}-\frac{q^3}{t}
-\frac{2 q^4}{t}\bigr)
+a \left(q+q^2-2 q^3-q^4-\frac{q^2}{t}
-\frac{q^3}{t}+q^2 t+q^3 t-q^4 t+q^3 t^2\right).
\)
}
\renewcommand{\baselinestretch}{1.2} 
\smallskip

Negative $m$ result in non-algebraic knots, however the
formula is uniform for any $m$. Negative
powers of $q,t$ cancel each other for $m<0$. 
A geometric 
interpretation
of non-algebraic torus iterated knots/links is an interesting
problem. See Section \ref{SEC:469} below. 
We note that Proposition \ref{FAM1-m} and its generalizations
provide important tools for counting points of Jacobian factors 
over finite fields in {\em families\,}, which is
parallel to Iwasawa theory in number theory.
See the end of Section \ref{SEC:COMM}.

Conjecture \ref{CONJ:HLF} is expected to follow from such
recurrence relations (in full generality), to be considered
elsewhere. We also hope that these relations will 
connect DAHA superpolynomials with 
the stable $KhR$\~polynomials 
(via the approach  based on Soergel modules).   

\Yboxdim7pt

There are no actual
problems with applying this formula to any $\ga_1$;
practically, its complexity is generally comparable with 
that of the starting $\hat{\h}_{\ga_1,\ga_2(0)}$.
The simplest example is for $\ga_1=\ga_{3,2}$ and 
the range $m\ge 0$.
Here the genus is $g(m)=8+m$ and we have:

\begin{align}\label{fam3-2}
&\left(\frac{t}{u}\right)^{8}
\hat{\mathbb H}_{\{3,2\},\{2,3+2\Z_+\}}=
\sum_{m=0}^\infty\left(\frac{u}{t}\right)^{m}
\hat{\h}_{\{3,2\},\{2,3+2m\}},\\
&(1-u/t)(1-u q)\,\left(\frac{t}{u}\right)^{8}
\hat{\mathbb H}_{\{3,2\},\{2,3+2\Z_+\}}=\notag
\end{align}
\renewcommand{\baselinestretch}{0.5} 
{\small
\(
1+q t+q^2 t+q^3 t+q^2 t^2+q^3 t^2+2 q^4 t^2+q^3 t^3+q^4 t^3
+2 q^5 t^3+q^4 t^4+q^5 t^4+2 q^6 t^4+q^5 t^5+q^6 t^5+q^7 t^5
+q^6 t^6+q^7 t^6+q^7 t^7+q^8 t^8+\bigl(-q-q^2 t-q^3 t-q^4 t
-q^3 t^2-q^4 t^2-2 q^5 t^2-q^4 t^3-q^5 t^3-2 q^6 t^3-q^5 t^4
-q^6 t^4-q^7 t^4-q^6 t^5-q^7 t^5-q^7 t^6-q^8 t^7\bigr) \mathbf u
+\mathbf{a}^3 \Bigl(q^6+q^7 t+q^8 t^2+\bigl(-q^7-q^8 t\bigr) 
\mathbf u\Bigr)+\mathbf{a}^2 \Bigl(q^3+q^4+q^5+q^4 t+2 q^5 t
+2 q^6 t+q^5 t^2+2 q^6 t^2+2 q^7 t^2+q^6 t^3+2 q^7 t^3+q^8 t^3
+q^7 t^4+q^8 t^4+q^8 t^5+\bigl(-q^4-q^5-q^6-q^5 t-2 q^6 t-2 q^7 t
-q^6 t^2-2 q^7 t^2-q^8 t^2-q^7 t^3-q^8 t^3-q^8 t^4\bigr) 
\mathbf u\Bigr)
+\mathbf{a} \Bigl(q+q^2+q^3+q^2 t+2 q^3 t+3 q^4 t+q^5 t+q^3 t^2
+2 q^4 t^2+4 q^5 t^2+q^6 t^2+q^4 t^3+2 q^5 t^3+4 q^6 t^3
+q^7 t^3+q^5 t^4+2 q^6 t^4+3 q^7 t^4+q^6 t^5+2 q^7 t^5
+q^8 t^5+q^7 t^6+q^8 t^6+q^8 t^7+\bigl(-q^2-q^3-q^4-q^3 t
-2 q^4 t-3 q^5 t-q^6 t-q^4 t^2-2 q^5 t^2-4 q^6 t^2-q^7 t^2
-q^5 t^3-2 q^6 t^3-3 q^7 t^3-q^6 t^4-2 q^7 t^4-q^8 t^4
-q^7 t^5-q^8 t^5-q^8 t^6\bigr) \mathbf u\Bigr).
\)
}
\renewcommand{\baselinestretch}{1.2} 
\smallskip

We note that when $a=0$, the coefficient of $u^1$
considered upon the
substitution $q\mapsto qt$ satisfies $RH$ (its $t$\~roots $\ze$
are complex such that $|\xi|=q^{1/2}$)
provided that $q<0.7562688464467736\ldots$. This is better
than the bound $q\le 1/2$ in Conjecture \ref{CON:om2} for
(individual) $m$.

\subsection{\bf Motivic superpolynomials}\label{SEC:MOTZ}
In the unibranch uncolored case, let us compare
the {\em motivic superpolynomials\,} from \cite{ChP1}
and the {\em Galkin-St\"ohr zeta-functions\,}.

\subsubsection{\sf Standard modules}
Let $\r\subset\mathbb \o\equal \mathbb F[[z]]$ be 
a Gorenstein ring over a finite field 
$\mathbb F\!=\!\mathbb F_{q}$ of cardinality $q$,  
$\de=dim_{\mathbb F}(\o/\r)$ (the arithmetic genus),
$\Ga_\r=val_z(\r)$ for the usual $z$\~valuation.
It is a semigroup and  $\de =|\Z_+\setminus \Ga|$.
For the later, we assume that $\r=\tilde{\r}_{\mathbb F_q}\equal
\tilde{\r}\otimes_{\Z}\mathbb F$, where 
$\tilde{\r}\subset \Z[[z]]$, and $\Ga_{\C}$ for
$\tilde{\r}\otimes_\Z \C$ coincides with $\Ga$ over $\mathbb F$. 

Given an $\r$\~module $M\subset \o$,  
$\De=\De_M\equal val_z(M)$ is $\Ga$\~module, i.e. 
$\Ga+\De=\De$. One has:
$
dev(M)\equal\de-dim_{\mathbb F}(\o/M)=\de-
|\Z_+\setminus\De|\equal dev(\De). 
$ If $0\in \De$, then $\De$ and the
corresponding $M$ are called {\em standard\,};
equivalently,  $M\cdot \o=\o$.
For any $\r$\~module $M$, let $M_{st}\equal z^{-m} M$
for $m=\min{\De_M}$, which is a standard module corresponding
to $\De_{st}=\De-\min\{\De\}$. Also, 
$M^\star\equal\{x\!\in\! 
\mathbb F((z))\mid xM\!\subset\! \r\}$ corresponds to
$\De^\star\equal\{n\in \Z_+
\,\mid\,n+\De\subset \Ga\}$ for any modules $M,\De$.

\begin{definition}\label{DEF:MOT}
For any $\mathbb F$\~subalgebra $\r\subset \o$ with 
$\mathbb F((z))$ as its field of fractions, let
$J\!=\!J_{\r}(\mathbb F)
\!\equal\!\{M=M_{st}\}$ be the {\sf Jacobian factor\,},
$J_\De\!\equal\!\{M=M_{st}, \De(M)=\De\},$ 
$\h_{mot}^0(q,t)=\sum_{M\in J}t^{dim_{\mathbb F}(\o/M)}$. For
a standard $\De$, we set
$\h_{mot}^0(\De,q,t)=\sum_{M\in J_{\!\De}}
t^{dim_{\mathbb F}(\o/M)}$. 
Following \cite{Sto} for Gorenstein $\r$, $M\subset \o$, 
 and standard $\De$:
$$
Z(\De,q,t)=
\sum_{M\supset \r}^{M_{st}\in J_\De } t^{dim_{\mathbb F}(M/\r)}=
\sum_{M\subset \r}^{(M^\star)_{st}\in J_\De} 
t^{dim_{\mathbb F}(\r /M)}.
$$
Accordingly, $Z(q,t)=\sum_{\De=\De_{st}}Z(\De,q,t).$
Also we set $$L(\De,q,t)\equal (1\!-\!t)\,Z(\De,q,t),\ 
L(q,t)\equal (1\!-\!t)\,Z(q,t).$$
\vskip -0.7cm\sq
\end{definition}

\subsubsection{\sf St\"ohr's formula}
We begin with the {\em functional equation\,}:
\begin{align}\label{feqze}
(qt^2)^\de\,L(\De,q,\frac{1}{qt})= L((\De^\star)_{st},q,t)
\hbox{\, for standard  } \De.
\end{align}
The natural setting
here is when the summation over $M$ from Definition \ref{DEF:MOT}
is reduced to (any) $\o^\ast$\~orbits;
$\o^\ast=\{x\in \o\mid x^{-1}\in \o\}$ acts in $J_\De$
by multiplication. 
See \cite{Sto},(3.10). We are going to state
a version of Theorem 3.1 there, which almost directly provides
(\ref{feqze}). Let
$$
g_\De(n)=\big| \{\N\ni m\not\in \De\,\mid\, m<n\} \big| 
\for n\in \Z_+,\hbox{\ so\ } g_\De(0)\!=\!0\!=\!g_\De(1).
$$
Using the relations $dev(M\cap (z^n\o))=dev(M)\!-\!n\!+\!g_\De(n)$,
one has:
\begin{align}\label{LDeqt}
Z(\De,q,t)=(qt)&^{dev(\De)}q^{-\de}\,\big|J_\De\big|
\sum_{n\in \De}t^n\,q^{g_{\!\De}(n)}\for \De=\De_{st},\notag\\
L(\De,q,t)=(qt)&^{dev(\De)}q^{-\de}\,\big|J_\De\big|
\notag\\
&\times\bigl(\sum_{n\in \De\not\ni n-1} t^n\,q^{g_{\!\De}(n)}
-\sum_{n\not\in\De\ni n-1}t^n\,q^{g_{\!\De}(n-1)}\bigr).
\end{align} 
Note that the second formula readily gives that 
$L(\De,q,t=1)=\big|J_\De\big|$.

\begin{conjecture}\label{CONJ:HL}
Let ${\r}_{\C}\subset \C[[z]]$ be the ring of
plane curve singularity (with $2$ generators and the same 
field of fractions). Within
its topological type, there exists ${\r}_{\Z}\subset 
\Z[[z]]$ such that ${\r}_{\C}$ and 
$\r_{\mathbb F}$ have coinciding $\Ga$. Let 
$\hat{\h}(q,t,a)$ be
the corresponding uncolored DAHA superpolynomial
for the link of $\r_{\C}$. Then 
$\h^0_{mot}(q,t)=\hat{\h}^0\equal\hat{\h}(q,t,a\!=\!0)$, 
which is the case $a=0$ of 
Conjecture 2.3,(iii) from \cite{ChP1}. Then one has:
\begin{align}\label{conzeconj}
&\h^0_{mot}(q\mapsto qt,t)\,=\, L(q,t)\,=\,
\hat{\h}(q\mapsto qt,t,a\!=\!0),\\
&L(\Ga,q,t)=
\hat{\h}(q,t,a\mapsto -(t/q))\big|_{q\mapsto qt}\,,\
\big|J_\Ga\big|=q^\de.\label{czej}
\end{align}\vskip -0.7cm\sq
\end{conjecture} 

This conjecture (and its $a$\~generalization,
see below), clarifies the substitution
$q\mapsto qt$; the DAHA super-duality then becomes the functional
equation. The DAHA super-duality  
(with $a$ and for any colors) is not very difficult
to justify, but some theory is needed; see \cite{GoN}
and \cite{CJJ} for a sketch of the proof via roots of unity.
(which can work beyond $A_n$).
The positivity of $L(q,t)$ seems new for the Galkin-St\"ohr zeta
(beyond some special values/ coefficients). The corresponding
cancelations in (\ref{LDeqt}) generally  
hold only for plane curve singularities. 


\begin{corollary}\label{SHUFF}
Provided (\ref{conzeconj}),
$\hat{\h}(t,t,a\!=\!0)\!=\!\sum_{\De=\De_{st}} 
t^{dev(\De)}\kapp(\De)$, where
$\kapp(\De)$ is the image of $J_\De$ considered as an
abstract (projective) variety in the quotient
of the Grothendieck ring $K_0(Var)$ by $1-\mathbb L$
for $\mathbb L=[\mathbb A^1]$,
e.g. $\kapp(\De)=1$ if  $J_\De$ is an 
affine space $\mathbb A^{\!N}$.\sq
\end{corollary}

\subsubsection{\sf A non-planar example}\label{SEC:469}
Relation (\ref{conzeconj}) was checked numerically for quite a 
few examples of plane curve singularities,
 including many torus knots, 
$\r=\mathbb F[[z^4,z^6+z^{7+2m}]]$ for many $m\ge 0$, 
$\mathbb F[[z^6,z^8+z^{9,\ldots,15}]]$, and
$\r=\mathbb F[[z^6,z^9+z^{10}]]$.

Interestingly, it fails for {\em non-planar\,}
Gorenstein singularities. Let $\r=\mathbb F[[z^4,z^6,z^9]]$.
Here $\de=6$ and it corresponds in some sense (which
we omit here) to the non-algebraic
knot $C\!ab(9,2)T(3,2)$, called {\em pseudo-algebraic\,} 
in \cite{ChD1}. See also Corollary 1.4 from \cite{Hed}. 
The positivity of the DAHA superpolynomial and its algebraic
similarity to that for $C\!ab(13,2)T(3,2)$ were the defining
features in \cite{ChD1}. One has:  $\hat{\h}=$
{\small
\(
1 + a^3 q^6 + q t + q^2 t + q^3 t + q^2 t^2 + q^3 t^2 + 
2 q^4 t^2 + q^3 t^3 +  q^4 t^3 + q^5 t^3 + q^4 t^4 + q^5 t^4 + 
q^5 t^5 + q^6 t^6 + a^2 (q^3 + q^4 + q^5 + q^4 t + 2 q^5 t + 
q^6 t + q^5 t^2 + 
q^6 t^2 +    q^6 t^3) + a (q + q^2 + q^3 + q^2 t + 2 q^3 t + 
3 q^4 t + q^5 t + q^3 t^2 + 2 q^4 t^2 + 3 q^5 t^2 + q^4 t^3 + 
2 q^5 t^3 + q^6 t^3 +    q^5 t^4 + q^6 t^4 + q^6 t^5).
\)
}
The smallest positive $\hat{\h}$ is for $C\!ab(7,2)T(3,2)$
associated with $\mathbb F[[z^4,z^6,z^7]]$, but its
$a$\~degree drops to $2$.

The motivic interpretation of $\hat{\h}$ 
above at $a=0$ is as follows.
We define 
 $'\h^0_{mot},\,'\!L$ by restricting 
the summation in $\h^0_{mot},Z$ to $\De$ such that  
$\De_{1,2}\setminus\Ga \neq \{2,11\},\{2,7,11\}$. 
We note that $\De_{1,2}$ satisfy  $(\De^\star)_{st}=\De$. 
Then
$'\h^0_{mot}(q\!\mapsto\! qt,t)\!=\! '\!L(q,t)$
(and with $a$, see (\ref{motpolfull}) below),
but $\h^0_{mot}(q\!\mapsto\! qt)\neq L$. 
We note that (\ref{czej}) holds without prime.

Namely, $\h^0_{mot}\!-'\!\!\h^0_{mot}=q^6t^4+q^5t^3$;
the latter monomials are contributions of $\De_1,\De_2$. 
One has:
$L(q,t)-\,'\!L(q,t)=q^2t^2-q^2t^4+q^3t^4-q^3t^6+q^4t^6
-q^4t^8+q^5t^8+q^6t^{10}$; by the way, $L(q,t)$ has positive
coefficients. Here $L(q,t))/t^6$ and $'\!L(q,t)/t^6$ 
satisfy the functional equation $t\mapsto 1/(qt)$. However, 
$\h^0_{mot}/t^6$ is not self-super-dual
for $q\leftrightarrow t^{-1}$. 
This failure is typical for motivic 
pseudo-algebraic superpolynomials. 

\vskip 0.2cm

{\sf Observation.}
For planar $\r$, let $\ddot{\r}$ be the
corresponding {\em quasi-homoge\-neous\,} ring, which is
$\ddot{\r}\equal\mathbb F[z^v, v\in \Ga]$, and 
$\ddot{\h}_{mot}(q,t,a)$ as in (\ref{motpolfull}).
The following holds for 
$\mathbb F [[z^4,z^6\!+\!z^{7+2m}]]$ we checked and for 
$\r\!=\!\mathbb F [[z^6,z^8\!+\!z^9]]$: 
{\em $\hat{\h}$ is $\ddot{\h}^{sym}_{mot}$, which is
$\ddot{\h}_{mot}$ minus
the sum 
of positive terms in $\ddot{\h}_{mot}(q,t,a)\!-\!
(qt)^\de\ddot{\h}_{mot}(\frac{1}{t},\frac{1}{q},a)$}.
The same tendency holds for (non-planar) $\r$ 
associated with {\em pseudo-algebraic\,} knots.

The smallest counterexample we found at $a\!=\!0$ is
$\mathbb F [[z^6,z^9\!+\!z^{10}]]$ ($\de\!=\!21$),
where $\ddot{\h}^{sym}_{mot}\!-\!\hat{\h}\big|_{a=0}=$ {\small
$q^{13}t^5+q^{14}t^5+2q^{14}t^6+2q^{15}t^6+2q^{15}t^7
+q^{16}t^7+q^{16}t^8$}. 
The number of $\De$ is $447$
for the corresponding $\ddot{\r}$.
By the way,  $J_{\!\ddot{R}}$ has a
component of $dim\!=\!\de\!+\!2$ and there exist $9$ non-affine
cells $\ddot J_\De$, 
but they do not contribute to
$\ddot{\h}^{sym}_{mot}\!-\!\hat{\h}\big|_{a\!=\!0}$
(all $J_\De$
are affine).

The rationale here is that $i\!-\!j\!=\!dim(J_\De)\!+\!dev(\De)
\!-\!\de$  is high
for $\De$ (with affine $J_\De$) 
and the corresponding $q^i t^j$ in $\hat{\h}$ if
$dim(\ddot J_\De)\!>\!dim(J_\De)$ 
($dim\,\emptyset\!=\!-\infty$).  
In the difference above, $i\!-\!j=8,9$, when 
$max\{i\!-\!j\}\!=\!9$ in $\hat{\h}^0$,
and only very few
$\De$ there are with $i\!-\!j\!\ge\! 8$.  

For instance, let $\r=\mathbb F[[z^4,z^6,z^{11}]]$.
Then the reduction above and the relations
$'\h_{mot}(q\!\mapsto\! qt,t,a)\!=\! '\!\l(q,t,a)\!=\!
\hat{\h}(q\!\mapsto\! qt,t,a)$ hold for non-admissible 
$\De\setminus\Ga=
\{2,13\},\{2,9,13\}$. For $z^7$ instead
of $z^{9,11}$ (not pseudo-algebraic according to 
\cite{ChD1}): $\hat{\h}\!=\!'\h_{mot}$
if  $\{2,9\},\{1,5,9\}$ are excluded,
but then
$\hat{\h}-\,'\!\l\!=\!q^2t^4(1\!-\!t)(1\!-\!qt)(1\!+\!aq\!+\!aqt)$,
which requires further adjustments. Also,
the observation above holds only for $a=0$. We note that
$\h^0_{mot}\!-\!\hat{\h}(a\!=\!0)$
for $11,9,7,5,3$ is uniform:\ 
$(qt)^\de(t^{-2}\!+\!t^{-3}q^{-1})$, where $\de\!=\!7,6,5,4,3$.


\subsubsection{\sf Flagged zeta-functions}\label{SEC:FLAGGED}
Following \cite{ChP1}, the {\em standard $\ell$\~flag\,} 
of $\Ga$\~modules is
a sequence
$\vec{\De}=
\{\De_0\!\subset\! \De_1\!=\!\De_0\cup\{g_1\}\subset
\!\cdots\!\subset
\De_i\!=\!\De_{i-1}\cup \{g_i\}\subset\! \De_\ell\}
$ of standard $\Ga$\~modules $\De_i$
such that $0\neq g_i\in \Z_+\setminus
\De_{i-1}$ and $g_{i-1}<g_{i}$ for $1\le \ell$.
Thus   $dev(D_i)=dev(D_{i-1})+1$.
We set $dev(\vec{\De})\equal dev(\De_m)$ and define the
{\em standardizable\,} flag of $\r$\~modules of 
length $\ell$ as the sequences
of $\r$\~modules 
$\vect{M}=\{M_0\!\subset\! M_1\!\cdots\!\subset  M_\ell\}$ 
such that $\De(\vect{M})\equal \{\De(M_i)\}$ becomes a standard
flag as above upon the subtraction of  $m=min(\De_{\ell})$,
i.e. for  $\vect{M}_{st}\equal z^{-m}\vect{M}$. Accordingly,

\begin{align}\label{zaqt}
\z(q,t,a)\equal
\sum_{\vect{M}\subset \r} 
a^\ell\, t^{dim_{\mathbb F}(\r /M_{\ell})},\ \,
\l(q,t,a)\equal (1\!-\!t)\,\z(q,t,a),
\end{align}
where the summation is over {\em standardizable\,}
flags. The (full) {\em motivic superpolynomial\,}
from \cite{ChP1} is as follows:
\begin{align} \label{motpolfull}
\h_{mot}(q,t,a)=
\sum_{\vect{M}\in \vec{J}}a^\ell\,
t^{dim_{\mathbb F}(\o/M_{\ell})},
\end{align}
where $\vec{J}=\vec{J}_{\!\!\r}$ is a scheme of all 
{\em standard\,} flags of
submodules in $\o$. 

\begin{conjecture}\label{CONJ:HLF}
In the case of plane
curve singularity, 
\begin{align}\notag
&\hat{\h}(q,t,a)\!=\!\h_{mot}(q,t,a)\!=\!\l(\frac{q}{t},t,a),
\,\, \h_{mot}(qt,t,a\!=\!-\frac{1}{q})\!=\! L(\Ga,q,t),
\end{align}
where the latter possibly holds for
any Gorenstein rings $\r$. \sq
\end{conjecture} 

It was checked in quite a few cases; as for   
$\h_{mot}(q,t,a)\!=\!\l(\frac{q}{t},t,a)$, including
$\mathbb F[[z^6,z^8\!+\!z^9]]$ for $\ell=0,1$ and 
$\mathbb F[[z^6,z^9\!+\!z^{10}]]$ for $\ell=0$. 
We will omit here flag generalizations of
formulas from (\ref{LDeqt}), which result in
a relatively straightforward proof of the flagged functional
equation for $t^{-\de}\l(q,t,a)$
upon the (same) transformation $t\mapsto 1/(qt)$. 
The generalization to algebraic (uncolored) links is
also known, as well as some steps of the
justification based on the analysis
of the monoidal transformations of singularities.
 (at least
for some families). 


\subsubsection{\sf Nested Hilbert schemes}
\label{SEC:HILB}
Our flags of $\r$\~modules can be interpreted via 
{\em restricted\,}
nested Hilbert schemes $Hilb^{\,m,m-\ell}_{\,res}$,
and nested Jacobian factors $J^{m,m-\ell}_{res}$:
\begin{align}
&Hilb^{\,m,m+\ell}_{\,\,res}=
\{ M\supset M'\supset \mathfrak{m} M\mid
M\subset \r,\\ 
\hbox{dim}(\r/&M)\!=\!m,\,  \ \hbox{dim}(M/M')\!=\!\ell,\ \, 
M\otimes_{\r} \o\!=\!
M'\otimes_{\r} \o \},\notag\\
&J^{\,m,m+\ell}_{\,\,res}=
\{ M\supset M'\supset \mathfrak{m} M\mid
M\subset \o,\\ 
\hbox{dim}&(\o/M)\!=\!m,\,  \ \hbox{dim}(M/M')\!=\!\ell,\ \, 
M'=M'_{st}\}.\notag
\end{align}
Here $\mathfrak{m}$ is the maximal ideal of $\r$.
If the condition $M\otimes_{\r} \o\!=\!
M'\otimes_{\r} \o$ or $M'=M'_{st}$ is omitted, one 
obtains nested $Hilb$ from \cite{ORS} or its variant
for $J$. As above, the pairs satisfying
this condition are called {\em standardizable\,};
also, $M'\!=\!M'_{st}$ obviously implies
$M\!=\!M_{st}$.

Using Proposition 2.3 from
\cite{ChP1}, 
we obtain that $\mathfrak{m}M_{\ell}\subset M_0$
for any standardizable flag $\vect{M}$. Moreover,
given any standardizable pair
$\{M\!=\!M_{\ell}\!\supset\!M'\!=\!M_0\}$, the
number of the corresponding 
standardizable $\ell$\~flags $\vect{M}$ is $q^{\ell(\ell-1)/2}$.
Furthermore, we have the following proposition.
\begin{proposition}\label{PROP:DIMFL}
For an $\r$\~module $M$, let
$
\De\!=\!\De(M)\!=\!\De(\mathfrak{m}M)\cup$ 
$\{d_1\!\!<\!d_2\!\!<,\ldots,<\!d_{r}\}$,
where $d_1\!=\!\min\{\De(M)\}$ and 
$ r\!=\!dim_{\r/\mathfrak{m}}
(M/\mathfrak{m}M)$. 
And let
$\De'\!=\!\De(M)\setminus
\{g_1\!\!<\!g_2\!\!<,\ldots,<\!\!g_{\ell}\}$
for $g_j$ taken from the set $\{d_i,\,i\!=\!2,\ldots,r\}$;
it is a $\Ga$\~module (containing $d_1$).
Setting 
$\{d_i\}\!\setminus\!\{g_i\}\!=\!
\{d_1\!=\!g^\circ_1\!\!<\!g^\circ_2\!\!<,\ldots,<
\!\!g^\circ_{r-\ell}\}$, 
the
number of standardizable $\ell$\~flags $\vect{M}=
\{M_i\}$ with $M_\ell=M$ and 
$\De(M_0)\!=\!\De'$ equals $q^N$ for $N=N(M,\De')=$
$$
\ \, |\{g_i\!>\!g^\circ_1\}|\!+\!|\{g_i\!>\!g^\circ_2\}|\!+\!\ldots
\!+\!
|\{g_i\!>\!g^\circ_{r-\ell}\}|\!+\! \frac{\ell(\ell\!-\!1)}{2},
\ \, 1\!\le\! i\!\le\! \ell\!<\!r.  \hbox{\hfill \sq}
$$
\end{proposition}

Geometrically, the flags with fixed
$M_\ell$ and $\De_0=\De(M_0)$ form an affine space 
$\mathbb A^{N(M_\ell,\De_0)}$.
We use here the Nakayama Lemma; cf. Section 2.1 from \cite{ORS}
and Section 9.1 from \cite{GORS}.
Thus, assuming that the stratification of 
$J_{\De_\ell}$ with respect to the Nakayama rank $r(M_\ell)$ 
is known, the calculation of $\h_{mot}(q,t,a)$ and
$\l(q,t,a)$ becomes in terms of
$\De_l,\De_0$. Using that  $r(M_\ell)=const$\,
within $J_{\Ga}$\~orbits (they are affine spaces), we
obtain the functional equation for 
$\l(q,t,a)$ (following St\"ohr). This proposition coupled
with our conjectures
provides far-reaching generalizations of the so-called
Shuffle Conjecture; see \cite{CaM}.  

Let us relax the definition of standardizable flags by
allowing $g_1$ to be $m\!=\!\min(D_0)$. Such flags can be
called (partial) full {\em 
gap-increasing\,} due to $g_{i-1}\!<\!g_i$.
We actually add to the standardizable $\ell$\~flags the
standardizable $(\ell-\!1)$\~flags from $M_0$ extended by
$M_{-1}\!=\!M_0\cap z^{m+1}\o$. This gives the following connection
with the usual nested Hilbert schemes:
\begin{align}\label{zaqtt}
&(1+a)\z(q,t,a)\,=\,\sum_{m,\ell=0}^\infty 
q^{(\ell-1)\ell/2}\,a^\ell t^m\, 
\bigl|Hilb^{m+\ell,m}(\mathbb F_q)\bigr|,
\end{align}
\vskip-0.2cm
\noindent
where $|\cdot|$ is the number of points, which is
directly related to the ORS conjecture. Namely, 
one replaces $\mathfrak{w}$ in (\ref{overlinep})
by the count of $\mathbb F_{q}$\~points (we will not comment on
that) and substitutes
$q^2_{st}\!\mapsto\! t,\, a_{st}^2 t_{st}\!\mapsto\! a,\,
t^2_{st}\!\mapsto\! q$; recall that $q$ in this section
 is $q/t$ 
via the DAHA parameters $q,t$.
\vskip 0.3cm

{\sf Using Proposition \ref{PROP:DIMFL}.} 
In fact, we gave a different definition
of $\h_{mot}$ in \cite{ChP1}:
$\h^\checkmark_{mot}\equal
t^{\de}\sum_{\vect{M}\in \vect{J}}\,a^\ell\,
q^{dev(M_{\ell})}$ for the field $\mathbb F_{1/t}$.
Accordingly, if $J_{\!\vec{\,\De}}$ is $\mathbb A^N$,
its contribution to $\h^\checkmark_{mot}$ is 
$a^\ell q^{dev(M_{\ell})} t^{\de-N}$. 
It coincides with $\h_{mot}$ assuming its self-duality, but 
can be more convenient.
 
Let $\r\!=\!\mathbb F[[z^4,z^6\!+\!z^7]]$. For 
$\,\ell\!=\!3\!=\!deg_a(\hat{\h})$, we use
 \cite{ChP1},(4.1):
\begin{align}\label{13-3rh} 
&D_0=[9,11,15],\, \vec{g}= (2,5,7),\, \hbox{dim}=8 
&\rightsquigarrow q^6t^0a^3,\\
&D_0=[7,9,11,15],\, \vec{g}=(2,3,5),\, \hbox{dim}=7
&\rightsquigarrow q^7t^1a^3,\notag\\
&D_0=[5,7,9,11,15],\, \vec{g}=(1,2,3),\, \hbox{dim}=6
&\rightsquigarrow q^8t^2a^3,\notag
\end{align}
where $D_i\!=\!\De_i\!\setminus\! \Ga$,\  $\vec{g}$
are the gaps (consecutively) added to $D_0$, and 
we show the contributions of the corresponding
cells to $\h^\checkmark_{mot}$. One has:
$$
D_3=[2,5,7,9,11,15],\ [2,3,5,7,9,11,15],\ 
[1,2,3,5,7,9,11,15].
$$
The corresponding $M_3/\mathfrak{m}M_3$ are all of rank
$r=4$ for any $M_3$ for these $D_3$ with an important
reservation. In the case of $D_3=[2,5,7,9,11,15]$, the
rank is $r=3$ for generic $M_3$, and an affine subspace
of codimension $1$ in $J_{\De_3}$ must be taken to
ensure $r=4$; otherwise 
$\De[\mathfrak{m}M_3]=[7,9,11,15]$. 
Thus the proposition gives that the dimensions
in (\ref{13-3rh}) are (indeed)
$\ell(\ell-1)/2+\ell=6$ plus
dim$J_{\De_3}=3,1,0$ minus $1$ for the first $D_3$.
See Table 1,\cite{ChP1}; all cells $J_\De$ are affine for 
$\r\!=\!\mathbb F[[z^4,z^6\!+\!z^7]]$.
\vskip 0.2cm

{\sf The case  $\ell\!=\!2$, $D_2$\!=\![2,5,7,9,11,15].}
There are now 3 possibilities for $D_0$ here:
$D_0=[2,9,11,15],D_0'=[7,9,11,15],D_0''=[5,9,11,15]$;
see now Table 3. The dimension of $J_{\De_2}$, which
is $3$,  must be
diminished by $1$ for $D_0$ and $D_0''$ due to the
absence of $7$ there (similar to the example above;
the same subspace serves $D_0$ and $D_0''$).
Then $r=4$ and $\{d_i\}=\{0,2,5,7\}.$ The summation 
in Proposition \ref{PROP:DIMFL} becomes:
$\frac{\ell(\ell\!-\!1)}{2}+\{4,2,3\}=\{7,6,6\}$
for $\{D_0,D_0',D_0''\}$.
The corresponding contribution to $\h_{mot}^\checkmark$
is $2q^6 t^2+q^6 t$. This matches Table 3 in \cite{ChP1}.
\vskip 0.2cm

{\sf Z(q,t,a) for trefoil.}
This proposition can be equally used for 
ideals in $\r$, though for $\r=\mathbb F[[z^2,z^3]]$
this is simple. The corresponding standardizable
ideals $M=\lan\cdot\ran_d\subset \r$
with $d=$dim$\r/M$ are as follows:
{\small 
$$\lan 1,z^2,z^3\ran_0,\ \lan z^2,z^3\ran_1,\ 
\lan z^2+\la z^3\ran_2,\, \lan z^3,z^4\ran_2,\ 
\lan z^3+\la z^4\ran_3,\,\lan z^4, z^5\ran_3,
\ldots,
$$
}
\!where $\la\in \mathbb F$. The standardizable pairs of 
$\ell\!=\!1$
are $\lan z^{i},z^{i+1}\ran_{i-1}\supset
\lan z^i+\la z^{i+1}\ran_i$ for $i\ge 2$. Thus
$\z=(1+qt^2+a q t )/(1-t)$.

The pair $\r\supset \mathfrak{m}$ is non-standardizable
of $\ell\!=\!1$; the other such pairs are
$\lan z^{i},z^{i+1}\ran_{i-1}\supset
\lan z^{i+1},z^{i+2}\ran_i$ and 
$\lan z^{i}+\la z^{i+1}\ran_{i}\supset
\lan z^{i+2},z^{i+3}\ran_{i+1}$. Also, there are
pairs $\lan z^{i},z^{i+1}\ran_{i-1}\supset
\lan z^{i+2},z^{i+3}\ran_{i+1}$ with $\ell\!=\!2$.
Formula (\ref{zaqtt}) gives then
 $(1\!+\!a)(1\!+\!qt^2\!+\!a q t )/(1\!-\!t)$,
which matches (\ref{overlinep}).

\comment{
Replacing $\mathfrak{w}$ in (\ref{overlinep})
by the count of $\mathbb F_{q}$\~points (we will not comment on
that), $q^2_{st}\!\mapsto\! t,\, a_{st}^2 t_{st}\!\mapsto\! a,\,
t^2_{st}\!\mapsto\! q$ (recall that $q$ now is $\frac{q}{t}$ 
via DAHA $q,t$), the ORS conjecture coupled with our
ones gives that the sum of 
$q^{\frac{\ell(\ell-1)}{2}}t^{d_{\ell}}a^\ell$ 
over all pairs $\{M_0\!\subset\! M_{\ell}\}$ above (not only 
standardizable) is $\l(q,t,a)\frac{1+a}{1-t}$, which holds. 
Generally, $\frac{\ell(\ell-1)}{2}$ is related
to that in the formula for $N$ in
Proposition \ref{PROP:DIMFL}, but there are other terms there 
(and it is for standardizable modules); we hope to address this
elsewhere."
}
\comment{
The corresponding sum is 
$\frac{a(1+qt^2)+a^2t}{1-t}$, and the sum for the
whole (non-restricted) nested $Hilb\,$
is not $(1\!+\!a)(1\!+\!qt^2\!+\!a q t )/(1\!-\!t)$, as
one would 
expect due to conjectural relations  the coefficient
of $a^2$ must be $qt$ for it
instead of $t$ we obtained above.
There is no discrepancy here in 
the HOMFLY-PT case proved in
\cite{Ma}, which is when $q=1$. 
}

\subsection{\bf On connection conjectures}
\subsubsection{\sf HOMFLY-PT polynomials}
Given a link colored by a set of Young diagrams,
let $\hbox{\small H\!O\!M}(q_{st},a_{st})$ be
the corresponding unreduced HOMFLY-PT polynomial. 
They can be defined via Quantum
Groups (in type $A$) or
using the corresponding {\em skein relations\,} and
Hecke algebras.  See e.g. \cite{QS} and references
there; we provide here only a sketchy discussion.

Recall that iterated torus links are determined by
the pairs of graphs $\{\l,\,'\!\l\}$ colored by  arbitrary
sequences $\{\la^j\},\{'\!\la^j\}$ of Young diagrams.
We need to switch to the {\em reduced\,} HOMFLY-PT polynomials
with respect to one of its components, say $j_\circ$,
and then perform the hat-normalization; the notation will be
$\hat{\hbox{\small H\!O\!M}}^{j_o}_{\l,\,'\!\l}\,
(q_{st}, a_{st})$. 
Let us mention that the $q$\~polynomiality of the 
{\em unreduced\,} $\hbox{\small H\!O\!M}(q_{st},a_{st})$  
generally does 
{\em not\,} hold for links. 

We put $\,\hat{\h}^{j_\circ}_{\l,\,'\!\l}\,(q,t,a)_{st}\,$
for \,$\hat{\h}^{j_\circ}_{\l,\,'\!\l}\,(q,t,a)\,$ from
Theorem \ref{STABILIZ} expressed
in terms of the {\em standard topological parameters\,}
(see \cite{CJ}, Section 1 in \cite{ORS}, and
(\ref{qtarelidaha}) above):
\begin{align}\label{qtareli}
&t=q^2_{st},\  q=(q_{st}t_{st})^2,\  a=a_{st}^2 t_{st},\notag\\
&q^2_{st}=t,\  t^2_{st}=q/t,\  a_{st}^2=a\sqrt{t/q}.
\end{align}
I.e. we use the substitutions from the first line here
in $\hat{\h}^{j_\circ}_{\l,\,'\l}\,(q,t,a).$ 
Taking $j_o=min$ here (the main setting of the paper) 
is generally ``non-topological", 
though very reasonable algebraically. 

The case $t=q$ results in HOMFLY-PT polynomials,
i.e. we set  $k\!=\!1$ in $t=q^k$. 
Recall that $J$\~polynomials must be replaced 
in the definition of $\hat{\h}$ by
$P_\la^{(\!k=1\!)}$ in this case, which are Schur functions.   
Finally:
\begin{align}\label{conjhomflyy}
&\hat{\h}^{j_o}_{\l,\,'\!\l}\,
(q,t\!\mapsto\! q,a\!\mapsto\! - a)_{st}=
\hat{\hbox{\small H\!O\!M}}^{j_o}_{\l,\,'\!\l}\,
(q_{st}, a_{st}),
\end{align}

A combination of \cite{ChD2} with Section 7.1 from \cite{MoS}
(the case of iterated knots)
proves (\ref{conjhomflyy}) for any iterated {\em links\,}.
Another way to justify this coincidence is via a relatively
straightforward generalization of
Proposition 2.3 from \cite{CJ} (for torus knots),
where we used \cite{Ste}; see also \cite{ChD1}.
This approach is based on the DAHA shift operators and Verlinde
algebras. Also, instead of using \cite{MoS} or
the knot operators from {\em CFT\,}
(and the Verlinde algebras), one can directly apply here
the {\em Rosso-Jones cabling formula} \cite{RJ,Mo,ChE}
upon its relatively straightforward adjustment to iterated 
torus links.

\subsubsection{\sf Khovanov-Rozansky theory}
Let us restrict ourselves to
$\la^i=\yng(1)=\,'\!\la^j$
for all $i,j$ (the uncolored case).  Then $\{J_{\Box}\}_{ev}=
t^{1/2}(1+a)/(a^2)^{1/4}$ and
we conjecture that for
the hat-normalization of the stable KhR polynomials:
\begin{align}\label{khrconj}
\Bigl(\frac{\hat{\h}^{\emptyset}_{\,\l,\,'\!\l}}
{(1-t)^{\kappa+'\!\kappa}}\Bigr)_{\!st}\!=
\Bigl(\frac{(1\!+\!a)\,\hat{\h}^{min}_{\,\l,\,'\!\l}}
{(1-t)^{\kappa+'\!\kappa}}\Bigr)_{\!st}\!=
\hat{K\!h\!R}{}^{\hbox{\tiny \,stab}}_{\l,\,'\!\l}
(q_{st},t_{st},a_{st}).
\end{align}
The topological setting is {\em unreduced\,} here; recall that
$j_o=\emptyset$ in the first term 
means that we do not divide by the evaluations of
Macdonald polynomials at $t^\rho$.
Unreduced $K\!h\!R{}^{\!\hbox{\tiny\, stab}}$ are polynomials in
terms of $a_{st},q_{st}$, but their coefficients are
generally infinite {\em $t_{st}$\~series\,}.

The stable {\em Khovanov-Rozansky homology\,} is the
$sl_N$ homology from \cite{KhR1,KhR2,Rou} in the range of $N$
where the isomorphism in Theorem 1 from \cite{Ras} holds
(see also \cite{Kh}). 
They can be obtained for any $N$ from the
{\em triply-graded HOMFLY-PT homology}, assuming that the
corresponding {\em differentials\,} are known, which are
generally involved.

Let us also mention the relation to the
{\em Heegaard-Floer homology\,}: $N=0$.
Also, the
{\em Alexander polynomial\,}
of the corresponding singularity
is $\hat{\h}^{min}_{\l}\,
(q,q,a=-1)/(1-q)^{\kappa-\de_{\kappa,1}}$
in the case of one {\em uncolored\,} graph $\l$
with $\kappa$ paths (the number of connected components in the
corresponding link). This is the {\em zeta-monodromy\,} from
\cite{DGPS} upon $t\mapsto q$ (unless for the unknot).
The DAHA parameters are used here.
\vskip 0.2cm

The ORS conjecture, namely Conjecture 2 from \cite{ORS}, states
(in the unreduced setting) that 
$
K\!h\!R^{\hbox{\tiny\, stab}}_{\l,\,'\!\l}=
\overline{\mathscr{P}}_{\! \hbox{\tiny alg}},
$
where the latter series is defined there
for (the germ of) the corresponding plane curve singularity $\c$ 
from (\ref{yxcurve})
in terms of the {\em weight filtration\,} in the cohomology of
its nested Hilbert scheme. This (conjecturally) connects
 the DAHA superpolynomial upon the division from
(\ref{khrconj}) with
$\overline{\mathscr{P}}_{\! \hbox{\tiny alg}}$ 
under the  hat-normalization. See Section \ref{SEC:HILB}
above. 

There are also other conjectures connecting KhR polynomials
with rational DAHA, Gorsky's combinatorial polynomials
(for torus knots), Hilbert scheme of $\C^2$ and physics
superpolynomials (the name ``superpolynomials"
came from \cite{DGR}). We will not discuss these and
other related directions in the present paper.

\subsection{\bf Riemann Hypothesis}
We are now ready to state $RH$ for DAHA superpolynomials.
The notation is from the previous section.
We will use Theorem \ref{STABILIZ}; see also formula
(\ref{deg-a-jj}) for deg$_a\hat{\h}$.

\subsubsection{\sf The RH-substitution}
For any positive pair of graphs
$\{\l,\,'\!\l\}$, let
\begin{align*}
&\hat{\h}(q,t,a)=\hat{\h}^{min}_{\,\l,\,'\!\l}(q,t,a)=
\sum_{i=0}^d \h^i(q,t) a^i
\for d= \hbox{deg}_a\hat{\h}(q,t,a),\\
&H(q,t;a)=\hat{\h}(q\mapsto qt, t, a),\ \,
H^i(q,t)=\h^i(q\mapsto qt, t),\ 0\le i\le d.
\end{align*}
We also set  
$\hat{H}^i(q,t)\equal q^{-m}t^{-n}H^i(q,t)$ for the
minimal degrees $m,n$ of $q,t$ in $H^i$. 
Switching from $\h^i$ to their {\em hat-normalizations\,}
 $\hat{\h}^i$ with the constant term $1$ as above, 
one has $\hat{H}^i(q,t)=\hat{\h}^i(qt,t)$.
They are considered
as polynomials in terms of $t$, and we will almost entirely
switch below to $\om\equal 1/q$ from $q$. The super-duality 
for $H(q,t,a)$ is now for the map  $q\mapsto q, t\mapsto 1/(qt)$.

For any uncolored algebraic links, the $t$\~degree of 
$\hat{H}^0(q,t)$ is conjecturally the sum of $\de$\~invariants 
of its components plus $(\kappa-1)$ for the number $\kappa$ of the 
components, which follows from Conjecture 2.4 of \cite{ChD1};
for any $a$, see Conjecture 2.3 from \cite{ChP1}.  
Indeed, the top term in $\h^i$ is ``diagonal" for
uncolored algebraic knots, i.e. of the form $(q t)^{m_i}$ 
for any $i$ due to the DAHA super-duality; the passage to 
links results from (\ref{q-1-prod}).  
\smallskip

We mostly stick to rectangle Young diagrams in this work; 
the square diagrams are
especially valuable for us because they are 
transposition-invariant.
The (expected) geometric interpretation of DAHA superpolynomials
is directly related to $RH$ below; this is known for $\yng(1)$
and any columns \cite{ChP1,ChP2}. 
For non-square rectangles,
$l\times m$ (corresponding to $m\om_l$),
one can do the following {\em symmetrization\,}:
\begin{align}\label{symhpol}
H^i\mapsto
H^i_{sym}\equal H^i(q,t)+ H^i(q,1/(qt))
\for 0\le i\le \hbox{deg}_a\hat{\h}.
\end{align}
Then we employ the hat-normalization:
$\hat{H}^i_{sym}\equal \widehat{H^i_{sym}}=
q^{-m}t^{-n}H^i_{sym}$ for the minimal degrees $m,n$ 
of $q,t$ in $H^i_{sym}$. 
We {\em always\,}
switch to $\hat{H}^i_{sym}$ below
if the transposition of all Young 
diagrams involved changes the isomorphism class of the 
diagram/link $\l$
or the pair $\{\l, '\!\!\l\}$.
\smallskip

\comment{
We must note that the DAHA superpolynomials
$\hat{\h}$  with \Yboxdim5pt \yng(2,2)\Yboxdim7pt\
can be huge and difficult to calculate practically,
even for simple torus knots. For instance, it took about
$4$ days to calculate $\hat{\h}$ with this diagram for
$T(5,3)$, which was one of the key ``colored" tests 
for us.}


\subsubsection{\sf Algebraic RH}
We assume that $\om=1/q$ is real positive.
Recall that $\hat{\h}=\hat{\h}^{min}_{\l,\,'\!\l}
(\,\lla,\,'\!\lla\,;\,q,t,a)\,=\sum_i\h^i a^i$
for $0\le i\le$ deg$_a \hat{\h}$. Also $\hat{H}^i$ is
the hat-normalization of $\h^i$ upon the substitution
$q\mapsto qt$; they  
are polynomial in terms of $q,t$ with the constant term $1$.
For the sake of definiteness, we assume that
deg$_q\!\le$ deg$_t$ in $\hat{\h}(a=0)$, employing the
super-duality if necessary. For instance, $l\ge m$
will be always imposed for {\em knots\,} colored by 
rectangles $l\times m$ (columns instead of rows). 

Furthermore, 
let $\h^i_\bullet=\h^i$,
 $\hat{H}^i_\bullet=\hat{H}^i$ if the corresponding
colored link is self-dual with respect to the transposition
of the diagrams in $\lla,\,'\!\lla$ (and the equivalence
of graphs). Otherwise 
$\h^i_\bullet=\h^i_{sym}$, 
$\hat{H}^i_\bullet=\hat{H}^i_{sym}$.

\begin{conjecture} \label{CONJLIM}
(i) For any pair of iterated links $\l,\,'\!\l$
(possibly non-algebraic), any colors
(Young diagrams) and arbitrary $0\le i\le $deg$_a\hat{\h}$,
the following limits exist and coincide:
\begin{align}\label{h-dag-0}
\hat{H}^i_{\dag}\!=\varsigma_i t^{\pi_i}\,S^i(t)\!\equal\! 
\lim_{q\to 0}\hat{H}^i_\bullet(q,t\!\mapsto\! 
\frac{t}{q^{\frac{1}{2}}})\!=\!t^{-\pi_i}
\lim_{q\to 0}\h^i_\bullet(q,t\!\mapsto\! \frac{t^2}{q}),
\end{align}
where $\varsigma_i=\pm$,
$\pi_i\ge 0$, $S^i(t\!=\!0)=1$. Moreover, the relation
$2(\pi_i\!+\!\si_i)\!=\!$ deg$_t \hat{H}^i$ always holds, 
where we set\, $2\si_i\equal$ deg$_tS^i$.

(ii) All zeros of the polynomials $S^i$ are roots of
unity if $\l,\,'\!\l$ are colored only by rectangle diagrams.
They are simple for 
any iterated torus links if either $i=0$ or 
they are uncolored. Furthermore,  
$S^i=(1-t^{2\si_i+2})/(1-t^2)$
for uncolored algebraic 
knots and  $S^i=(1+t^{2\si_i})$  for those colored by 
columns upon the symmetrization
(including \,\Yboxdim7pt$\yng(1)$\ \Yboxdim7pt); $i$ are
arbitrary and $\pi_i=0$ in such cases. For uncolored
algebraic links,  $\pi_0$ coincides with the number of components
minus $1$. 
\sq
\end{conjecture}

The passage to $t^2/q$ in $(i)$
is not difficult to justify; thus $S^i$ are
actually polynomials in terms of $t^2$. This is somewhat
parallel to \cite{OnRS} and a recent work \cite{GORZ}
of Griffin, Ono, Rolen, and Zagier. 
We expect that at least in the uncolored case, the invariants
$\pi_i$ and $S^i(t)$ have topological meaning 
beyond iterated torus links.
For instance, if the stable, reduced and hat-normalized
KhR polynomials at $a=0$  are
used instead of  $\hat{\h}$, then $(\pi_i\!+\!\si_i)$ 
is presumably the {\em 4-genus\,} for any 
{\em positively iterated torus knot\,} 
(not only algebraic). See 
\cite{RasG} for the connection to the
Khovanov polynomial and other invariants and also \cite{Hed}. 
For instance, assuming that $i\!\le\! 2j$ in
$\hat{K\!h\!R}_{\hbox{\tiny red}}^{\hbox{\tiny stab}}(a_{st}\!=
\!0)=\sum_{i,j} c_{i,j} q_{st}^{2i} t_{st}^{2j}$ (possibly
upon some adjustment of $i,j$), $\pi_0$ and $S^0(t)$
can exists for the sum of borderline 
terms \ $\sum_i c_{2i,i} q_{st}^{4i} t_{st}^{2i}$ not only for
algebraic knots.

For our {\em geometric superpolynomials\,} of
algebraic knots, $\pi_i\!+\!\si_i=\de$, 
e.g. $=(\rr-1)(\ss-1)/2$ for $T(\rr,\ss)$. This
becomes significantly more involved for non-rectangle
Young diagrams. It seems that
{\em non-cyclotomic factors of $S^i$ appear only for 
non-rectangle diagrams\,} for any knots.

The diagram \Yboxdim5pt$\yng(2,2)$\ \Yboxdim7pt deserves
a special comment; no symmetrization is necessary here. 
In this case, there are no multiple zeros of 
$S^i$ for algebraic knots in the examples we reached
(though not many)
for any $i\ge 0$. Also, $\pi_i=0$ for even $i$ and  
$\pi_i=1$ for odd $i$, when $\hat{H}^i$ always
has two {\em trivial\,} irregular zeros $\{-1,-\om\}$. 
Correspondingly (conjecturally) for even and odd $i$:
$\varsigma_i t^{\pi_i}S^i=\frac{1-t^{2\si_i+4}}{1-t^4}$ and
$\varsigma_i t^{\pi_i}S^i=t\frac{1-t^{2\si_i+2}}{1-t^4}$. 

Part $(i)$ is expected to hold under 
{\em antisymmetrization\,}
instead of symmetrization, which is 
$H^i_{a\!sym}\!\equal H^i(q,t)- H^i(q,1/(qt))$ followed
by $H^i_{a\!sym}\!\mapsto\! \hat{H}^i_{a\!sym}$. As above,
we assume that deg$_q\!\le$ deg$_t$ in $\hat{\h}(a=0)$; 
in particular, $l> m$ for non-square diagrams $l\times m$. 
Then the properties of $\pi_i,S^i$ generally
become somewhat ``better" vs. the symmetrization.
Extra {\em trivial\,} zeros of $\hat{H}^i$ now can emerge, 
which are $t=\pm\sqrt{\om}$. Presumably
$S^i=t^{2\si_i}-1$ in this case for any
algebraic knots colored by columns and for any $i$.  

\subsubsection{\sf Analytic aspects}
Conjecture \ref{CONJLIM} gives that
the number of
non-$RH$ (irregular) zeros of $\hat{\h}^i$ in the
vicinity of $q=1/\om=0$ is no greater that $2\pi_i$ 
plus the number of non-unimodular zeros of $S^i$ 
and multiplicities of its multiple zeros. It can
be smaller when some multiple
zeros of $S^i$ are unimodular (not always).
 We arrive at the following.

\begin{conjecture} \label{CONJRH}
(i) {\sf Strong RH\,}.
For an arbitrary uncolored algebraic knot
and any given $0\le i\le $deg$_a\hat{\h}$,
there exists $\om'=\om'_i>0$ such
that for all $\,\om>\om_i'$, 
the $t$\~zeros $\,\xi$ (if any) of $\hat{H}^i$ are all
complex, simple
and satisfying the {\sf $RH$\~equality}: $|\xi|=\sqrt{\om}$.
 
(ii) The same holds
for algebraic knots colored by 
\Yboxdim5pt$\yng(2,2)$\ \Yboxdim7pt if the
{\sf trivial zeros\,}  $\xi\!=\!\! -1,-\om$ are omitted,
which occur in $\hat{H}^i$ if and only if  $\,i\,$ is odd. 
Also, $\hat{H}_{sym}^{i}$ satisfies
$(i)$ for columns (any $i\ge 0$) and
non-square rectangles $l\times m$ with $l>m$ for $i=0$. 

(iii) {\sf Weak RH\,}.
Given an uncolored algebraic link, there exists $\om'_0$ such
that the number of pairs of {\sf stable irregular zeros\,}
of $\hat{H}^{i=0}$
(satisfying $|\xi|\neq\sqrt{\om}$) for $\om>\om'_0$ equals
the number of components minus $1$. \sq 
\end{conjecture}

\comment{
(iv) RH from (i) or (iii)
holds for $\om\ge 2$ for uncolored algebraic links when $i\!=\!0$
or for any $i\ge 0$ in the case of uncolored torus knots. For 
the latter, $\hat{H}^i$ at $\om\!=\!1$ are products of 
cyclotomic polynomials for any $i\ge 0$ 
and $RH$ holds in a neighborhood of $\om=1$.

(iv) Upon the {\sf antisymmetrization\,} 
$H^i_{a\!sym}\!\equal H^i(q,t)- H^i(q,1/(qt))$ followed
by $H^i_{a\!sym}\!\mapsto\! \hat{H}^i_{a\!sym}$
for rectangles $l\times m$ with $l>m$,
$RH$ holds for algebraic knots (any $i$)
and for algebraic links when
$i=0$ and all components have the 
same color ($\xi=\pm\sqrt{\om}$ are disregarded).
}

A connection is expected with the 
{\em spectral zeta-functions\,}, especially in
the case of Schottky uniformization of Riemann surfaces; 
see e.g. \cite{CM}. The following conjecture is of
this nature.


\begin{conjecture}\label{CON:om2}
For any uncolored algebraic {\sf knot}, 
$\varpi_i\equal \hbox{\rm inf } \om'_i$ for $i=0$
is smaller than $2$. Moreover, 
$\lim_{m\to\infty}\varpi_0=2$ for  
$\r=\C[[z^4,z^6\!+\!z^{7+2m}]]$ corresponding to 
the cables $C\!ab(13+2m,2)T(3,2)$, which
sequence of $\varpi_0$ is actually increasing when 
$\,m\hbox{\rm\,  mod }4$ is fixed. Also,
$\sup_{m}\varpi_1=${\small $2.2132458\ldots$},
$\sup_{m}\varpi_2=1.$
 \sq
\end{conjecture}

This is actually the only ``family" with $\varpi_0$
near $2$ we
found. Let us provide  $\varpi_i$ for
$i=0,1,2$ in the following $3$ cases:
{\small
\begin{align*}
&C\!ab(13\ ,2)T(3,2): \varpi_i= &1.495583269,\ 2.176487419,\
0.9430445115, \\
&C\!ab(113,2)T(3,2): \varpi_i= &1.993679388, \ 2.134669951,\
0.9955504853, \\
&C\!ab(313,2)T(3,2): \varpi_i= &1.997705951, \ 2.210868584,\
0.9992171315.
\end{align*}
}
Using formula (\ref{fam3-2}) for $m=110,210,310,510,1010,1510$,
we have:
{\small
$\varpi_0=  1.996921, 1.998340, 1.998863, 1.999303, 1.999645,
1.99976261$}
and the corresponding 
{\small
$q\!=\!\frac{1}{\varpi_0}\!=\!0.500771, 0.500415, 0.500284, 
0.500174, 0.500088$}, {\small $0.50005935$}. 
The convergence to $\om=2$ is the best for\,
$m\!\!\mod\! 4=1,2$.
Such limits can be
generally calculated using the conjectural coincidence 
of $\hat{\h}$ with the geometric superpolynomials from
\cite{ChP1}. 


This family also has the largest $\varpi_i$ for $i=1,2$
among uncolored algebraic knots we considered. Here 
deg$_a\hat{\h}\!=\!3$, 
and $\hat{H}^3\!=\!1$; $\varpi_i(i>0)$ exceed $2$ only 
for this family in our ``database". 

The greatest $\varpi_0$ we found so far among uncolored 
algebraic {\em links\,} is $2.062433590332$ for
$(C\!ab(5,3),C\!ab(4,3))T(1,1)$ corresponding to the
coinvariant 
$\{\ga[1,1](\ga[3,1]P(1)\ga[3,2]P(1))\}$ in the notation
from the table. Here Weak $RH$ holds for all $i$. The
corresponding singularity is
$(x^5\!-\!y^3)(x^3\!-\!y^4)\!=\!0$ with 
$Z=${\small $1+q^7+q^8+q^{14}+q^{15}+q^{16}+
q^{22}+q^{23}+q^{30}$} and $lk=9$ (in this table,
$\varpi_0<2$ for uncolored algebraic links). 

\subsubsection{\sf Comments}\label{SEC:COMM}
For algebraic knots colored by columns, 
the statements of Conjecture \ref{CONJLIM} can be verified if one
switches from DAHA superpolynomials $\hat{\h}$ to 
geometric (motivic) ones. The proof goes as follows. Only $J_\De$  
with dimensions no smaller than $dev(\De)$ contribute
to the limit  of $\h_{mot}$ from
Definition \ref{DEF:MOT}.
 Given $dev(\De)$, such $\De$  occurs only for a {\em single\,}
$\De$, namely, the one obtained by adding {\em consecutive\,} 
gaps to $\Ga$ starting with the top one, which is
$\max\{\Z_+\setminus \Ga\}$. 

The generalization to any powers of $a$ and columns 
follows from \cite{ChP1,ChP2}. Part $(i)$ of
Conjecture \ref{CONJLIM}
for (at least) algebraic knots upon $a=0$ and in the
case of columns can be managed within the DAHA 
theory. We use that only positive powers of $\tau_{\pm}$
appear in the formulas and that $X_{\om_i}$ are 
{\em nonsymmetric\,} Macdonald polynomials.

An extension of Part $(ii)$ there to
{\em rectangle\,} Young diagrams larger that
\Yboxdim5pt$\yng(2,2)$\ is likely, but there is a lack
of (numerical) examples. Also, the restriction 
$i\!=\!0\,$ in our conjectures can presumably be replaced  
by ``for any even $i$", but we will stick to $i\!=\!0$
in this paper. 
\vskip 0.3cm

{\sf Torus knots.}
The polynomials $\hat{H}^i$ at $\om\!=\!1$ are products of 
cyclotomic polynomials for {\em uncolored torus knots\,} and  
any $i\ge 0$. They can be calculated
explicitly; for $a=0$, the formula can be deduced from
the {\em shuffle conjecture\,} (proven in \cite{CaM}). In the 
absence of multiple zeros at $\om\!=\!1$ (which is always true for 
$i\!=\!0$), this implies $RH$ in some neighborhood of $\om\!=\!1$. 
However this is not generally the case for $i>0$ (for
torus knots); multiple zeros do occur. We note that 
$\varpi_0$ can be beyond $1$ for uncolored torus knots
for sufficiently general torus knots. 
Then  $RH$ fails somewhere between $\om=1$ 
and $\varpi_0$.
\vskip 0.1cm

We mention that there are sufficiently
explicit known/conjectured formulas for
uncolored torus knots $T(m,km+1)$. 
See \cite{GM,Mel1,ChP1} and also 
\cite{ORS,DMMSS,FGS}; the proof of the formula for DAHA
superpolynomials colored by any rows for $T(2m+ 1,2)$ is in
\cite{CJJ}. Even when explicit formulas are known/conjectured,
they are generally involved for explicit finding $\varpi_i$,
but can be helpful (at least) to examine the point $\om=1$
for torus knots and obtain {\em family superpolynomials\,}
(actually rational functions)
as was done in (\ref{2z+1},\ref{3z+1},\ref{fam3-2}).
\vskip 0.2865cm
\vfil

Note that $H^i(\om\!=\!1)$ upon $t\mapsto q_{st}^2$ is
the $a_{st}^{2i}$\~coefficient of  the corresponding HOMFLY-PT 
polynomial up to $\pm q_{st}^{\bullet}$.  
See (\ref{conjhomflyy}). This coincidence for torus knots
was justified in \cite{CJ} using \cite{Ste}. See also
\cite{MoS,ChD2} for any
iterated torus knots and links. 
\vskip 0.3cm

{\sf Asymptotic class numbers.}
A natural application of Conjecture \ref{CON:om2}
is to the growth estimates for $|J(\mathbb F_q)|$ for
the Jacobian factors $J=J_\r$ from Definition \ref{DEF:MOT}.
This is a classical track (for any curves). 
The conjecture that $RH$ holds for  
$\om =1/q\le 2$ provides some estimates for  
$|J(\mathbb F_q)|$ with any $q$. 
Here we obtain ``pure singular"
contributions; the smooth case is covered by the Weil $RH$.  

\vskip 0.2cm
The switch to the {\em families\,} and the family polynomials 
$\hat{\mathbb{H}}$ from (\ref{hsfam}) seems the most
relevant here. This can potentially clarify the 
parallelism between the Alexander knot polynomials and
{\em Iwasawa theory\,} observed by B.Mazur; see e.g. \cite{Mor}
and around (\ref{Iwasa}) below.
The Iwasawa polynomials  
give exact formulas for the growth of ideal class groups in
$\Gamma$-extensions; we do the towers of 
Puiseux field extensions in the {\em families} for $i=1$.  
This is connected with the so-called 
{\it Drinfeld-Vladut bound},
but superpolynomials 
provide {\em exact\,} formulas
(as in the Iwasawa theory); see e.g. \cite{GaS}.
 
\subsection{\bf Non-classical features}
Similar to the classical Riemann zeta,
all $t^{\bullet}\hat{H}^i$ are real-valued at
$U_{\sqrt{\om}}\!\equal\!\{z\!\in\! \C\mid
|z|\!=\!\sqrt{\om}\}$  for proper powers
of $t\,$ due to the super-symmetry (our functional equation)
and the reality of $\hat{H}^i$. Also, if $\xi$ is a zero of
$\hat{H}^i$, then
$\sqrt{\om}/\,\overline{\xi}$ is its zero with the same angle, 
where $\overline{\xi}$
is the complex conjugation of $\xi$.

\subsubsection{\sf The range of 
\texorpdfstring{{\mathversion{bold}$\om$}}{omega}}
There are new opportunities here
vs. the classical theory,
since $\om$ is arbitrary for us
(a free parameter). 

\begin{lemma}\label{LEMOM}
Assuming that $\om'_i$ from 
Conjecture \ref{CONJRH},(iii) exists for
some $i\ge 0$, let $\varpi_i$ be (as above) the lowest
such $\om'_i$. Then $\varpi_i$ is a (real) root
of the {\sf reduced $t$\~discriminant\,} $D^i$ of $\hat{H}^i$,
which is the product of all simple factors in the actual 
discriminant of $\hat{H}^i$. 
For instance, $\varpi_i$ is an algebraic number; it  
coincides with the greatest real zero 
$\om^{top}_i$ of  $D^i$ if we add the simplicity of zeros
$\xi$ for $\om>\om'_i$ to the definition of $\om'_i$. 

Also, assuming that $\hat{H}^i$ satisfies Weak
$RH$ in an interval beyond (greater than) $\om_i^{top}$,
Weak $RH$ holds then for all $\om\ge \om_i^{top}$. The same is 
true if Weak $RH$ holds in the interval 
$\om_i^{top}-\ep<\om<\om_i^{top}$ for some $\ep>0$ and
$\hat{H}^i$ has no multiple roots at\ \,
$\om=\om^{top}_i$ of norm $\sqrt{\om}$.
\end{lemma}
{\it Proof.} Here we use that the zeros apart from
$U_{\sqrt{\om}}$ vanish or emerge only at 
(real) zeros of $D^i$.
Indeed, they appear in pairs $\{z,z'\}$ with coinciding
angles and therefore create multiple zeros of $\hat{H}^i$
when approaching $U_{\sqrt{\om}}$.

In particular, if Weak $RH$ holds for $\om$ in an interval 
greater than  $\om^{top}$, 
then the formation
of a non-$RH$ pair (i.e. that apart from $U_{\sqrt{\om}}$)
at some  $\om>\om_i^{top}$ results in a multiple zero
of $\hat{H}^i$ beyond $\om_i^{top}$, which is impossible.
Similarly, if Weak $RH$ holds for
$\om_i^{top}-\ep<\om<\om_i^{top}$, then the multiple
zeros at $\om=\om_i^{top}$ can emerge only from some pairs of
zeros of norm $\sqrt{\om}$. \sq

\subsubsection{\sf Non-RH zeros}
By $RH$ for links, we will always mean Weak $RH$ 
from Conjecture \ref{CONJRH}, $(iii)$, allowing multiple zeros and
$\kappa-1$ super-dual pairs of (stable)
irregular zeros for $\kappa$ branches. Accordingly, 
\vskip -0.6cm
\begin{align}\label{varpi-noprime}
\varpi_i\equal\hbox{inf}\, \{\om'\mid\!
\hbox{\, Weak RH holds for\, }\hat{H}^i
\hbox{\, for\, }
\om\ge \om'\};
\end{align}
\vskip -0.1cm
\noindent
this is a real zero of $D^i$ from the lemma. 

Thus, $\varpi_i$ conjecturally exist for  
uncolored algebraic {\em knots\,} for any $i\ge 0$.
Also, they exist for algebraic knots colored by
$(a)\,$ 
\Yboxdim5pt$\yng(2,2)$\Yboxdim7pt\,,
where {\em trivial zeros\,} $-1,-\om$ are excluded,
$(b)\,$ columns $\om_l$ upon the symmetrization,
$(c)\,$ any rectangles upon the
symmetrization and for $i=0$.


Concerning (\ref{varpi-noprime}), it is
possible that $\varpi_i<\om_i^{top}$; then multiple
zeros appear after $\varpi_i$; this occurs only twice
in the table.
Also, there are no reasons for ``nice" formulas for
$\varpi_i,\om_i^{top}$, but at least they are algebraic numbers
(zeros
of $D^i$) and can be calculated as exactly as necessary.

\comment{
The claim from $(iv)$ that $\varpi'_i<2$
for uncolored algebraic knots may reflect the
fact that plane curve singularities do not have places
(prime numbers) of {\em bad reduction\,}, when considered up
to {\em topological\,} equivalence. This argument is
not rigorous and we have some understanding of
bad reductions only for uncolored knots and when $i=0$.}

If $\varpi_i=0$, then  Weak $RH$ holds
for any $\om>0$. This happens for any uncolored 
$T(2m\!+\!1,2)$ and for iterated Hopf-type links
$2T(m,1)$. Here deg$_a=2$ and $i=0,1$. For the Hopf $2$-link,
$\varsigma_0 t^{\pi_0}S^0=
t\frac{1-t^{2m}}{1-t^2}$. It corresponds to 
the $2$\~branch singularity $(x^m-y)(x^m+y)=0$ with the 
linking number $lk=m$. One pair of irregular (real) zeros
of $H^0$ approaching $\{1,\om\}$ as $\om\to\infty$
occurs here after $\om^{top}_0=(m+1)^2/m^2$;
the other zeros for any $\om>0$ are of norm $\sqrt{\om}$.
Thus $\om_0^{top}>\varpi_0$ in this case. There is only one 
instance of non-Hopf link in the table when this happens: 
entry $35$ ($N_z^0=26$).

\section{\sc RH numerically}
We calculated quite a few examples in the range deg$_a\le 8$
(and sometimes beyond). The most instructional ones are collected
in the table below, though we provide many examples beyond it. 

\subsection{\bf Table organization} It is based on the 
$\hat{\h}$\~polynomials above with the following deviation:
we consider $\hat{H}^i$ {\em not\,} only
for algebraic links/knots. The DAHA construction is
fully applicable without any positivity conditions
for $\rr,\ss$, though the connection with plane
singularities will be lost then.
Numerical experiments clearly indicate that the class
of iterated links with $\hat{\h}$\~polynomials satisfying
$RH$ is wider than algebraic links only
(combinatorially, {\em positive\,} pairs of graphs).
\vskip 0.2cm

\subsubsection{\sf Main notations}
{\em By $RH$, we will mean below Weak RH\,} from
$(iii)$ of Conjecture \ref{CONJRH}; multiple zeros
will be allowed.
The number of pairs $\{\xi,\om/\xi\}\,$
of {\em irregular\,} zeros of $\hat{H}^i$ is shown in the last
column of the table below after "$-$".
Since $|\xi|>\sqrt{\om}$ (assuming the simplicity) for
one of the zeros of a non-$RH$ pair, they significantly influence
expected counterparts of Weil-style estimates
(those in his proof of $RH$).

{\sf Status.}
For knots and links we simply put ``knot" or ``link" in
the corresponding entry; ``alg" or  ``alg\kern-13pt{\bf=}"
are naturally stand for algebraic and non-algebraic
knots/links. The first column provides the status of
Weak $RH$
with the following 3 options:

``{\em HOLDS\,}", if it holds for all $a^i$,\ \,\,
``{\em FAILS\,}", when it fails at $a\!=\!0$,

``{\em OK $a\!=\!0$\,}", \,if it holds for $a\!=\!0$ but
fails for some other $a^i\, (i\!>\!0)$.
\smallskip

In the ``{\sl OK\,}"-case, we give the range of $a^i$ when $RH$
holds; the total number of polynomials
$\hat{H}^i$ is also provided, which is (deg$_a+1$). We
calculate only until the first failure of $RH$ in this table;
then we stop.
We numerate the examples and provide
he corresponding numbers $N^0_z$ of all zeros of $\hat{H}^0$
(only for $i\!=\!0$) after ``No". See the first column in the
table below. Note that $N^0_z=2\de$ for algebraic knots.
\vskip 0.2cm

The last two columns contain $\varpi_0$ and the maximum, 
denoted by
$\varpi_{1\ldots}^{max}$,
of $\varpi_1,\varpi_2,\ldots$ calculated till the first
failure of $RH$. They
can be zero  when all $\om>0$ satisfy $RH$.
We put "$-$" for $\varpi_0$ if $RH$ fails at $a=0$,
and make $\varpi_{1\ldots}^{max}=$"$-$" if
$RH$ fails for $a=0$ or fails for $a^1$.

\comment{
The last column contains
the number of pairs of
{\em irregular\,} zeros at $i\!=\!0$ and the
corresponding maximum calculated for $i\!=\!1,\!2,\ldots$
until $RH$ fails for the first time. We provide $N_{irz}$
calculated right after $\varpi_i$ and then (after "$-$")
its stable value for $\om>>0$.
}
\vskip 0.2cm
\Yboxdim5pt

\subsubsection{\sf DAHA formulas}\label{SEC:ABREV}
The 4{\tiny th} column contains the
DAHA formula used to calculate the corresponding $\hat{\h}$,
where $\{\cdots\}$ means the coinvariant and $\ga[\rr,\ss]$
is understood as the lift of
$\ga_{\rr,\ss}\in PSL_{\,2}(\Z)$ to
$\hat{\ga}_{\rr,\ss}\in PSL_{\,2}^{\wedge}(\Z)$. We
somewhat abuse the notations by omitting
$H\!\Downarrow=H(1)$; the $a$\~stabilization is also
assumed.
For instance, $\{\ga[3,2]\ga[2,1](P)\}$
must be understood as
$\{\hat{\ga}_{3,2}(\hat{\ga}_{2,1}(P_{\yng(1)})\!\Downarrow)\}$,
{\em not\,} as a product of two $\ga$ inside the coinvariant,
and upon the further $a$\~stabilization.
\vskip 0.2cm

Also, we actually use the $J$\~polynomials (not $P$\~polynomials)
for links and always do division by the $LCM$ of
the evaluations at $t^{\rho}$; see the definition of
the {\em\  min\,}-normalization of $\hat{\h}$. In the table,
we simply put $P$  assuming the rest. There are the
following abbreviations:  $P$ stands for $P_{\yng(1)}=P_{\om_1}$,
$P(1+1),P(2+0), P(2+1), P(2+2)$, and $P(3+3)$ stand for
$2\om_1=\yng(2)$ and $\om_2=\yng(1,1)$\,,
$\om_1+\om_2=\yng(2,1)$ (the hook), $2\om_2=\yng(2,2)$\,,
 $2\om_3=\yng(2,2,2)$\,.
Recall that if some of $\rr,\ss$ are
negative then the corresponding link/knot
is non-algebraic. 
\Yboxdim7pt

\vskip 0.2cm
\subsubsection{\sf Basic cables}
The cab-presentations in the table are partial:
we omit the {\em arrows\,}
and do not show at which vertex the
corresponding polynomials are inserted. This can be
seen from the DAHA presentations. For instance, entry
$46$ with the number of zeros $N_z^0=72$ is represented by the 
graph
$\twoone$ with both arrows colored by 
$\yng(1)\,$ (omitted)
and the $\tax$\~labels $[5,2]$ (for the left vertex)
and $[2,1]$. Also, applying $P(Y)$ in the DAHA-formula
means that we consider a
pair of graphs: the one
shown in the table and
$'\!\l$ that is a pure arrow $\rightarrow$ colored by $\yng(1)\,$
(without vertices).

To give another example,
$\bigl(C\!ab(11,3),C\!ab(11,3)\bigr)T(3,1)$ from entry
No$=45$ (with
$94$ zeros) means the tree {\large $\twotwo$},
where the first vertex is labeled by $[3,1]$ and the other
two by (coinciding) labels $[3,2]$ with the arrows colored
by $\yng(1)$\,. $RH$ holds here for all $\hat{H}^i$ (from $i=0$
through $i=$deg$_a=5$). This link is algebraic, $\varpi_0=
1.39031$ (we calculate them with much greater accuracy),
$\varpi^{max}_{1,\ldots}=1.55923$, which is the maximum of
$\varpi_i$ from $i\!=\!1$ to $i\!=\!6$ in this case. Generally,
this
maximum is taken from $\hat{H}^1$ till the last $\hat{H}^i$ 
where $RH$ holds.

The last column gives the corresponding
number of pairs of {\em irregular\,} zeros,
where $N_{irz}^{max}$ is
calculated in the same manner; the first  $N_{irz}$
is calculated right after $\varpi$ and the second one
(after "$-$") is the stable number
of such pairs (for any large $\om$).

\subsection{\bf On Conjecture \ref{CONJLIM}}
The presence of {\em non-unimodular\,} roots
of polynomials $S^i(t)$ from
Conjecture \ref{CONJLIM} is interesting. The only examples
we found so far are only for
non-rectangle diagrams, which includes links (possibly
non-algebraic). Many examples were considered beyond the
table and for all $i$, not only those till the first failure 
of $RH$, which are presented in the table below. 

\subsubsection{\sf The case of 3-hook} We will give here
all links from the table where $S^i$ are not products
of cyclotomic polynomials, also providing the corresponding 
$\pi_0,S^0$, even if $S^0$ are such products. Actually all
$S^0$ are products of cyclotomic
polynomials within the table (so only $S^{i>0}$ 
can be non-cyclotomic), but 
is not always the case. We will provide below an example
when $S^{0}$ is non-cyclotomic. Also, {\em In all examples
we considered  $\ \varsigma_i t^{\pi_i}=-t^3$}, so we will omit it. 
This factor  contributes $3$ to the number of stable irregular 
(non-$RH$) pairs of zeros of $\hat{H}^i$. Thus, let us focus 
on $S^i$.

\vskip 0.4cm
(1) $T(5,2), \{\ga[5,2](P(2+1))\}\, :$
\vskip 0.2cm

\noindent
($i=0$) \ \  $\varsigma_0 t^{\pi_0}S^0=-t^3 
\bigl(\frac{1-t^6}{1-t^2}\bigr)\frac{1-t^{12}}{1-t^2},$
where the quantity $\bigl(\ldots\bigr)$ is the multiple
part of $S^0$. So there are potentially $4$ pairs of 
non-$RH$ zeros due to $S^0$ and $3$ pairs due to $t^3$;
this matches entry $No=57$ with $N_z^0=20$ from the table.

\vskip 0.2cm
\noindent
($i=1$)\ \  $\varsigma_1 t^{\pi_1}S^1= 
-t^3 (1 + t^2) 
(1 + 2 t^2 + 2 t^4 + 3 t^6 + 2 t^8 + 2 t^{10} + t^{12})$,
where the latter factor is irreducible and {\em non-cyclotomic\,}. 
It has $4$ non-unimodular zeros, which matches the total number
of $6$ non-$RH$ pairs for $\hat{H}^1$ (this is not in the table).


\vskip 0.4cm
\vfil
(2) $T(4,3), \{\ga[4,3](P(2+1))\}\, :$
\vskip 0.2cm

\noindent
($i=0$)\ \  $S^0= 
 \frac{1-t^{20}}{1-t^2}\frac{1-t^{8}}{1-t^2},$
where the multiple zeros come from $\frac{1-t^{4}}{1-t^2}.$
The expected number of non-$RH$ pairs is therefore $5$,
which matches entry $59$ with $N_z^0=30$ in the table.


\vskip 0.2cm
\noindent
($i=1$) \ \ $S^1= \frac{1-t^{9}}{1-t^2}
(1 + 2 t^2 + 2 t^4 + t^6 + t^8 + t^{10} + 2 t^{12} + 
2 t^{14} + t^{16})$,
where the latter factor is irreducible and {\em non-cyclotomic\,}.


\vskip 0.4cm
(3) $T(5,3), \{\ga[5,3](P(2+1))\}\, :$
\vskip 0.2cm

\noindent
($i=0$)\ \  $S^0= \frac{1-t^{28}}{1-t^2}
\frac{1-t^{10}}{1-t^2},$
where there are no multiple roots of $S^0$; this
matches the number of non-$RH$ pairs, which is $3$,
in entry $60$ with $N_z^0=40$.

\vskip 0.2cm
\noindent
($i=1$)\ \  $S^1=\frac{1-t^{4}}{1-t^2}
\frac{1-t^{6}}{1-t^2}
\frac{1-t^{8}}{1-t^4}\frac{1-t^{10}}{1-t^2}
(1 - t^4 + t^6 + t^{10} - t^{12} +  t^{16})$,
where there are no multiple roots and the last factor
is irreducible {\em non-cyclotomic\,}. 


\vskip 0.2cm
\noindent
($i=3$)\ \   $S^3= (1 + t^2)^2 \frac{1-t^{6}}{1-t^2}
(1 + 2 t^2 + 2 t^4 + 2 t^6 + t^8 + 
   t^{10} + t^{12} + 2 t^{14} + 2 t^{16} + 2 t^{18} + t^{20})$,
where the latter factor is {\em non-cyclotomic\,};
we have totally $2$ multiple roots and $4$ non-unimodular
zeros, so the expected number of non-$RH$ zeros is
$4+4+6=14$. This corresponds to the actual number
for $\om>3.402358077$ (not in the table).


\vskip 0.4cm
(4) $T(7,3), \{\ga[7,3](P(2+1))\}\, :$
\vskip 0.2cm

\noindent
$(i=0)$\ \  $S^0=\frac{1-t^{44}}{1-t^2}
\frac{1-t^{14}}{1-t^2}$, which has no multiple
nonzero zeros, and the number of pairs
of non-$RH$ zeros in entry $61$ with $N_z^0=60$ is
$3$ indeed.  

\vskip 0.2cm
\noindent
$(i=1)$ $S^1=
(1+t^2) (1+2 t^2+3 t^4+4 t^6+5 t^8
+6 t^{10}+6 t^{12}+7 t^{14}+6 t^{16}+7 t^{18}+6 t^{20}+7 t^{22}
+6 t^{24}+7 t^{26}+6 t^{28}+7 t^{30}+6 t^{32}+7 t^{34}+6 t^{36}
+7 t^{38}+6 t^{40}+
6 t^{42}+5 t^{44}+4 t^{46}+3 t^{48}+2 t^{50}+t^{52}),$
where the last factor is irreducible {\em non-cyclotomic\,}.

\vskip 0.2cm
\noindent
$(i=3)$\ \  $S^3=
(1+t^2)^2 \frac{1-t^{14}}{1-t^2}
\frac{1+t^{10}}{1+t^2} 
(1+t^2+t^4+t^6-t^{12}
+t^{18}+t^{20}+t^{22}+t^{24})$,
with {\em non-cyclotomic\,} last factor,
$4$ non-unimodular zeros and $2$ multiple zeros 
$\pm\, \imath$
due to $(1\!+\!t^2)^2$. The total number of non-$RH$ pairs
is expected $2+(4/2)+3=7$, which is actually greater
than the actual number $5$; these $5$ pairs occur 
after a huge $\varpi_3=159557.4798$. 
Thus, the $q$\~deformations of $\pm\, \imath$ are $RH$\~zeros
for any $\om>\!>0$ in this case.


\subsubsection{\sf Non-cyclotomic 
\texorpdfstring{{\mathversion{bold}$
S^{i=0}$}}{{\em S } for zero}}
\label{SEC:2-1}
The last example we will provide is for
the $3$\~hook \Yboxdim5pt \yng(2,1)\,\Yboxdim7pt\ 
where $S^0$ is not a product of cyclotomic polynomials
(even for $i=0$). This is not in the table.
\vskip 0.2cm

$T(7,2), \{\ga[7,2](P(2+1))\}\, :$ 
\vskip 0.1cm

\noindent
$(i=0)$\ \  $\varsigma_0 t^{\pi_0}S^0=
-t^3 \frac{1-t^{10}}{1-t^2}
(1 + t^2 + 2t^4 + t^6 + t^8 + t^{10} + 2t^{12} + t^{14} + 
t^{16})$,
where the last factor has $4$ pairs of non-unimodular zeros, 
which results in the total of
$7$ pairs of non-$RH$ pairs for $\hat{H}^0$.
These $7$ pairs occur after $\varpi_0=6522.513197.$ The total
number of zeros here is $N_z^0=30$. 

This is actually unexpected because there is quite a regular
behavior of $S^0$ for the family $T(3m\pm 1,3)$. Namely,
the following formula for \Yboxdim5pt \yng(2,1)\,\Yboxdim7pt\ 
is likely to hold for this family:
$$
\varsigma_0 t^{\pi_0}S^0=-t^3\ \frac{1-t^{4(2n-3)}}{1-t^2}\ \, 
\frac{1-t^{2n}}{1-t^2}
\for T(n,3),\, n=3m\pm 1>2.
$$
It was checked in the examples above and for
$T(8,3)$.
It obviously collapses for $T(3,2)$, where the actually one is
$S^0=\frac{1-t^{6}}{1-t^2}$. 

Surprisingly, 
a similar formula for the $3$\~hook is more 
involved for $T(2m+1,2)$. Namely, non-cyclotomic 
factors occur in $S^0$ for $m\ge 3$. For instance for $T(11,2)$:\,
$\varsigma_0 t^{\pi_0}S^0\!=
${\small $-t^3\, \frac{1-t^{10}}{1-t^2}\,
(1+t^2+2 t^4+2 t^6+3 t^8+3 t^{10}+3 t^{12}+3 t^{14}
+4 t^{16}+3 t^{18}+4 t^{20}+3 t^{22}
+3 t^{24}+3 t^{26}+3 t^{28}+2 t^{30}+2 t^{32}+t^{34}+t^{36}).$}
There is some pattern here, but not too simple.



\subsubsection{\sf Multiple zeros}\label{SEC:MULTZ}
Let us provide $\varsigma_i t^{\pi_i}S^i$
for all links from the table colored by rectangles
where at least one $S^i$ has multiple zeros. We give them for 
all $0\le i\le$ deg$_{a} \hat{\h}$ (not only when multiple 
roots occur). 

\vskip 0.1cm

(1) \ \ $\{\ga[3,2](P(3+3)\}$, entry $58 (N_z^0=60) :$
{\footnotesize
\begin{align*}
&\varsigma_0 t^{\pi_0}S^0=\,\frac{1-t^{120}}{1-t^{60}},\ \ \ \ 
\varsigma_1 t^{\pi_1}S^1=t\frac{1-t^{109}}{1-t^{54}},\ \ \ \ 
\varsigma_2 t^{\pi_2}S^2=(1+t^2)(1+t^{48}),\\
&\varsigma_3 t^{\pi_3}S^3=t\frac{1-t^{84}}{1-t^{42}},\ \ \ \
\varsigma_4 t^{\pi_4}S^4=(1+t^2)(1+t^{34}), \\
&\varsigma_5 t^{\pi_5}S^5=t\frac{1-t^{52}}{1-t^{26}},\ \ \ \
\varsigma_6 t^{\pi_6}S^6=(1-t^{36})(1-t^{18}).
\end{align*}
}
The failure of Strong $RH$ at $i=1$ is due to the 
factor $t$. The multiple factor is $(1+t^2)$ at $i=4$, but it 
does not result in non-$RH$ zeros. 

\vskip 0.1cm
(2) \ \ $\{\ga[3,2](P(2+0)P(1+1)\}$, entry $63 (N_z^0=36) :$
{\footnotesize
\begin{align*}
&\varsigma_0 t^{\pi_0}S^0\!=\!-t^2(1\!+\!t^2)(1\!+\!t^{30}),\ \ 
\varsigma_1 t^{\pi_1}S^1\!=\!-t^2(1\!+\!t^2)(1\!+\!t^{28}),\ \ 
\varsigma_2 t^{\pi_2}S^2\!=\!-t^2\frac{1\!-\!t^{56}}{1\!-\!t^{28}},\\
&\varsigma_3 t^{\pi_3}S^3\!=\!-t^2\frac{1-t^{44}}{1-t^{22}},\ \ \ 
\varsigma_4 t^{\pi_4}S^4=-t^2\frac{1-t^{28}}{1-t^{14}},\ \ \ 
\varsigma_5 t^{\pi_5}S^5=-t^2\frac{1-t^{8}}{1-t^{4}}.
\end{align*}
}
The expected number of pairs of non-$RH$ zeros at $i=0$
is $4$: $2$ because of $t^2$ plus $2$ due to the multiple
$(1+t^2)$, which matches the table. Here
the $q$\~deformations of multiple $\pm \imath$ are {\em not\,} 
of the $RH$\~ type. 

\vskip 0.1cm
(3)\ \ $\{\ga[1,1](\ga[2,1](P)(1)\ga[3,2](P)(1))\}$, 
entry $37 (N_z^0=22) :$

{\footnotesize
\begin{align*}
&\varsigma_0 t^{\pi_0}S^0=-t\frac{1-t^{22}}{1-t^{2}},\ \ \ 
\varsigma_1 t^{\pi_1}S^1=-t\frac{1-t^{20}}{1-t^{2}},\ \ \
\varsigma_2 t^{\pi_2}S^2=-t\frac{1-t^{16}}{1-t^{2}},\\
&\varsigma_3 t^{\pi_3}S^3=-t(1+t^2)\frac{1-t^{8}}{1-t^{2}},\ \ \ \  
\varsigma_4 t^{\pi_4}S^4=1+t^{2}.
\end{align*}
}
Here the multiple factor is $(1+t^2)$ for $i=3$; however
Weak $RH$ holds (with one pair of real non-$RH$ zeros due
to $t=0$). The $q$\~deformations of zeros of $1+t^2=0$ are
(remain)  of the $RH$\~type.


\vskip 0.1cm
(4)\ \ $\{\ga[1,1](\ga[3,2](P)(1)\ga[3,2](P)(1))\}$, 
entry $64 (N_z^0=34) :$

{\footnotesize
\begin{align*}
&\varsigma_0 t^{\pi_0}S^0=-t\frac{1-t^{34}}{1-t^{2}},\ \ \
\varsigma_1 t^{\pi_1}S^1=-t\frac{1-t^{32}}{1-t^{2}},\ \ \
\varsigma_2 t^{\pi_2}S^2=-t\frac{1-t^{28}}{1-t^{2}},\\
&\varsigma_3 t^{\pi_3}S^3=-t\frac{1-t^{22}}{1-t^{2}},\ \ 
\varsigma_4 t^{\pi_4}S^4=-t(1+t^2)^2(1+t^4)^2, \ \ 
\varsigma_5 t^{\pi_5}S^5=(1+t^2)^2, 
\end{align*}
}
Here Weak $RH$ fails at $i=4$ ($N_z^4=14$) with
$5$ pairs of nonzero non-$RH$ zeros; the expected number
is $7$, i.e. {\em all\,} of them (including the contribution
of $t=0$). However the $q$\~deformations of the roots of
$(1+t^2)=0$ remain of $RH$\~type. Similarly, 
Weak $RH$ holds for $i=5$.

\comment{
{"super2-3+3-3", 4, 1, cycle, 0, 1 + t^2 + t^34 + t^36}

{"super6-4+2-0&1-1", 0, 1, cycle, 2, -t^2 - t^4 - t^32 - t^34}

{"super4-2+1-0&&1-0!+1-1", 3, 1, cycle, 1, -t - 2*t^3 - 
2*t^5 - 2*t^7 - t^9}

{"super6-4+1-0&&1-0!+1-1", 4, 1, cycle, 1, -t - 2*t^3 - 
3*t^5 - 4*t^7 - 3*t^9 - 2*t^11 - t^13}

{"super6-4+1-0&&1-0!+1-1", 5, 1, cycle, 0, 1 + 2*t^2 + t^4}
}

\subsection{\bf The table}
It is focused on the validity of Weak $RH$ for all
$0\!\le\! i\!\le\!$ deg$_a$, especially in the case $i\!=\!0$.
See below some review of our calculations presented
in the table (and beyond). 
\comment{
If ``{\sl HOLDS\,}" is not the case, but $RH$
is still valid for $i=0$, we consider $i=1,i=2$ and so on
till the last $i$, when $RH$ holds. We put
``{\sl OK $a\!=\!0\,$}" in this case. The uncolored
case $a=0$ is actually
important for $RH$; it is directly connected
with zeta-functions of singular curves.
If $RH$ fails even at $a=0$, then we put ``{\sl FAILS\,}".

Actually, there is only one
$OK$ (a failure of Weak $RH$ for $i\!>\!0$)
among uncolored algebraic links in the table, which is
at $i=4$ for entry No=$64$ (with $N^0_z=34$).
This limits our expectations for the validity of $RH$ for
flags ($i\!>\!0$) in the case of (algebraic)
links. 

It is not impossible 
that any DAHA superpolynomials for {\em colored\,}
 algebraic links are 
of geometric origin, but we do not know how far 
this can go. 
Anyway, $RH$ with colors and flags is not actually connected to 
the classical zeta-functions, even heuristically.
}

{\footnotesize
\begin{longtable}{|c|c|c|c|c|c|} 
 \hline 
No ($N^0_z$) & GOOD & ALG & DAHA-formula & 
$\varpi_0\,$ & $N^0_{irz}$ \\ 
RH-type  & $deg_a\!+\!1$ & type & CABLE(basic) & 
$\varpi_{1,\ldots}^{max}$ 
& $N^{max}_{irz}$ \\ 
\hline 
\hline 
1 (12) & 0$\le\!i\!\le$2 & alg & \{$\ga$[7,3](P)\} & 0.95272 & 0-0 \\ 
HOLDS & {\tiny all= }3 & knot & T(7,3) & 1. & 0-0 \\ 
\hline 
2 (14) & 0$\le\!i\!\le$2 & alg & \{$\ga$[8,3](P)\} & 0.96465 & 0-0 \\ 
HOLDS & {\tiny all= }3 & knot & T(8,3) & 0.962706 & 0-0 \\ 
\hline 
3 (18) & 0$\le\!i\!\le$2 & alg & \{$\ga$[10,3](P)\} & 0.980586 & 0-0 \\ 
HOLDS & {\tiny all= }3 & knot & T(10,3) & 0.95272 & 0-0 \\ 
\hline 
4 (20) & 0$\le\!i\!\le$2 & alg & \{$\ga$[11,3](P)\} & 0.984635 & 0-0 \\ 
HOLDS & {\tiny all= }3 & knot & T(11,3) & 1. & 0-0 \\ 
\hline 
5 (360) & 0$\le\!i\!\le$2 & alg & \{$\ga$[181,3](P)\} & 0.999995 & 0-0 \\ 
HOLDS & {\tiny all= }3 & knot & T(181,3) & 1. & 0-0 \\ 
\hline 
6 (366) & 0$\le\!i\!\le$2 & alg & \{$\ga$[184,3](P)\} & 0.999995 & 0-0 \\ 
HOLDS & {\tiny all= }3 & knot & T(184,3) & 0.999995 & 0-0 \\ 
\hline 
7 (16) & 0$\le\!i\!\le$3 & alg & \{$\ga$[3,2]$\ga$[2,1](P)\} & 1.49558 & 0-0 \\ 
HOLDS & {\tiny all= }4 & knot & Cab(13,2)T(3,2) & 2.17649 & 0-0 \\ 
\hline 
8 (12) & 0$\le\!i\!\le$3 & alg & \{$\ga$[5,4](P)\} & 0.945441 & 0-0 \\ 
HOLDS & {\tiny all= }4 & knot & T(5,4) & 1. & 0-0 \\ 
\hline 
9 (24) & 0$\le\!i\!\le$3 & alg & \{$\ga$[9,4](P)\} & 0.987206 & 0-0 \\ 
HOLDS & {\tiny all= }4 & knot & T(9,4) & 0.978578 & 0-0 \\ 
\hline 
10 (32) & 0$\le\!i\!\le$4 & alg & \{$\ga$[5,2](P(2+2))\} & 2.17837 & 0-0 \\ 
HOLDS & {\tiny all= }5 & knot & T(5,2) & 3.72664 & 0-0 \\ 
\hline 
11 (36) & 0$\le\!i\!\le$4 & alg & \{$\ga$[4,3](P(2+0))\} & 1.69034 & 0-0 \\ 
HOLDS & {\tiny all= }5 & knot & T(4,3) & 1.79669 & 0-0 \\ 
\hline 
12 (20) & 0$\le\!i\!\le$4 & alg & \{$\ga$[6,5](P)\} & 0.98107 & 0-0 \\ 
HOLDS & {\tiny all= }5 & knot & T(6,5) & 1. & 0-0 \\ 
\hline 
13 (32) & 0$\le\!i\!\le$4 & alg & \{$\ga$[9,5](P)\} & 0.993845 & 0-0 \\ 
HOLDS & {\tiny all= }5 & knot & T(9,5) & 1. & 0-0 \\ 
\hline 
14 (74) & 0$\le\!i\!\le$5 & alg & \{$\ga$[5,2]$\ga$[3,2](P)\} & 1.36037 & 0-0 \\ 
HOLDS & {\tiny all= }6 & knot & Cab(32,3)T(5,2) & 1.57212 & 0-0 \\ 
\hline 
15 (36) & 0$\le\!i\!\le$5 & alg & \{$\ga$[4,3]$\ga$[2,1](P)\} & 1.42228 & 0-0 \\ 
HOLDS & {\tiny all= }6 & knot & Cab(25,2)T(4,3) & 1.6464 & 0-0 \\ 
\hline 
16 (30) & 0$\le\!i\!\le$5 & alg & \{$\ga$[7,6](P)\} & 0.993335 & 0-0 \\ 
HOLDS & {\tiny all= }6 & knot & T(7,6) & 1. & 0-0 \\ 
\hline 
17 (50) & 0$\le\!i\!\le$5 & alg & \{$\ga$[11,6](P)\} & 1.14474 & 0-0 \\ 
HOLDS & {\tiny all= }6 & knot & T(11,6) & 1.13811 & 0-0 \\ 
\hline 
18 (96) & 0$\le\!i\!\le$6 & alg & \{$\ga$[3,2]$\ga$[2,1](P(2+0))\} & 1.49797 & 0-0 \\ 
HOLDS & {\tiny all= }7 & knot & Cab(13,2)T(3,2) & 1.65491 & 0-0 \\ 
\hline 
19 (108) & 0$\le\!i\!\le$6 & alg & \{$\ga$[3,2]$\ga$[2,3](P(2+0))\} & 1.45667 & 0-0 \\ 
HOLDS & {\tiny all= }7 & knot & Cab(15,2)T(3,2) & 1.56317 & 0-0 \\ 
\hline 
20 (66) & 0$\le\!i\!\le$6 & alg & \{$\ga$[12,7](P)\} & 1.119 & 0-0 \\ 
HOLDS & {\tiny all= }7 & knot & T(12,7) & 1.11544 & 0-0 \\ 
\hline 
21 (72) & 0$\le\!i\!\le$6 & alg & \{$\ga$[13,7](P)\} & 1.16157 & 0-0 \\ 
HOLDS & {\tiny all= }7 & knot & T(13,7) & 1.17036 & 0-0 \\ 
\hline 
22 (84) & 0$\le\!i\!\le$7 & alg & \{$\ga$[3,2]$\ga$[2,1]$\ga$[2,1](P)\} & 1.4672 & 0-0 \\ 
HOLDS & {\tiny all= }8 & knot & Cab(53,2)Cab(13,2)T(3,2) & 1.56196 & 0-0 \\ 
\hline 
23 (80) & 0$\le\!i\!\le$7 & alg & \{$\ga$[3,2]$\ga$[4,1](P)\} & 1.37538 & 0-0 \\ 
HOLDS & {\tiny all= }8 & knot & Cab(25,4)T(3,2) & 1.46552 & 0-0 \\ 
\hline 
24 (66) & 0$\le\!i\!\le$7 & alg & \{$\ga$[5,4]$\ga$[2,3](P)\} & 1.51732 & 0-0 \\ 
HOLDS & {\tiny all= }8 & knot & Cab(43,2)T(5,4) & 1.58227 & 0-0 \\ 
\hline 
25 (48) & 0$\le\!i\!\le$8 & alg & \{$\ga$[4,3](P(2+2))\} & 2.73447 & 0-0 \\ 
HOLDS & {\tiny all= }9 & knot & T(4,3) & 6.01964 & 0-0 \\ 
\hline 
26 (90) & 0$\le\!i\!\le$8 & alg & \{$\ga$[4,3]$\ga$[3,1](P)\} & 1.31174 & 0-0 \\ 
HOLDS & {\tiny all= }9 & knot & Cab(37,3)T(4,3) & 1.42245 & 0-0 \\ 
\hline 
27 (64) & 0$\le\!i\!\le$8 & alg & \{$\ga$[5,3](P(2+2))\} & 2.1449 & 0-0 \\ 
HOLDS & {\tiny all= }9 & knot & T(5,3) & 21.1022 & 0-0 \\ 
\hline 
28 (116) & 0$\le\!i\!\le$8 & alg & \{$\ga$[5,3]$\ga$[3,2](P)\} & 1.31641 & 0-0 \\ 
HOLDS & {\tiny all= }9 & knot & Cab(47,3)T(5,3) & 1.57212 & 0-0 \\ 
\hline 
29 (72) & 0$\le\!i\!\le$8 & alg & \{$\ga$[10,9](P)\} & 1.19316 & 0-0 \\ 
HOLDS & {\tiny all= }9 & knot & T(10,9) & 1.21884 & 0-0 \\ 
\hline 
30 (80) & 0$\le\!i\!\le$8 & alg & \{$\ga$[11,9](P)\} & 1.20386 & 0-0 \\ 
HOLDS & {\tiny all= }9 & knot & T(11,9) & 1.26173 & 0-0 \\ 
\hline 
31 (96) & 0$\le\!i\!\le$8 & alg & \{$\ga$[13,9](P)\} & 1.18157 & 0-0 \\ 
HOLDS & {\tiny all= }9 & knot & T(13,9) & 1.20367 & 0-0 \\ 
\hline 
32 (6) & 0$\le\!i\!\le$1 & alg & \{$\ga$[3,1](PP)\} & 0 & 0-1 \\ 
HOLDS & {\tiny all= }2 & link & T(3,1) & 0 & 0-1 \\ 
\hline 
33 (24) & 0$\le\!i\!\le$3 & alg & \{$\ga$[1,1]($\ga$[1,1](PPPP)(1))\} & 1.87601 & 3-3 \\ 
HOLDS & {\tiny all= }4 & link & (4Cab(1,1))(T(1,1)) & 1.89877 & 3-3 \\ 
\hline 
34 (28) & 0$\le\!i\!\le$3 & alg & \{$\ga$[5,2](PP)\} & 1.60724 & 1-1 \\ 
HOLDS & {\tiny all= }4 & link & T(5,2) & 1.94384 & 1-1 \\ 
\hline 
35 (26) & 0$\le\!i\!\le$3 & alg & \{$\ga$[2,1](P $\ga$[3,2](P)(1))\} & 1.21155 & 0-1 \\ 
HOLDS & {\tiny all= }4 & link & T(2,1) & 1.54031 & 0-1 \\ 
\hline 
36 (52) & 0$\le\!i\!\le$4 & alg & \{$\ga$[3,2](P(2+0)P)\} & 1.61361 & 1-1 \\ 
HOLDS & {\tiny all= }5 & link & T(3,2) & 1.81723 & 1-1 \\ 
\hline 
37 (22) & 0$\le\!i\!\le$4 & alg & \{$\ga$[1,1]($\ga$[2,1](P)(1)$\ga$[3,2](P)(1))\} & 1.91393 & 1-1 \\ 
HOLDS & {\tiny all= }5 & link & (Cab(5,3),Cab(3,2))(T(1,1)) & 2.00929 & 1-1 \\ 
\hline 
38 (42) & 0$\le\!i\!\le$4 & alg & \{$\ga$[2,1]($\ga$[2,1](P)(1)$\ga$[3,2](P)(1))\} & 1.37608 & 1-1 \\ 
HOLDS & {\tiny all= }5 & link & (Cab(7,3),Cab(5,2))(T(2,1)) & 1.91393 & 1-1 \\ 
\hline 
39 (24) & 0$\le\!i\!\le$4 & alg & \{P(Y)($\ga$[2,3](PP)(1))\} & 0.736757 & 2-2 \\ 
HOLDS & {\tiny all= }5 & link & T(2,3) & 1.2963 & 2-2 \\ 
\hline 
40 (46) & 0$\le\!i\!\le$5 & alg & \{$\ga$[1,1]$\ga$[3,2](PP)\} & 1.61257 & 1-1 \\ 
HOLDS & {\tiny all= }6 & link & Cab(5,3)T(1,1) & 1.94384 & 1-1 \\ 
\hline 
41 (76) & 0$\le\!i\!\le$5 & alg & \{$\ga$[2,1]$\ga$[3,2](PP)\} & 1.4269 & 1-1 \\ 
HOLDS & {\tiny all= }6 & link & Cab(8,3)T(2,1) & 1.61257 & 1-1 \\ 
\hline 
42 (42) & 0$\le\!i\!\le$5 & alg & \{$\ga$[3,2](PPP)\} & 1.59651 & 2-2 \\ 
HOLDS & {\tiny all= }6 & link & 3T(3,2) & 1.99525 & 2-2 \\ 
\hline 
43 (36) & 0$\le\!i\!\le$5 & alg & \{$\ga$[4,3](PP)\} & 1.69177 & 1-1 \\ 
HOLDS & {\tiny all= }6 & link & T(4,3) & 1.72808 & 1-1 \\ 
\hline 
44 (64) & 0$\le\!i\!\le$5 & alg & \{$\ga$[2,1]($\ga$[3,2](P)(1)$\ga$[3,2](P)(1))\} & 1.55923 & 1-1 \\ 
HOLDS & {\tiny all= }6 & link & (Cab(8,3),Cab(8,3))(T(2,1)) & 2.18422 & 1-1 \\ 
\hline 
45 (94) & 0$\le\!i\!\le$5 & alg & \{$\ga$[3,1]($\ga$[3,2](P)(1)$\ga$[3,2](P)(1))\} & 1.39031 & 1-1 \\ 
HOLDS & {\tiny all= }6 & link & (Cab(11,3),Cab(11,3))(T(3,1)) & 1.55923 & 1-1 \\ 
\hline 
46 (72) & 0$\le\!i\!\le$5 & alg & \{$\ga$[5,2](P $\ga$[2,1](P))\} & 1.39868 & 1-1 \\ 
HOLDS & {\tiny all= }6 & link & Cab(21,2)T(5,2) & 1.60752 & 1-1 \\ 
\hline 
47 (100) & 0$\le\!i\!\le$6 & alg & \{P(Y)($\ga$[2,1]$\ga$[3,2](PP)(1))\} & 1.35321 & 2-2 \\ 
HOLDS & {\tiny all= }7 & link & Cab(8,3)T(2,1) & 1.52945 & 2-2 \\ 
\hline 
48 (80) & 0$\le\!i\!\le$7 & alg & \{$\ga$[3,2](PPPP)\} & 1.59986 & 3-3 \\ 
HOLDS & {\tiny all= }8 & link & 4T(3,2) & 1.912 & 3-3 \\ 
\hline 
49 (84) & 0$\le\!i\!\le$7 & alg & \{$\ga$[3,2]$\ga$[2,1](PP)\} & 1.47591 & 1-1 \\ 
HOLDS & {\tiny all= }8 & link & Cab(13,2)T(3,2) & 1.60724 & 1-1 \\ 
\hline 
50 (62) & 0$\le\!i\!\le$7 & alg & \{$\ga$[1,1]($\ga$[5,2](P)(1)$\ga$[3,2](P)(1))\} & 1.59801 & 1-1 \\ 
HOLDS & {\tiny all= }8 & link & (Cab(7,5),Cab(5,3))(T(1,1)) & 2.06334 & 1-1 \\ 
\hline 
51 (82) & 0$\le\!i\!\le$7 & alg & \{$\ga$[3,2](P $\ga$[3,2](P))\} & 1.35745 & 1-1 \\ 
HOLDS & {\tiny all= }8 & link & Cab(20,3)T(3,2) & 1.59033 & 1-1 \\ 
\hline 
52 (90) & 0$\le\!i\!\le$8 & alg & \{$\ga$[4,3](PPP)\} & 1.55728 & 2-2 \\ 
HOLDS & {\tiny all= }9 & link & 3T(4,3) & 1.6737 & 2-2 \\ 
\hline 
53 (10) & 0$\le\!i\!\le$2 & alg\kern-13pt{\bf=} & \{$\ga$[3,2]$\ga$[-2,5](P)\} & 1.8969 & 0-0 \\ 
HOLDS & {\tiny all= }3 & knot & Cab(7,2)T(3,2) & 0.75 & 0-0 \\ 
\hline 
54 (8) & 0$\le\!i\!\le$3 & alg\kern-13pt{\bf=} & \{$\ga$[3,2]$\ga$[-2,7](P)\} & 0.88617 & 0-0 \\ 
HOLDS & {\tiny all= }4 & knot & Cab(5,2)T(3,2) & 0.5 & 0-0 \\ 
\hline 
55 (14) & 0$\le\!i\!\le$3 & alg\kern-13pt{\bf=} & \{$\ga$[3,2]$\ga$[-2,1](P)\} & 1.42745 & 0-0 \\ 
HOLDS & {\tiny all= }4 & knot & Cab(11,2)T(3,2) & 1.88749 & 0-0 \\ 
\hline 
56 (26) & 0$\le\!i\!\le$5 & alg\kern-13pt{\bf=} & \{$\ga$[4,3]$\ga$[2,-9](P)\} & 1.6649 & 0-0 \\ 
HOLDS & {\tiny all= }6 & knot & Cab(15,2)T(4,3) & 1.61525 & 0-0 \\ 
\hline 
\hline 
57 (20) & 0$\le\!i\!\le$-1 & alg & \{$\ga$[5,2](P(2+1))\} & -- & --5 \\ 
FAILS & {\tiny all= }4 & knot & T(5,2) & -- & --- \\ 
\hline 
58 (60) & 0$\le\!i\!\le$0 & alg & \{$\ga$[3,2](P(3+3))\} & 2.15141 & 0-0 \\ 
OK {\tiny$a\!=\!\!0$} & {\tiny all= }7 & knot & T(3,2) & -- & --- \\ 
\hline 
59 (30) & 0$\le\!i\!\le$-1 & alg & \{$\ga$[4,3](P(2+1))\} & -- & --5 \\ 
FAILS & {\tiny all= }7 & knot & T(4,3) & -- & --- \\ 
\hline 
60 (40) & 0$\le\!i\!\le$-1 & alg & \{$\ga$[5,3](P(2+1))\} & -- & --3 \\ 
FAILS & {\tiny all= }7 & knot & T(5,3) & -- & --- \\ 
\hline 
61 (60) & 0$\le\!i\!\le$-1 & alg & \{$\ga$[7,3](P(2+1))\} & -- & --3 \\ 
FAILS & {\tiny all= }7 & knot & T(7,3) & -- & --- \\ 
\hline 
62 (8) & 0$\le\!i\!\le$-1 & alg & \{$\ga$[2,1](P(2+2)P)\} & -- & --2 \\ 
FAILS & {\tiny all= }2 & link & T(2,1) & -- & --- \\ 
\hline 
63 (36) & 0$\le\!i\!\le$-1 & alg & \{$\ga$[3,2](P(2+0)P(1+1))\} & -- & --4 \\ 
FAILS & {\tiny all= }6 & link & T(3,2) & -- & --- \\ 
\hline 
64 (34) & 0$\le\!i\!\le$3 & alg & \{$\ga$[1,1]($\ga$[3,2](P)(1)$\ga$[3,2](P)(1))\} & 2.18422 & 1-1 \\ 
OK {\tiny$a\!=\!\!0$} & {\tiny all= }6 & link & (Cab(5,3),Cab(5,3))(T(1,1)) & 2.61867 & 1-1 \\ 
\hline 
65 (42) & 0$\le\!i\!\le$-1 & alg & \{$\ga$[3,2](P(2+2)P)\} & -- & --5 \\ 
FAILS & {\tiny all= }7 & link & T(3,2) & -- & --- \\ 
\hline 
66 (6) & 0$\le\!i\!\le$-1 & alg\kern-13pt{\bf=} & \{$\ga$[3,2]$\ga$[-2,13](P)\} & -- & --2 \\ 
FAILS & {\tiny all= }4 & knot & Cab(1,-2)T(3,2) & -- & --- \\ 
\hline 
67 (8) & 0$\le\!i\!\le$-1 & alg\kern-13pt{\bf=} & \{$\ga$[3,2]$\ga$[-2,15](P)\} & -- & --4 \\ 
FAILS & {\tiny all= }4 & knot & Cab(3,-2)T(3,2) & -- & --- \\ 
\hline 
68 (10) & 0$\le\!i\!\le$-1 & alg\kern-13pt{\bf=} & \{$\ga$[3,2]$\ga$[-2,17](P)\} & -- & --3 \\ 
FAILS & {\tiny all= }4 & knot & Cab(5,-2)T(3,2) & -- & --- \\ 
\hline 
69 (12) & 0$\le\!i\!\le$-1 & alg\kern-13pt{\bf=} & \{$\ga$[3,2]$\ga$[-2,19](P)\} & -- & --6 \\ 
FAILS & {\tiny all= }4 & knot & Cab(7,-2)T(3,2) & -- & --- \\ 
\hline 
70 (4) & 0$\le\!i\!\le$0 & alg\kern-13pt{\bf=} & \{$\ga$[3,2]$\ga$[-2,11](P)\} & 1.33333 & 0-0 \\ 
OK {\tiny$a\!=\!\!0$} & {\tiny all= }4 & knot & Cab(1,2)T(3,2) & -- & --- \\ 
\hline 
71 (22) & 0$\le\!i\!\le$-1 & alg\kern-13pt{\bf=} & \{$\ga$[4,3]$\ga$[-2,25](P)\} & -- & --6 \\ 
FAILS & {\tiny all= }6 & knot & Cab(-1,2)T(4,3) & -- & --- \\ 
\hline 
72 (20) & 0$\le\!i\!\le$-1 & alg\kern-13pt{\bf=} & \{$\ga$[4,3]$\ga$[-2,23](P)\} & -- & --4 \\ 
FAILS & {\tiny all= }6 & knot & Cab(1,2)T(4,3) & -- & --- \\ 
\hline 
73 (18) & 0$\le\!i\!\le$-1 & alg\kern-13pt{\bf=} & \{$\ga$[4,3]$\ga$[-2,21](P)\} & -- & --2 \\ 
FAILS & {\tiny all= }6 & knot & Cab(3,2)T(4,3) & -- & --- \\ 
\hline 
74 (16) & 0$\le\!i\!\le$0 & alg\kern-13pt{\bf=} & \{$\ga$[4,3]$\ga$[-2,19](P)\} & 1.29137 & 0-0 \\ 
OK {\tiny$a\!=\!\!0$} & {\tiny all= }6 & knot & Cab(5,2)T(4,3) & -- & --- \\ 
\hline 
75 (18) & 0$\le\!i\!\le$1 & alg\kern-13pt{\bf=} & \{$\ga$[4,3]$\ga$[-2,17](P)\} & 1.06963 & 0-0 \\ 
OK {\tiny$a\!=\!\!0$} & {\tiny all= }6 & knot & Cab(7,2)T(4,3) & 1.38291 & 0-0 \\ 
\hline 
76 (22) & 0$\le\!i\!\le$3 & alg\kern-13pt{\bf=} & \{$\ga$[4,3]$\ga$[-2,13](P)\} & 1.69298 & 0-0 \\ 
OK {\tiny$a\!=\!\!0$} & {\tiny all= }6 & knot & Cab(11,2)T(4,3) & 1.51094 & 0-0 \\ 
\hline 
\end{longtable} 

}

\vskip -1cm
\subsection{\bf Brief analysis}
We will focus here on the most instructional cases and
features. Only a few formulas for DAHA
superpolynomials are provided in this work; the
files are available and see our prior works, which
contain many.


\subsubsection{\sf Non-algebraic links}
The validity of $RH$ is certainly
most likely for {\em algebraic\,} knots/links. However,
the table and calculations we performed show
that it holds significantly beyond this class. For instance,
the whole series of cables $C\!ab(2m+1,2)T(3,2)$ from
$C\!ab(13,2)T(3,2)$ (which is the smallest algebraic)
down to $C\!ab(3,2)T(3,2)$ satisfies $RH$ for all $i\ge 0$.
If $a=0$ and one continues to diminish  $m$, then
$RH$ fails at $C\!ab(-1,2)T(3,2)$ (and further on).

We note that the DAHA superpolynomials remain positive
in this series till $C\!ab(7,2)T(3,2)$ (all
their $q,t,a$\~coefficients are positive).  We called them
{\em pseudo-algebraic\,} in \cite{ChD1}, where this series
was considered in detail.  They do resemble
algebraic knots (especially for $2m+1=11,9$).
Interestingly, the positivity of the superpolynomial
recovers starting with
$C\!ab(-7,2)T(3,2)$, but with a different pattern
of superpolynomials  (not like those for $m\ge 4$).
However, generally,
$RH$ has a clear tendency to
fail when ``approaching" non-algebraic links.

Let us consider
entry No$=64$ with $34$ zeros (at $i\!=\!0$), already discussed
above. This is an
{\em algebraic link\,} where $RH$ is valid only
partially (till $a^3$, where deg$_a=5$). The outer $\ga$ is
$\ga[1,1]$ here, i.e. minimal possible to make it algebraic.
Using $\ga[2,1]$ and $\ga[3,1]$ here instead of $\ga[1,1]$
results in $RH$ for {\em all\,} $a^i$
(the entries are $44,45$ with $N_z^0=64,94$). This can be
informally considered
as ``moving away" from non-algebraic ones. By the way,
uncolored $C\!ab(13,2)T(3,2)$, the smallest non-torus cable,
has $\varpi_3>2$, the only uncolored algebraic knot in the table
with $\varpi_i>2$ for $i>0$; we conjecture that
it is below $2$ for $i=0$.

\subsubsection{\sf Non-square diagrams} 
Strong $RH$ always fails for
\Yboxdim5pt \yng(2,1)\,\Yboxdim7pt\,
in the examples we calculated. This diagram is self-dual with
respect to transposition, which is valuable to us.
See entries $57,59,60,61$ (with $N_z^0=20,30,40,60$)
and Section \ref{SEC:2-1}.

For the entry $60$, there are $3$ pair of irregular
zeros of $\hat{H}^0$ (of $t$\~degree $40$)
that tend to $\{\om^{1/3},\om^{2/3}\}$ for the
$3$ different values of the cube root as $\om\to\infty$.
This is directly
related to the expansion $\hat{H}^0= \om+2t^2-t^3+O(\om^2)$,
so a variant of Weak $RH$ holds in this case. 

In the remaining $3$ cases, there are two more complex
(conjugated to each other) irregular pairs.  See the last column 
of the table, which provides the (actual) number of pairs of 
irregular zeros. Recall that we provide this number only for 
$i=0$ for entries with  "{\em FAILS\,}". 
These additional pairs $\{\xi',\om/\xi'\}$
stay in the vicinity of
$U_{\sqrt{\om}}$. More exactly,
$\xi'=\pm(\sqrt{\om}+C)\imath+(C+o(1))/\sqrt{\om}$, 
where $C\in \R$ 
tends to some limit as $\om\to \infty$. 
Similar $\xi'$ occur for entries
$62,63,65$ (with $N_z^0=8,36,42$).

\comment{
We note that 
$\om_0^{top}$ (the top real root of the reduced
discriminant $D^0$) can be large for
$\Yboxdim5pt \yng(2,1)\,$;\, for instance,
$\om_0^{top}\!=\!190,\, 335\,$ for entries
 $60,61$ (with $N_z^0\!=\!40,60$).}
\Yboxdim5pt
\vskip 0.2cm

For {\em non-square\,}
rectangles,  $RH$ for $\hat{H}^i_{sym}$ can hold 
beyond $i=0$, where this
was conjectured. For instance, it holds for any $i$
for entries $18,19$ ($N_z^0=96,108$):\,
$\{\hat{\ga}_{3,2}(\hat{\ga}_{2,1}(P(\yng(1,1)\,)))\}$,
$\{\hat{\ga}_{3,2}(\hat{\ga}_{2,3}(P(\yng(1,1)\,)))\}$.

For $\{\hat{\ga}_{3,2}(P(\yng(2,2,2)\,))\}$,
only one  {\em non-RH\,} pair appears for
entry $58$ ($N_z^0=60$, deg$_a=6$) when $i=1$;
it is real and quickly approaches $\{-1,-\om\}$ as 
$\om\to \infty$. See Section \ref{SEC:MULTZ}.
$RH$ also fails here for $i=1,3,5$ upon the
{\em antisymmetrization\,}; the corresponding 
$\varsigma_i t^{\pi_i}S^i$ are $-t(1+t^{54,42,26})$.

For the two $2$\~rectangle and $2T(3,2)$, i.e. for
$\{\hat{\ga}_{3,2}(P(\yng(1,1)\,)^2)\}$, there are
$2$ real pairs of such zeros when $i=0$
(from $48$ pairs)
approaching $\{1,\om\}$
and $\{\pm\sqrt{\om}\}$. One has:
$\varsigma_0 t^{\pi_0}S^0=-t\frac{1-t^{92}}{1-t^4}(1-t^2)^2(1+t^2)$
in this case. The multiple roots $\pm 1$ 
and those due to $t$ give the total number of pairs
of non-$RH$ zeros $3$ (all are real). 
Upon the antisymmetrization here,  
$\varsigma_0 t^{\pi_0}S^0$ becomes
$t(1-t^2)(1+t^{92})$; Weak $RH$ is satisfied now.

Also,  $RH$ holds for 
$\{\hat{\ga}_{3,2}(P(\yng(1,1)\,)P(\yng(1))\}$
upon the symmetrization for all $0\le i\le 4$. One 
respectively has
$\varsigma_i t^{\pi_i}S^i=-t(1+t^{50,46,40,32,22})$.

\comment{
Anyway using colors more general than columns and the
symmetrization is not supported geometrically at
the moment;
 we do not know which ``Jacobian factors"
can be associated with them (if any) beyond \cite{ChP2}.}

\comment{
Using various rectangles for links is not sufficiently supported
numerically. The failure at entry No$=\!\!63$ ($N_z^0\!=\!36$) for
$\{\hat{\ga}_{3,2}(P(\yng(1,1)\,)P(\yng(2)\,))\}$
is interesting. No symmetrization is necessary here,
since the resulting superpolynomial satisfies super-duality.
$RH$ fails in this example (even for $a=0$), as well as for
$\{\hat{\ga}_{2,1}(P(\yng(1,1)\,)P(\yng(2)\,)\}$=
$\{\hat{\ga}_{2,1}(P(\yng(2,2)\,)P(\yng(1)\,))\}$, which
is an interesting (but special) identity.
Entry $65$ ($N_z^0\!=\!42$) does not satisfy $RH$ too;  
the corresponding DAHA line is 
$\{\hat{\ga}_{3,2}(P(\yng(2,2)\,)P(\yng(1)\,))\}$.
With $\hat{\ga}_{2,1}$ here, $RH$ fails as well (entry $62$ with
$N_z^0\!=\!8$).
}
\vskip 0.2cm

\subsubsection{\sf Large-small 
\texorpdfstring{{\mathversion{bold}$\varpi$}}{omega}}
When $RH$ holds, which is a {\em qualitative\,} property,
the ``actual" $RH$ is in finding $\varpi$ {\em quantitatively\,}. 
In the following example of
$2\times 2$\~diagram, $RH$ ``almost" fails.
The  $\varpi$\~number becomes much larger than ``usual"
$1\sim 2$ due to the color and $i>0$.
This is entry No=$27$ with $64$ zeros and
deg$_a=8$ for  $T(5,3)$ colored by
\Yboxdim5pt \yng(2,2)\Yboxdim7pt\ .
The corresponding $\varpi_4$ is $21.1022$; it is
$\om_4^{top}$ for $D^4$ (and can be calculated as exactly as
necessary). By contrast, $\varpi_0=2.1449$, not too large.

This is similar to entry $25$ (with $48$ zeros). The latter is for
$T(4,3)$ and \Yboxdim5pt \yng(2,2)\Yboxdim7pt\ . Its 
superpolynomial is
significantly simpler to calculate than that for $T(5,3)$ (which 
took about $4$ days). One has $\varpi_4=6.01964$ for $T(4,3)$,
which is large but not too much.

\vskip 0.2cm

{\sf The inequality $\varpi<2$.} This seems a counterpart
of the bounds for the spectrum of Laplace-Dirac operators
in the theory of {\em spectral zeta-functions}.
We conjecture that $\varpi_0<2$ upon  
$a=0$ for {\em uncolored\,} algebraic knots. For $i>0$,
entry $7$ (with
$16$ zeros) is an example when $\varpi_1>2$. This is for
$C\!ab(13,2)T(3,2)$. Concerning algebraic uncolored {\em links\,}
for $i=0$, the
greatest $\varpi_0$ in the table is $1.919393$ for the entry
$37$ (with $N_z^0=22$), but it can go beyond $2$.

{\sf Cyclotomic polynomials at $\om=1$.} 
Let us discuss  $T(11,6)$; its superpolynomial was posted
in the online version of \cite{CJ}. See entry $17$ 
($N_z^0=50$); strong $RH$ for $i\!=\!0$
holds for $\varpi_0=1.1447417735112874\ldots$ and for 
$0.9985190700739621\!<\!\om\!<\!1.0021178996517260.$  
We expect that $\varpi_0=1+O(1/\rr)$ for {\em knots}
$T(\rr\!\to\!\infty,\ss\!>\!2)$, which is based on numerical
evidence for $\ss=3,6,7,9$ and the considerations
of {\em geometric superpolynomials} 
\cite{ChP1} as $\rr\to\infty$.   
One has:
$\hat{H}^0(\om\!=\!1)=\prod_{j=1}^{\ss-1}\frac{1-t^{\rr+\ss-j}}
{1-t^{1+j}}$ for any $T(\rr,\ss)$,  which follows from the 
Shuffle Conjecture proved in \cite{CaM}.  

This can be deduced from Corollary \ref{SHUFF};
here $\hat{H}^0(\om\!=\!1)$ is 
the {\em rational slope $t$\~Catalan number\,}. 
The following formula can be checked for
{\em geometric superpolynomials\,},
conjecturally coinciding with DAHA ones, and the 
$a$\~version of Corollary \ref{SHUFF}.

\comment{
$T(rrr,sss), rrr>sss here; zzzx is the deg=da component$.
\begin{align*}
zzzx=\prod_{ib=1}^{sss-1}
\frac{(1-t^{rrr+sss-ib})}{(1-t^{sss-ib+1})}\times
\frac{(1-t^{rrr-1})(1-t^{rrr-2})\ldots(1-t^{rrr-da})*
(1-t^{sss-1})(1-t^{sss-2})\ldots(1-t^{sss-da})}
{(1-t)(1-t^2)\ldots(1-t^da)
(1-t^{rrr+sss-1})(1-t^{rrr+sss-2})\ldots(1-t^{rrr+sss-da})}.
\end{align*}
}

\vskip -0.2cm
{\small
\begin{align*}
&\ \ \ \ \ \ \ \hbox{For \,} \rr>\ss,\ 0\le i \le \ss-1, \ \ 
\hbox{\normalsize $\hat{H}^i(\om=1)$}=
\prod_{j=1}^{\ss-1}
\frac{(1-t^{\rr+\ss-j})}{(1-t^{1+\ss-j})}\times\\
&\frac{\bigl((1-t^{\rr-1})(1-t^{\rr-2})\ldots(1-t^{\rr-i})\bigr)
\bigl((1-t^{\ss-1})(1-t^{\ss-2})\ldots(1-t^{\ss-i})\bigr)}
{\bigl((1-t)(1-t^2)\ldots(1-t^i)\bigr)\ \,
\bigl((1-t^{\rr+\ss-1})(1-t^{\rr+\ss-2})\ldots
(1-t^{\rr+\ss-i})\bigr)}\,.  
\end{align*}
}

The second line here is trivial when $i=0$.
Similar formulas exist for 
(at least) the family
$C\!ab(13+2m,2)T(3,2)$ with $m\in\Z_+$, where all 
Piontkowski cells $J_\De$ are affine spaces.

 
We note that $\hat{H}^i(\om\!=\!1)$ for $i>0$ generally
have multiple zeros. For instance, in the case of $T(11,6)$:
$\hat{H}^1(\om\!=\!1)\!=\!\Phi_5^2 \Phi_7 \Phi_{10} \Phi_{12} 
\Phi_{13} \Phi_{14} \Phi_{15}$ for cyclotomic polynomials
$\Phi_m$. It has $4$ pairs of irregular (non-RH) zeros in 
any punctured neighborhood of $\om\!=\!1$. Strong $RH$ begins 
only after $\varpi_1=1.1381148969721394\ldots=\om^{top}_1$. 


\subsubsection{\sf Irregular zeros}
We provide the number of super-dual-invariant pairs  
of {\em irregular zeros\,} in the last column;
the upper two numbers are for $a=0$, the lower ones
are for the maximum among $i>0$. The first gives the
number of pairs right after $\varpi$; the second, after
"$-$", is for the stable number for $\om>\!>0$. Mostly
they coincide in this table.
\vskip 0.2cm

{\em Trivial zeros\,}, $-1,-\om$, are not
counted in the table.
For $\om\!>\!\om^{top}_i$ such zeros appear
for  \Yboxdim5pt$\yng(2,2)\ \,$\Yboxdim7pt
for odd $i$. They are likely to reflect some symmetries
of the DAHA construction, and presumably can be interpreted
geometrically.
We note that more general squares 
result in huge superpolynomials, which we do not have
by now.
\vskip 0.2cm

If the number of pairs of {\em irregular zeros\,} is
$\{\cdot\}-0$ in the last column, then Strong $RH$ holds
for $\om>\!>0$; multiple $RH$\~zeros appear after $\varpi$
in this case. 
If it is $0-0$, then $RH$ holds starting with $\varpi$ from
the table. Thus  Strong $RH$
always holds in the table for uncolored algebraic knots
and those colored by \Yboxdim5pt$\yng(2,2)\,$\Yboxdim7pt
(where we disregard trivial zeros).

Recall that Weak $RH$ allows $\kappa-1$
pairs of irregular zeros
for {\em uncolored\,} algebraic links, 
 where $\kappa$ is the number
of components of a link. 
Practically, the number of allowed
irregular pairs (disregarding trivial ones)
equals the total number of symbols $P$ in the
DAHA presentation of a link minus $1$.
Mostly this difference is $0$ or $1$ in the table. It is $1$ 
for $2$\~links and 
the corresponding irregular pair
is automatically real (the number of irregular {\em complex\,}
pairs is always even). For  entry $33$ (with $N_z^0=24$), there
are four $P$ in its DAHA presentation and, indeed, $\hat{H}^i$
have $3$ irregular
pairs of zeros. Similarly, $RH$ holds for
$\{\hat{\ga}(3,2)(P^3)\}$ (entry $52$ with $N_z^0=90$), 
with $2$ pairs of (automatically complex) irregular roots.
We provide one algebraic {\em link\,} (entry $64$ with $N_z^0=34$)
where Weak $RH$ fails for 
$i=4$ (only for such $i$). See below.


\comment{
{\sf Family $T(3m\pm 1,3)$.} 
For any uncolored algebraic {\em knots\,}, we conjecture that
irregular zeros do not appear, i.e. there are
no {\em stable\,} irregular pairs.
Actually {\em unstable\,} pairs, which 
disappear for large $\om$, are of interest too.
For $T(3m\!\pm\! 1,3)$, exactly
one pair of unstable irregular
(real) zeros appears before $\varpi_i$
for odd powers $a^i$
and for $T(6\pm 1 + 6k,3)$ for $k\ge 0$. The whole series
$T(3m\pm 1,3)$ has $\varpi<1$.
}

\comment{
We mention that there are general conjectured formulas for
uncolored torus knots $T(m,km+1)$; they can be justified in the
combinatorial setting and/or generalized using the dimension
formulas from \cite{Pi} as in \cite{GM,Mel1,ChP1}. See also
\cite{ORS,DMMSS,FGS}; the proof of the formula for DAHA
superpolynomials colored by any rows for $T(2m+ 1,2)$ is in
\cite{CJJ}. Even when explicit formulas are known/conjectured,
they seem generally involved to analyze $RH$
algebraically, but can be helpful to examine the point $\om=1$
for torus knots and for similar questions.
}

\subsection{\bf Composite theory, etc}
We do not discuss in this work
other root systems. The $a$\~stabilization was conjectured
for classical series in \cite{CJ}; the corresponding 
polynomials are called DAHA {\em hyperpolynomials}. 
The {\em hyper-duality\,} is expected to hold too;
one can expect the corresponding $RH$ for algebraic knots/links.
However, the corresponding DAHA
hyperpolynomials are known so far only for small knots.
We think it makes some sense to provide at least one example
from the {\em composite theory\,} for the exceptional
series from \cite{DG}, which
is topologically for the annulus multiplied by $\R^1$
instead of $\S^3$. 
Algebraically,
this is the case of $a$\~stabilization when {\em two\,} Young
diagrams are placed at the opposite ends of the (nonaffine)
Dynkin diagram of type $A$
and the stabilization is with respect to the distance between
these diagrams.

We will consider only the case of uncolored $T(4,3)$ from
\cite{ChE}. Then
$RH$ is ``$OK$"  for the corresponding
{\em composite superpolynomial\,}. Namely, it holds for $i=0,1$
but fails for $i=2$ (deg$_a=5$ in this case). One has:
$\varpi_0=0.84405$,  $\varpi_1=0.6874328$. The first
coincides with $\om_0^{top}$, the second is among the
toots $\om=1/q$ of $D^1$ (which is always true for $\varpi$),
but smaller than   $\om_1^{top}=2$. The
factor of $D^1$ corresponding to $\varpi_1$ is the square of:
$$-9248 + 12492 q - 14345 q^2  + 844 q^3  + 6308 q^4  -
1608 q^5  + 112 q^6.$$

Some recalculation
is necessary from the setting of \cite{ChE}
to make the super-duality exactly as in the present work. So 
we will provide the
corresponding  $\hat{\h}(q,t,a)$:

\renewcommand{\baselinestretch}{0.5}
{\small
\(
1+2 q t+2 q^2 t+3 q^2 t^2+2 q^3 t^2+q^4 t^2+4 q^3 t^3+2 q^4 t^3
+3 q^4 t^4+2 q^5 t^4+2 q^5 t^5+q^6 t^6+a^5 \bigl(-q^5+q^6-q^4 t
+2 q^5 t-q^6 t+q^4 t^2-q^5 t^2\bigr)+a^4 \bigl(-q^3+2 q^4-q^5
+q^2 t-4 q^4 t+4 q^5 t-q^6 t+q^3 t^2+q^4 t^2-4 q^5 t^2+2 q^6 t^2
+q^4 t^3-q^6 t^3+q^5 t^4\bigr)+a^3 \bigl(q-q^2-3 q^3+3 q^4+q^5
-q^6+q t+3 q^2 t-q^3 t-8 q^4 t+4 q^5 t+q^6 t+2 q^2 t^2+5 q^3 t^2
-2 q^4 t^2-8 q^5 t^2+3 q^6 t^2+3 q^3 t^3+5 q^4 t^3-q^5 t^3
-3 q^6 t^3+2 q^4 t^4+3 q^5 t^4-q^6 t^4+q^5 t^5+q^6 t^5\bigr)
+a^2 \bigl(1+2 q-2 q^2-3 q^3+q^5+q^6+4 q t+5 q^2 t-2 q^3 t
-8 q^4 t+q^6 t+7 q^2 t^2+9 q^3 t^2-2 q^4 t^2-8 q^5 t^2+8 q^3 t^3
+9 q^4 t^3-2 q^5 t^3-3 q^6 t^3+7 q^4 t^4+5 q^5 t^4-2 q^6 t^4
+4 q^5 t^5+2 q^6 t^5+q^6 t^6\bigr)+a \bigl(2+q-q^2-q^3-q^4
+5 q t+5 q^2 t-q^3 t-3 q^4 t-2 q^5 t+8 q^2 t^2+7 q^3 t^2+q^4 t^2
-3 q^5 t^2-q^6 t^2+9 q^3 t^3+7 q^4 t^3-q^5 t^3-q^6 t^3
+8 q^4 t^4+5 q^5 t^4-q^6 t^4+5 q^5 t^5+q^6 t^5+2 q^6 t^6\bigr).
\)
}
\renewcommand{\baselinestretch}{1.2}
\vskip 0.2cm

Let us also mention here that {\em Heegaard-Floer homology\,} is
the specialization of Khovanov-Rozansky link homology at for
the differential at $a=-1$
(in DAHA parameters). This specialization
preserves the super-duality, in
contrast to the differentials at $a=-t^{n+1}$ (in DAHA parameters)
to the $A_n$\~theories.
Practically everything we conjecture for $\hat{H}^i$ is 
applicable
when $a=-1$. In this paper, we do not discuss polynomials/series
that are sums of $\h^i$ over all $i$, which is a clear
potential of the theory. The specialization $a\mapsto -(t/q)$
from (\ref{conzeconj}) is of obvious importance too. 
Also, employing the connection 
conjectures,
one has an opportunity to interpret $RH$ for
the topological and other superpolynomials, including the
HOMFLY-PT polynomials, KhR polynomials, physical ones and
those associated with the rational DAHA and the Hilbert scheme of 
$\C^2$.

\section{\sc Some formulas, conclusion}
\subsection{\bf A non-RH link for
\texorpdfstring{{\mathversion{bold}$i\!>\!0$}}{nonzero i}}
\label{SEC:FAIL}
The only counterexample to the total (all $a^i$)
Weak $RH$ among algebraic knots/links in the table
is for the entry $64$ ($N_z^0=34$),
where it fails for $a^4$ (but holds for all other $a^i$).
See also Section \ref{SEC:MULTZ},
(4) and Section \ref{SEC:2TREF} below. Let us discuss it.

The corresponding cable is $\l_{1}=
(C\!ab(5,3),C\!ab(5,3))T(1,1)$,
which is the link of the singularity 
$\c_1=\{(x^5-y^3)(x^3-y^5)=0\}$
at $x\!=\!0\!=\!y$.
The linking number is $lk=9$ for its two components, which are
$T(5,3)$. The zeta-monodromy from \cite{DGPS}
upon $t\mapsto q$,
essentially the Alexander polynomial, is
$Z_1=q^{32}+2q^{24}+3q^{16}+2q^8+1$. According to
Section 5.4 from \cite{ChD2}, the following connection with our
superpolynomials is expected (unless for the unknot):
\begin{align}\label{zet-alex}
Z=\hat{\h}^{min}_{\l}\,
(q,q,a=-1)/(1-q)^{\kappa-\de_{\kappa,1}}.
\end{align}
This is for {\em uncolored\,} graph $\l$ (without $'\!\l$)
with $\kappa$ paths (the number of connected components in the
corresponding cable). The linking number is then $Z(q\!=\!1)$.
Recall that we always impose the minimal normalization
$\hat{\h}^{min}$
in the present paper. 

Let us mention that the link
$2\,T(5,3)=(C\!ab(1,1),C\!ab(1,1))T(5,3)$, corresponding
to $\c'_1=\{(x^5\!-\!y^3)(x^5\!+\!y^3)=0\}$ with the 
linking number $Z_1'(1)=15$,
satisfies $RH$. One has: $Z'_1=
q^{44}+q^{38}+q^{34}+q^{32}+q^{28}+q^{26}+q^{24}+q^{22}
+q^{20}+q^{18}+q^{16}+q^{12}+q^{10}+q^6+1$ in this case.
This is actually entry No$=40$ (with $N_z^0=46$), because
the cables $C\!ab(1,1)T(3,2)$ and $T(5,3)$ are isotopic.  
Note that the degree of $Z$ is $N_z^0-2$.
 
Also, $RH$ holds for all $i$ for the following direct
modifications of $\l_{1}$:
$\l_{2}=(C\!ab(8,3),C\!ab(8,3))T(2,1)$,
$\l_3=(C\!ab(11,3),C\!ab(11,3))T(3,1)$,
which are entries $44,45$ $(N_z^0=64,94)$.
They are correspondingly unions of two copies of 
$T(8,3)$ and $T(11,3)$ with linking numbers $18$ and $27$.

Finding the equations of the corresponding plane
curve singularities is more involved in these examples.
Say for $\c_2$, we begin with $(x^8-y^3)(x^3-y^5)=0$ and replace
$x^3=y^5$ by ``its double", which is $x^3=y^8$, provided
that the corresponding link is in the vicinity of
$T(5,3)$ (then the resulting linking number becomes $18$).
This is the meaning of cabling in this case. 
The equations and  $Z$\~polynomials are as follows:
{\small
\begin{align*}
\c_2=&\,\{(x^8\!-\!y^3)((y\!+\!x^2)^3\!+\!x^8)\!=\!0\},\ \ 
\c_3=\{(x^{11}\!-\!y^3)((y\!+\!x^3)^3\!+\!x^{11})\!=\!0\},\\
Z_2=\,&1\!+\!q^6\!+\!2 q^{14}\!+\!2 q^{20}\!+\!3 q^{28}\!
+\!3 q^{34}
\!+\!2 q^{42}\!+\!2 q^{48}\!+\!q^{56}\!+\!q^{62},\ \ \,Z_3=
1\!+\!q^6 +\notag\\
\!+q^{12}\!+&2 q^{20}\!+\!2 q^{26}\!+\!2 q^{32}\!+\!3 q^{40}\!
+\!3 q^{46}\!+\!3 q^{52}
\!+\!2 q^{60}\!+\!2 q^{66}\!+\!2 q^{72}\!+\!q^{80}\!+\!q^{86}\!
+\!q^{92}.\notag
\end{align*}
}

\subsubsection{\sf Non-RH superpolynomial} 
The failure of Weak $RH$ for the link
of $\c_1=\{(x^5-y^3)(x^3-y^5)=0\}$
above is
at $a^4$; let us provide $\hat{\h}^{min}_{\l}\,:$
\vskip 0.1cm

\renewcommand{\baselinestretch}{0.5}
{\small
\(
1-t+q t+q^2 t+q^3 t+q^4 t+q^5 t-q t^2+q^4 t^2+2 q^5 t^2+4 q^6 t^2
+q^7 t^2+q^8 t^2-q^2 t^3-q^4 t^3-q^5 t^3+5 q^7 t^3+4 q^8 t^3
+4 q^9 t^3-q^3 t^4-q^5 t^4-2 q^6 t^4-2 q^7 t^4+3 q^8 t^4
+5 q^9 t^4+7 q^{10} t^4+q^{11} t^4-q^4 t^5-q^6 t^5-2 q^7 t^5
-4 q^8 t^5+q^9 t^5+5 q^{10} t^5+8 q^{11} t^5+2 q^{12} t^5-q^5 t^6
-q^7 t^6-2 q^8 t^6-5 q^9 t^6+5 q^{11} t^6+8 q^{12} t^6+q^{13} t^6
-q^6 t^7-q^8 t^7-2 q^9 t^7-6 q^{10} t^7+5 q^{12} t^7+7 q^{13} t^7
-q^7 t^8-q^9 t^8-2 q^{10} t^8-5 q^{11} t^8+q^{12} t^8
+5 q^{13} t^8+4 q^{14} t^8-q^8 t^9-q^{10} t^9-2 q^{11} t^9
-4 q^{12} t^9+3 q^{13} t^9+4 q^{14} t^9+q^{15} t^9-q^9 t^{10}
-q^{11} t^{10}-2 q^{12} t^{10}-2 q^{13} t^{10}+5 q^{14} t^{10}
+q^{15} t^{10}-q^{10} t^{11}-q^{12} t^{11}-2 q^{13} t^{11}
+4 q^{15} t^{11}-q^{11} t^{12}-q^{13} t^{12}-q^{14} t^{12}
+2 q^{15} t^{12}+q^{16} t^{12}-q^{12} t^{13}-q^{14} t^{13}
+q^{15} t^{13}+q^{16} t^{13}-q^{13} t^{14}+q^{16} t^{14}
-q^{14} t^{15}+q^{16} t^{15}-q^{15} t^{16}+q^{16} t^{16}
-q^{16} t^{17}+q^{17} t^{17}+a^5 \bigl(q^{15}+2 q^{16} t
+q^{17} t^2\bigr)
\)
\vfil
\noindent
\(
+\mathbf{ a^4 \bigl(q^{10}+q^{11}}
+q^{12}+q^{13}+q^{14}
-q^{10} t+q^{11} t+3 q^{12} t+3 q^{13} t+3 q^{14} t+3 q^{15} t
-2 q^{11} t^2-q^{12} t^2+3 q^{13} t^2+4 q^{14} t^2+4 q^{15} t^2
+3 q^{16} t^2-3 q^{12} t^3-q^{13} t^3+2 q^{14} t^3+4 q^{15} t^3
+3 q^{16} t^3+q^{17} t^3-4 q^{13} t^4-q^{14} t^4+3 q^{15} t^4
+3 q^{16} t^4+q^{17} t^4-3 q^{14} t^5-q^{15} t^5+3 q^{16} t^5
+q^{17} t^5-2 q^{15} t^6+q^{16} t^6+q^{17} t^6-q^{16} t^7
+ \mathbf{q^{17}t^7\bigr)}
\)
\vfil
\noindent
\(
+a^3 \bigl(q^6+q^7+2 q^8+2 q^9
+2 q^{10}+q^{11}
+q^{12}-q^6 t+q^8 t+4 q^9 t+6 q^{10} t+8 q^{11} t+5 q^{12} t
+4 q^{13} t+q^{14} t-q^7 t^2-2 q^8 t^2-3 q^9 t^2+q^{10} t^2
+7 q^{11} t^2+14 q^{12} t^2+10 q^{13} t^2+7 q^{14} t^2
+2 q^{15} t^2-q^8 t^3-2 q^9 t^3-6 q^{10} t^3-4 q^{11} t^3
+5 q^{12} t^3+16 q^{13} t^3+12 q^{14} t^3+7 q^{15} t^3
+q^{16} t^3-q^9 t^4-2 q^{10} t^4-8 q^{11} t^4-8 q^{12} t^4
+5 q^{13} t^4+16 q^{14} t^4+10 q^{15} t^4+4 q^{16} t^4-q^{10} t^5
-2 q^{11} t^5-8 q^{12} t^5-8 q^{13} t^5+5 q^{14} t^5
+14 q^{15} t^5+5 q^{16} t^5+q^{17} t^5-q^{11} t^6-2 q^{12} t^6
-8 q^{13} t^6-4 q^{14} t^6+7 q^{15} t^6+8 q^{16} t^6+q^{17} t^6
-q^{12} t^7-2 q^{13} t^7-6 q^{14} t^7+q^{15} t^7+6 q^{16} t^7
+2 q^{17} t^7-q^{13} t^8-2 q^{14} t^8-3 q^{15} t^8+4 q^{16} t^8
+2 q^{17} t^8-q^{14} t^9-2 q^{15} t^9+q^{16} t^9+2 q^{17} t^9
-q^{15} t^{10}+q^{17} t^{10}-q^{16} t^{11}+q^{17} t^{11}\bigr)
+a^2 \bigl(q^3+q^4+2 q^5+2 q^6+2 q^7+q^8+q^9-q^3 t+3 q^6 t
+6 q^7 t+8 q^8 t+7 q^9 t+6 q^{10} t+2 q^{11} t+q^{12} t-q^4 t^2
-q^5 t^2-3 q^6 t^2-2 q^7 t^2+3 q^8 t^2+12 q^9 t^2+16 q^{10} t^2
+15 q^{11} t^2+7 q^{12} t^2+3 q^{13} t^2-q^5 t^3-q^6 t^3
-4 q^7 t^3-6 q^8 t^3-6 q^9 t^3+8 q^{10} t^3+20 q^{11} t^3
+23 q^{12} t^3+11 q^{13} t^3+4 q^{14} t^3-q^6 t^4-q^7 t^4
-4 q^8 t^4-8 q^9 t^4-13 q^{10} t^4+2 q^{11} t^4+21 q^{12} t^4
+26 q^{13} t^4+11 q^{14} t^4+3 q^{15} t^4-q^7 t^5-q^8 t^5
-4 q^9 t^5-9 q^{10} t^5-16 q^{11} t^5-q^{12} t^5+21 q^{13} t^5
+23 q^{14} t^5+7 q^{15} t^5+q^{16} t^5-q^8 t^6-q^9 t^6
-4 q^{10} t^6-10 q^{11} t^6-16 q^{12} t^6+2 q^{13} t^6
+20 q^{14} t^6+15 q^{15} t^6+2 q^{16} t^6-q^9 t^7-q^{10} t^7
-4 q^{11} t^7-9 q^{12} t^7-13 q^{13} t^7+8 q^{14} t^7
+16 q^{15} t^7+6 q^{16} t^7-q^{10} t^8-q^{11} t^8-4 q^{12} t^8
-8 q^{13} t^8-6 q^{14} t^8+12 q^{15} t^8+7 q^{16} t^8+q^{17} t^8
-q^{11} t^9-q^{12} t^9-4 q^{13} t^9-6 q^{14} t^9+3 q^{15} t^9
+8 q^{16} t^9+q^{17} t^9-q^{12} t^{10}-q^{13} t^{10}
-4 q^{14} t^{10}-2 q^{15} t^{10}+6 q^{16} t^{10}+2 q^{17} t^{10}
-q^{13} t^{11}-q^{14} t^{11}-3 q^{15} t^{11}+3 q^{16} t^{11}
+2 q^{17} t^{11}-q^{14} t^{12}-q^{15} t^{12}+2 q^{17} t^{12}
-q^{15} t^{13}+q^{17} t^{13}-q^{16} t^{14}+q^{17} t^{14}\bigr)
+a \bigl(q+q^2+q^3+q^4+q^5-q t+q^3 t+2 q^4 t+4 q^5 t+6 q^6 t
+3 q^7 t+2 q^8 t+q^9 t-q^2 t^2-q^3 t^2-q^4 t^2-q^5 t^2+2 q^6 t^2
+10 q^7 t^2+10 q^8 t^2+9 q^9 t^2+4 q^{10} t^2+q^{11} t^2-q^3 t^3
-q^4 t^3-2 q^5 t^3-4 q^6 t^3-4 q^7 t^3+5 q^8 t^3+13 q^9 t^3
+18 q^{10} t^3+9 q^{11} t^3+3 q^{12} t^3-q^4 t^4-q^5 t^4
-2 q^6 t^4-5 q^7 t^4-8 q^8 t^4-2 q^9 t^4+11 q^{10} t^4
+22 q^{11} t^4+14 q^{12} t^4+4 q^{13} t^4-q^5 t^5-q^6 t^5
-2 q^7 t^5-5 q^8 t^5-11 q^9 t^5-7 q^{10} t^5+9 q^{11} t^5
+23 q^{12} t^5+14 q^{13} t^5+3 q^{14} t^5-q^6 t^6-q^7 t^6
-2 q^8 t^6-5 q^9 t^6-13 q^{10} t^6-8 q^{11} t^6+9 q^{12} t^6
+22 q^{13} t^6+9 q^{14} t^6+q^{15} t^6-q^7 t^7-q^8 t^7-2 q^9 t^7
-5 q^{10} t^7-13 q^{11} t^7-7 q^{12} t^7+11 q^{13} t^7
+18 q^{14} t^7+4 q^{15} t^7-q^8 t^8-q^9 t^8-2 q^{10} t^8
-5 q^{11} t^8-11 q^{12} t^8-2 q^{13} t^8+13 q^{14} t^8
+9 q^{15} t^8+q^{16} t^8-q^9 t^9-q^{10} t^9-2 q^{11} t^9
-5 q^{12} t^9-8 q^{13} t^9+5 q^{14} t^9+10 q^{15} t^9
+2 q^{16} t^9-q^{10} t^{10}-q^{11} t^{10}-2 q^{12} t^{10}
-5 q^{13} t^{10}-4 q^{14} t^{10}+10 q^{15} t^{10}+3 q^{16} t^{10}
-q^{11} t^{11}-q^{12} t^{11}-2 q^{13} t^{11}-4 q^{14} t^{11}
+2 q^{15} t^{11}+6 q^{16} t^{11}-q^{12} t^{12}-q^{13} t^{12}
-2 q^{14} t^{12}-q^{15} t^{12}+4 q^{16} t^{12}+q^{17} t^{12}
-q^{13} t^{13}-q^{14} t^{13}-q^{15} t^{13}+2 q^{16} t^{13}
+q^{17} t^{13}-q^{14} t^{14}-q^{15} t^{14}+q^{16} t^{14}
+q^{17} t^{14}
-q^{15} t^{15}+q^{17} t^{15}-q^{16} t^{16}+q^{17} t^{16}\bigr).
\)
}
\renewcommand{\baselinestretch}{1.2}

\vskip 0.2cm

It makes sense to provide $\hat{H}^4$, responsible for the
failure of $RH$. It is
$\h^4$ upon the substitution 
$q\mapsto qt$ and under the hat-normalization: 

\renewcommand{\baselinestretch}{0.5}
{\small

\(\hat{H}^4(1+qt^2)^{-2}=
1 - t + q t - q t^2 + q^2 t^2 + q^2 t^3 + q^3 t^3 + q^3 t^4 +
 q^4 t^4 - 2 q^2 t^5 + q^4 t^5 + q^4 t^6 + q^5 t^6 + q^4 t^7 +
 q^5 t^7 - q^4 t^8 + q^5 t^8 - q^4 t^9 + q^5 t^9 + q^5 t^{10}.
\)
}
\renewcommand{\baselinestretch}{1.2}
We note that $\varsigma_4 t^{\pi_4}S^4\!=\!-t\bigl((1\!+\!t^2)
(1\!+\!t^4)\bigr)^2.$ 
\comment{
which has essentially random norms of the zeros (though
all in the vicinity of $\sqrt{\om}$).
An extension of \cite{ChP1}
to links is in progress; a geometric
interpretation of such $H$ is likely to exist, but
this failure  
can be due to complicated nature of the corresponding
flagged variety.} 
In the case of
{\em Heegaard-Floer substitution\,}, which  
is $a=-1$ in $\hat{\h}_\l$, Weak $RH$ holds for this link.

\subsection{\bf An example of 3-link}
The appearance of $\kappa-1$ (super-dual)
pairs of zeros certainly deserves a comment.
Non-real non-$RH$ roots can occur only when the number of
branches is $\kappa>2$ and they do appear. Let 
$\l=3\,T(4,3)=(C\!ab(1,1),C\!ab(1,1),C\!ab(1,1))T(4,3)$ be
the link of 
$(x^4-y^3)(x^4+y^3)(x^4+2y^3)=0$, which has $3$ components
$T(4,3)$ and the pairwise linking numbers $12$. The corresponding
zeta-monodromy from (\ref{zet-alex}) is
{\small
$Z\!=\!1+q^9+q^{12}+q^{18}+q^{21}+q^{24}+q^{27}+q^{30}+q^{33}
-q^{36}+q^{39}+q^{42}-q^{45}-q^{48}+q^{51}-q^{54}-
q^{57}-q^{60}-q^{63}-q^{66}-q^{69}-q^{75}-q^{78}-q^{87}.$
}
\vskip 0.1cm

This example is entry No$=52$ with $N_z^0=90$
from the table. The corresponding
$\hat{\h}$ is large ($4308$ different $a,q,t$\~monomials).
The polynomial $\hat{H}^0$ is of degree $90$ with respect
to $t$ (and $45$ with respect to $q$). Weak $RH$ holds
after $\varpi_0=1.55727521033844259502499...$, coinciding
with the top real zero $\om_0^{top}$
of the reduced discriminant $D^0$.  Let us make $\om=2$.
Then there is only one non-$RH$ pair of zeros up to
the complex conjugation:
\begin{align*}
&\xi\,=\,1.99999963844688175553480 + 0.00001272650573633190499 
\,\imath,\\
&\frac{\om}{\xi}=
1.00000018073610079338078 - 6.363255168562838442\times
10^{-6}\imath.
\end{align*}
Their product is $2$.
The first quickly approaches the corresponding $\om$ as
$\om\to\infty$
(the difference is approximately some power 
of $1/\om$).
These zeros cannot become real for any $\om>\om^{top}_0$
because this would result in multiple roots after $\om^{top}_0$.
Weak $RH$ holds for $\om=2$  and any
$\hat{H}^i,\, 0\le i\le 8=$deg$_{a}$
with $2$ pairs of complex irregular zeros.
The number of
zeros $\xi$ is correspondingly $90,88,84,78,70,60,48,34,18$
for $i=0,\ldots,8$.

\subsubsection{\sf Superpolynomial at
\texorpdfstring{{\mathversion{bold}$a=0$}}{zero}}
Let us provide $\h^0=\hat{\h}(a=0)$ (not $\hat{H}^0$,
i.e. without the substitution $q\mapsto q t$)
for this $3$\~link:

\renewcommand{\baselinestretch}{0.5}
{\footnotesize
\(
1-2 t+q t+q^2 t+q^3 t+q^4 t+q^5 t+q^6 t+q^7 t+q^8 t+t^2-2 q t^2
-q^2 t^2-q^3 t^2+q^6 t^2+q^7 t^2+2 q^8 t^2+4 q^9 t^2+4 q^{10} t^2
+2 q^{11} t^2+2 q^{12} t^2+q^{13} t^2+q^{14} t^2+q t^3-q^2 t^3
-2 q^4 t^3-q^5 t^3-2 q^6 t^3-q^7 t^3-2 q^8 t^3-q^9 t^3
+5 q^{11} t^3+5 q^{12} t^3+6 q^{13} t^3+5 q^{14} t^3+6 q^{15} t^3
+3 q^{16} t^3+q^{17} t^3+q^{18} t^3+q^2 t^4-q^3 t^4+q^4 t^4
-q^5 t^4-q^6 t^4-2 q^7 t^4-q^8 t^4-4 q^9 t^4-3 q^{10} t^4
-5 q^{11} t^4-2 q^{12} t^4+4 q^{14} t^4+4 q^{15} t^4
+11 q^{16} t^4+11 q^{17} t^4+7 q^{18} t^4+5 q^{19} t^4
+3 q^{20} t^4+q^3 t^5-q^4 t^5+q^5 t^5-2 q^8 t^5-3 q^{10} t^5
-3 q^{11} t^5-6 q^{12} t^5-5 q^{13} t^5-7 q^{14} t^5-2 q^{15} t^5
-5 q^{16} t^5+3 q^{17} t^5+10 q^{18} t^5+14 q^{19} t^5
+12 q^{20} t^5+11 q^{21} t^5+5 q^{22} t^5+q^{23} t^5+q^4 t^6
-q^5 t^6+q^6 t^6+q^8 t^6-q^9 t^6-2 q^{11} t^6-q^{12} t^6
-5 q^{13} t^6-4 q^{14} t^6-9 q^{15} t^6-6 q^{16} t^6
-10 q^{17} t^6-6 q^{18} t^6-4 q^{19} t^6+7 q^{20} t^6
+14 q^{21} t^6+18 q^{22} t^6+14 q^{23} t^6+8 q^{24} t^6
+q^{25} t^6+q^5 t^7-q^6 t^7+q^7 t^7+q^9 t^7+q^{11} t^7
-2 q^{12} t^7-3 q^{14} t^7-2 q^{15} t^7-7 q^{16} t^7-6 q^{17} t^7
-13 q^{18} t^7-8 q^{19} t^7-12 q^{20} t^7-8 q^{21} t^7
+2 q^{22} t^7+17 q^{23} t^7+21 q^{24} t^7+18 q^{25} t^7
+8 q^{26} t^7+q^{27} t^7+q^6 t^8-q^7 t^8+q^8 t^8+q^{10} t^8
+2 q^{12} t^8-q^{13} t^8-2 q^{15} t^8-5 q^{17} t^8-3 q^{18} t^8
-12 q^{19} t^8-9 q^{20} t^8-13 q^{21} t^8-14 q^{22} t^8
-14 q^{23} t^8+3 q^{24} t^8+21 q^{25} t^8+26 q^{26} t^8
+15 q^{27} t^8+6 q^{28} t^8+q^{29} t^8+q^7 t^9-q^8 t^9+q^9 t^9
+q^{11} t^9+2 q^{13} t^9+q^{15} t^9-2 q^{16} t^9+q^{17} t^9
-3 q^{18} t^9-q^{19} t^9-9 q^{20} t^9-7 q^{21} t^9-13 q^{22} t^9
-13 q^{23} t^9-21 q^{24} t^9-15 q^{25} t^9+10 q^{26} t^9
+31 q^{27} t^9+23 q^{28} t^9+12 q^{29} t^9+3 q^{30} t^9
+q^8 t^{10}-q^9 t^{10}+q^{10} t^{10}+q^{12} t^{10}
+2 q^{14} t^{10}+2 q^{16} t^{10}-q^{17} t^{10}+q^{18} t^{10}
-2 q^{19} t^{10}+q^{20} t^{10}-7 q^{21} t^{10}-4 q^{22} t^{10}
-11 q^{23} t^{10}-12 q^{24} t^{10}-20 q^{25} t^{10}
-25 q^{26} t^{10}-5 q^{27} t^{10}+25 q^{28} t^{10}
+30 q^{29} t^{10}+18 q^{30} t^{10}+5 q^{31} t^{10}+q^{32} t^{10}
+q^9 t^{11}-q^{10} t^{11}+q^{11} t^{11}+q^{13} t^{11}
+2 q^{15} t^{11}+2 q^{17} t^{11}+2 q^{19} t^{11}-2 q^{20} t^{11}
+2 q^{21} t^{11}-5 q^{22} t^{11}-2 q^{23} t^{11}-8 q^{24} t^{11}
-10 q^{25} t^{11}-20 q^{26} t^{11}-27 q^{27} t^{11}
-14 q^{28} t^{11}+13 q^{29} t^{11}+34 q^{30} t^{11}
+23 q^{31} t^{11}+7 q^{32} t^{11}+q^{33} t^{11}+q^{10} t^{12}
-q^{11} t^{12}+q^{12} t^{12}+q^{14} t^{12}+2 q^{16} t^{12}
+2 q^{18} t^{12}+3 q^{20} t^{12}-q^{21} t^{12}+2 q^{22} t^{12}
-4 q^{23} t^{12}-6 q^{25} t^{12}-7 q^{26} t^{12}-20 q^{27} t^{12}
-29 q^{28} t^{12}-17 q^{29} t^{12}+7 q^{30} t^{12}
+33 q^{31} t^{12}+23 q^{32} t^{12}+9 q^{33} t^{12}+q^{34} t^{12}
+q^{11} t^{13}-q^{12} t^{13}+q^{13} t^{13}+q^{15} t^{13}
+2 q^{17} t^{13}+2 q^{19} t^{13}+3 q^{21} t^{13}+3 q^{23} t^{13}
-4 q^{24} t^{13}+q^{25} t^{13}-4 q^{26} t^{13}-5 q^{27} t^{13}
-18 q^{28} t^{13}-30 q^{29} t^{13}-24 q^{30} t^{13}
+6 q^{31} t^{13}+35 q^{32} t^{13}+23 q^{33} t^{13}
+7 q^{34} t^{13}+q^{35} t^{13}+q^{12} t^{14}-q^{13} t^{14}
+q^{14} t^{14}+q^{16} t^{14}+2 q^{18} t^{14}+2 q^{20} t^{14}
+3 q^{22} t^{14}+4 q^{24} t^{14}-3 q^{25} t^{14}+q^{26} t^{14}
-3 q^{27} t^{14}-4 q^{28} t^{14}-16 q^{29} t^{14}
-29 q^{30} t^{14}-26 q^{31} t^{14}+6 q^{32} t^{14}
+33 q^{33} t^{14}+23 q^{34} t^{14}+5 q^{35} t^{14}+q^{13} t^{15}
-q^{14} t^{15}+q^{15} t^{15}+q^{17} t^{15}+2 q^{19} t^{15}
+2 q^{21} t^{15}+3 q^{23} t^{15}+4 q^{25} t^{15}-2 q^{26} t^{15}
+2 q^{27} t^{15}-3 q^{28} t^{15}-4 q^{29} t^{15}-15 q^{30} t^{15}
-29 q^{31} t^{15}-24 q^{32} t^{15}
+7 q^{33} t^{15}
+34 q^{34} t^{15}+18 q^{35} t^{15}+3 q^{36} t^{15}+q^{14} t^{16}
-q^{15} t^{16}+q^{16} t^{16}+q^{18} t^{16}+2 q^{20} t^{16}
+2 q^{22} t^{16}+3 q^{24} t^{16}+4 q^{26} t^{16}-2 q^{27} t^{16}
+3 q^{28} t^{16}-3 q^{29} t^{16}-4 q^{30} t^{16}-16 q^{31} t^{16}
-30 q^{32} t^{16}-17 q^{33} t^{16}+13 q^{34} t^{16}
+30 q^{35} t^{16}+12 q^{36} t^{16}+q^{37} t^{16}+q^{15} t^{17}
-q^{16} t^{17}+q^{17} t^{17}+q^{19} t^{17}+2 q^{21} t^{17}
+2 q^{23} t^{17}+3 q^{25} t^{17}+4 q^{27} t^{17}-2 q^{28} t^{17}
+3 q^{29} t^{17}-3 q^{30} t^{17}-4 q^{31} t^{17}-18 q^{32} t^{17}
-29 q^{33} t^{17}-14 q^{34} t^{17}+25 q^{35} t^{17}
+23 q^{36} t^{17}+6 q^{37} t^{17}+q^{16} t^{18}-q^{17} t^{18}
+q^{18} t^{18}+q^{20} t^{18}+2 q^{22} t^{18}+2 q^{24} t^{18}
+3 q^{26} t^{18}+4 q^{28} t^{18}-2 q^{29} t^{18}+2 q^{30} t^{18}
-3 q^{31} t^{18}-5 q^{32} t^{18}-20 q^{33} t^{18}-27 q^{34} t^{18}
-5 q^{35} t^{18}+31 q^{36} t^{18}+15 q^{37} t^{18}+q^{38} t^{18}
+q^{17} t^{19}-q^{18} t^{19}+q^{19} t^{19}+q^{21} t^{19}
+2 q^{23} t^{19}+2 q^{25} t^{19}+3 q^{27} t^{19}+4 q^{29} t^{19}
-2 q^{30} t^{19}+q^{31} t^{19}-4 q^{32} t^{19}-7 q^{33} t^{19}
-20 q^{34} t^{19}-25 q^{35} t^{19}+10 q^{36} t^{19}
+26 q^{37} t^{19}+8 q^{38} t^{19}+q^{18} t^{20}-q^{19} t^{20}
+q^{20} t^{20}+q^{22} t^{20}+2 q^{24} t^{20}+2 q^{26} t^{20}
+3 q^{28} t^{20}+4 q^{30} t^{20}-3 q^{31} t^{20}+q^{32} t^{20}
-6 q^{33} t^{20}-10 q^{34} t^{20}-20 q^{35} t^{20}
-15 q^{36} t^{20}+21 q^{37} t^{20}+18 q^{38} t^{20}+q^{39} t^{20}
+q^{19} t^{21}-q^{20} t^{21}+q^{21} t^{21}+q^{23} t^{21}
+2 q^{25} t^{21}+2 q^{27} t^{21}+3 q^{29} t^{21}+4 q^{31} t^{21}
-4 q^{32} t^{21}-8 q^{34} t^{21}-12 q^{35} t^{21}
-21 q^{36} t^{21}+3 q^{37} t^{21}+21 q^{38} t^{21}
+8 q^{39} t^{21}+q^{20} t^{22}-q^{21} t^{22}+q^{22} t^{22}
+q^{24} t^{22}+2 q^{26} t^{22}+2 q^{28} t^{22}+3 q^{30} t^{22}
+3 q^{32} t^{22}-4 q^{33} t^{22}-2 q^{34} t^{22}-11 q^{35} t^{22}
-13 q^{36} t^{22}-14 q^{37} t^{22}+17 q^{38} t^{22}
+14 q^{39} t^{22}+q^{40} t^{22}+q^{21} t^{23}-q^{22} t^{23}
+q^{23} t^{23}+q^{25} t^{23}+2 q^{27} t^{23}+2 q^{29} t^{23}
+3 q^{31} t^{23}+2 q^{33} t^{23}-5 q^{34} t^{23}-4 q^{35} t^{23}
-13 q^{36} t^{23}-14 q^{37} t^{23}+2 q^{38} t^{23}
+18 q^{39} t^{23}+5 q^{40} t^{23}+q^{22} t^{24}-q^{23} t^{24}
+q^{24} t^{24}+q^{26} t^{24}+2 q^{28} t^{24}+2 q^{30} t^{24}
+3 q^{32} t^{24}-q^{33} t^{24}+2 q^{34} t^{24}-7 q^{35} t^{24}
-7 q^{36} t^{24}-13 q^{37} t^{24}-8 q^{38} t^{24}
+14 q^{39} t^{24}+11 q^{40} t^{24}+q^{23} t^{25}-q^{24} t^{25}
+q^{25} t^{25}+q^{27} t^{25}+2 q^{29} t^{25}+2 q^{31} t^{25}
+3 q^{33} t^{25}-2 q^{34} t^{25}+q^{35} t^{25}-9 q^{36} t^{25}
-9 q^{37} t^{25}-12 q^{38} t^{25}+7 q^{39} t^{25}
+12 q^{40} t^{25}+3 q^{41} t^{25}
+q^{24} t^{26}-q^{25} t^{26}
+q^{26} t^{26}+q^{28} t^{26}+2 q^{30} t^{26}+2 q^{32} t^{26}
+2 q^{34} t^{26}-2 q^{35} t^{26}-q^{36} t^{26}-12 q^{37} t^{26}
-8 q^{38} t^{26}-4 q^{39} t^{26}+14 q^{40} t^{26}+5 q^{41} t^{26}
+q^{25} t^{27}-q^{26} t^{27}+q^{27} t^{27}+q^{29} t^{27}
+2 q^{31} t^{27}+2 q^{33} t^{27}+q^{35} t^{27}-3 q^{36} t^{27}
-3 q^{37} t^{27}-13 q^{38} t^{27}-6 q^{39} t^{27}
+10 q^{40} t^{27}+7 q^{41} t^{27}+q^{42} t^{27}+q^{26} t^{28}
-q^{27} t^{28}+q^{28} t^{28}+q^{30} t^{28}+2 q^{32} t^{28}
+2 q^{34} t^{28}-q^{35} t^{28}+q^{36} t^{28}-5 q^{37} t^{28}
-6 q^{38} t^{28}-10 q^{39} t^{28}+3 q^{40} t^{28}
+11 q^{41} t^{28}+q^{42} t^{28}+q^{27} t^{29}-q^{28} t^{29}
+q^{29} t^{29}+q^{31} t^{29}+2 q^{33} t^{29}+2 q^{35} t^{29}
-2 q^{36} t^{29}-7 q^{38} t^{29}-6 q^{39} t^{29}-5 q^{40} t^{29}
+11 q^{41} t^{29}+3 q^{42} t^{29}+q^{28} t^{30}-q^{29} t^{30}
+q^{30} t^{30}+q^{32} t^{30}+2 q^{34} t^{30}+q^{36} t^{30}
-2 q^{37} t^{30}-2 q^{38} t^{30}-9 q^{39} t^{30}-2 q^{40} t^{30}
+4 q^{41} t^{30}+6 q^{42} t^{30}+q^{29} t^{31}-q^{30} t^{31}
+q^{31} t^{31}+q^{33} t^{31}+2 q^{35} t^{31}-3 q^{38} t^{31}
-4 q^{39} t^{31}-7 q^{40} t^{31}+4 q^{41} t^{31}+5 q^{42} t^{31}
+q^{43} t^{31}+q^{30} t^{32}-q^{31} t^{32}+q^{32} t^{32}
+q^{34} t^{32}+2 q^{36} t^{32}-q^{37} t^{32}-5 q^{39} t^{32}
-5 q^{40} t^{32}+6 q^{42} t^{32}+q^{43} t^{32}+q^{31} t^{33}
-q^{32} t^{33}+q^{33} t^{33}+q^{35} t^{33}+2 q^{37} t^{33}
-2 q^{38} t^{33}-q^{39} t^{33}-6 q^{40} t^{33}-2 q^{41} t^{33}
+5 q^{42} t^{33}+2 q^{43} t^{33}+q^{32} t^{34}-q^{33} t^{34}
+q^{34} t^{34}+q^{36} t^{34}+q^{38} t^{34}-2 q^{39} t^{34}
-3 q^{40} t^{34}-5 q^{41} t^{34}+5 q^{42} t^{34}+2 q^{43} t^{34}
+q^{33} t^{35}-q^{34} t^{35}+q^{35} t^{35}+q^{37} t^{35}
-3 q^{40} t^{35}-3 q^{41} t^{35}+4 q^{43} t^{35}+q^{34} t^{36}
-q^{35} t^{36}+q^{36} t^{36}+q^{38} t^{36}-q^{39} t^{36}
-4 q^{41} t^{36}-q^{42} t^{36}+4 q^{43} t^{36}+q^{35} t^{37}
-q^{36} t^{37}+q^{37} t^{37}+q^{39} t^{37}-2 q^{40} t^{37}
-q^{41} t^{37}-2 q^{42} t^{37}+2 q^{43} t^{37}+q^{44} t^{37}
+q^{36} t^{38}-q^{37} t^{38}+q^{38} t^{38}-2 q^{41} t^{38}
-q^{42} t^{38}+q^{43} t^{38}+q^{44} t^{38}+q^{37} t^{39}
-q^{38} t^{39}+q^{39} t^{39}-q^{41} t^{39}-2 q^{42} t^{39}
+q^{43} t^{39}+q^{44} t^{39}+q^{38} t^{40}-q^{39} t^{40}
+q^{40} t^{40}-q^{41} t^{40}-q^{42} t^{40}+q^{44} t^{40}
+q^{39} t^{41}-q^{40} t^{41}+q^{41} t^{41}-2 q^{42} t^{41}
+q^{44} t^{41}+q^{40} t^{42}-q^{41} t^{42}-q^{43} t^{42}
+q^{44} t^{42}+q^{41} t^{43}-q^{42} t^{43}-q^{43} t^{43}
+q^{44} t^{43}+q^{42} t^{44}
-2 q^{43} t^{44}+q^{44} t^{44}+q^{43} t^{45}-2 q^{44} t^{45}
+q^{45} t^{45}.
\)
}
\renewcommand{\baselinestretch}{1.2}
\vskip 0.2cm

The irregular zeros become more distant
from $\om$ for $H^i$ with $i$ close to deg$_a$, but the
tendency remains the same. The counterpart
of irregular $\xi$ above for $i=$deg$_a=8$ is
$\xi= 1.973849767 + 0.055623630$. 

\comment{
The tendency for  irregular zeros to be close
to $1,\om$, which is getting sharper for $\om>>0$, remains
the same in other examples we considered.}
\vskip 0.2cm

Let us very briefly discuss entry No=$48$ ($N_z^0=80$)
with the link $4\,T(3,2)$, corresponding to the $4$\~branch
plane curve singularity
$(x^3-y^2)(x^3+y^2)(x^3-2y^2)(x^3+2y^2)=0$.
It has $1$ real and $2$ complex pairs of irregular zeros.
Up to the complex conjugation and $\xi\mapsto \om/\xi$,
they are for $\om=2$:
{\small
$1.999451149,\, 2.000252243 + 0.000499389\imath$},
and for $\om=20$:
{\small $19.999999999999995,\, 
20.000000000000003+\!4\!\times\!10^{-15}\imath$}.

\subsection{\bf Some simple cases}
Let us provide the simplest algebraic
uncolored knots, links, and discuss the simplest
non-algebraic cable where Weak $RH$ fails, which are
$C\!ab(-1\!-\!2m)T(3,2)$ for $m\ge 0$. 

\subsubsection{\sf Trefoil, Hopf link}
For the simplest unibranch plane curve singularities 
$\c_{32}=\{x^3\!=\!y^2\}$ at 
$x\!=\!0\!=\!y$ and
$\c_{52}=\{x^5\!=\!y^2\}$:
\begin{align*}
&\hat{\h}_{32}=1+qt+aq,\ \ \ \hat{\h}_{52}=1 + qt + q^2t^2 +
a(q + q^2t),\\
&H_{32}=1+qt^2+aqt,\ H_{52}=1 + qt^2 + q^2t^4 +
a(qt + q^2t^3).
\end{align*}
The corresponding $\hat{H}^i$ obviously have only (complex) 
zeros satisfying $RH$; note that   $\hat{H}_{32}^1=1$ and
$\hat{H}_{52}^1=1+q t^2$.
\vskip 0.2cm

For $2$\~branch $\c_{22}\!=\!\{(x\!+\!y)(x\!-\!y)\!=\!0\}$,
$\c_{42}\!=\!\{(x^2\!+\!y)(x^2\!-\!y)\!=\!0\}$:
\begin{align*}
& \hat{\h}_{22}=1 - t + q t +aq, \ \
\hat{\h}_{42}=1 - t + qt - qt^2 + q^2t^2 + a(q - qt + q^2t),\\
& H_{22}\!=\!1\! -\! t\! +\! q t^2 + aqt,\ \ \
H_{42}=1\!-\! t\! +\! qt^2\! -\! qt^3\! +\! q^2t^4\,+\,
a(qt\! -\! qt^2\! +\! q^2t^3).
\end{align*}
The zeros are obviously real irregular if $\om=1/q>4$
for $\hat{H}_{22}^0$ and $\hat{H}_{42}^1$. One
(real) pair of {\em irregular zeros\,}
occurs if $\om>2.25$ for $\hat{H}_{42}^0$; this pair
approaches $\{1,\om\}$ as $\om\to \infty$, which is
obvious from the formula. Otherwise their
norms are $\sqrt{\om}$.

\subsubsection{\sf Adding 
\texorpdfstring{\mathversion{normal}$Y$}{Y}, colors}
Let us provide $2$ examples in the case of the non-trivial
pairs $\{\l,'\!\!\l\}$. In the notation from
the table, they are 
$\{P(Y)(\ga[2,3](P))\}$
and $\{P(Y)(\ga[3,2](P))\}$. The corresponding
singularities are $\c_{1,23}=\{(x^3\!-\!y^2)x\!=\!0\},
\c_{1,32}=\{(x^3\!-\!y^2)y\!=\!0\}$; their links are 
$\,T(3,2)\cup \unknot\,\,$ with the linking numbers
$2,3$. One has:
{\footnotesize\small
\begin{align*}
&\hat{\h}_{1,23}=
1 - t + q t + q^2 t - q t^2 + q^2 t^2 - q^2 t^3 
+ q^3 t^3  + a^2 q^3\\ 
&+a(q + q^2 - q t + q^2 t + q^3 t - q^2 t^2 + q^3 t^2),\\
&H_{1,23}=
1 - t + q t^2 + q^2 t^3 - q t^3 + q^2 t^4 - q^2 t^5 
+ q^3 t^6  + a^2 q^3 t^3\\
&+a(q t + q^2 t^2 - q t^2 + q^2 t^3 + q^3 t^4 - q^2 t^4 + q^3 t^5),
\end{align*}
\vskip -1.0cm
\begin{align*}
&\hat{\h}_{1,32}=
1 \!-\! t \!+\! q t \!+\! q^2 t \!-\! q t^2 \!+\! q^3 t^2 \!
-\! q^2 t^3 \!+\! q^3 t^3 
\!-\! q^3 t^4 \!+\!  q^4 t^4 \\
&\!+\! a^2 (q^3 \!-\! q^3 t \!+\! q^4 t) \!+\! 
a(q \!+\! q^2 \!-\! q t \!+\! 2 q^3 t \!-\! q^2 t^2 \!
+\! q^4 t^2 \!-\! q^3 t^3 \!+\! q^4 t^3),
\\
&H_{1,32}=
1 \!-\! t \!+\! q t^2 \!+\! q^2 t^3 \!-\! q t^3 \!+\! q^3 t^5 
\!-\! q^2 t^5 \!+\! q^3 t^6 
\!-\! q^3 t^7 \!+\!  q^4 t^8 \\
&\!+\! a^2 (q^3 t^3 \!-\! q^3 t^4 \!+\! q^4 t^5) \!+\! 
a(q t \!+\! q^2 t^2 \!-\! q t^2 \!+\! 2 q^3 t^4 \!-\! q^2 t^4 
\!+\! q^4 t^6 \!-\! 
q^3 t^6 \!+\! q^4 t^7).
\end{align*}
} 
They satisfy Weak $RH$ with one pair of stable real irregular 
zeros, approaching $1,\om$ for $\om\!\to\! \infty$. For instance, 
$\hat{H}_{1,23}^1\!=\!(1\! -\! t\! +\! q t^2 ) 
(1\! +\! q t\! +\! 
q t^2 )$. 
\vskip 0.2cm

{\sf Using colors.}
The simplest colored superpolynomials is for $T(3,2)$
colored by $\om_2$:
$\hat{\h}\!=\!$ {\footnotesize
$1 + a^2\frac{q^2}{t} + q t + q t^2 + q^2 t^4 + 
a(q + \frac{q}{t} + q^2 t + q^2 t^2)$}. Accordingly,
$\hat{H}_{sym}^0\!=\!$
{\footnotesize 
$1 + q^2 t^3 + q^3 t^4 + 2 q^4 t^6 + q^5 t^8 + q^5 t^9 + 
q^6 t^{12}$}, which has $12$ $RH$\~zeros for
$\om=1/q>\varpi=1.464541725162\ldots$\,.

Let us also provide $\h^0$ for 
$\{\hat{\ga}_{2,1}(P(\yng(1,1)\,)P(\yng(2)\,)\}$=
$\{\hat{\ga}_{2,1}(P(\yng(2,2)\,)P(\yng(1)\,))\}$
(they coincide!). It is  
$1-t^2+q^2t^2-q^2t^4+q^4t^4$. The corresponding 
$\varsigma_0 t^{\pi_0}S^0$ is $-t^2(1+t^4)$, so Weak $RH$ fails
in this case with $2$ pairs of non-$RH$ zeros due to $t^2$.

\subsubsection{\sf A failure at 
\texorpdfstring{$i>0$}{zero}}
\label{SEC:2TREF}
The superpolynomials $\hat{\h}_{1,23}, \hat{\h}_{1,32}$
can be also obtained as 
$\{(\ga[2,3](P)\!\Downarrow)(Y)(P)\}$ and
$\{(\ga[3,2](P)\!\Downarrow)(Y)(P)\}$, i.e. using the
pairs $\{\l,'\!\!\l\}$
with non-trivial $'\!\l$. We use here that the DAHA construction
is isotopy-invariant. A similar one is 
$\{(\ga[2,3](P)\!\Downarrow)(Y)(\ga[3,2](P))\}$, corresponding
to $\c_{32,32}\!=\!\{(x^3\!-\!y^2)(x^3\!+\!y^2)=0\}$ with the link 
$2\,T(3,2)$ and $lk=6$. 

Transposing $3$ and $2$ in the second factor of the last
equation, the singularity $\c_{32,23}\!=\!
\{(x^3\!-\!y^2)(x^2\!-\!y^3)\!=\!0\}$ with $Z=1+2q^5+q^{10}$ 
and $lk=4$ provides a 
counterexample to Weak $RH$ with $i>0$ among uncolored
algebraic  {\em links\,}. 
Here deg$_a=3$ and the failure of $RH$ is only at $i=2$. 
This is a simplification of  
the counterexample from Section \ref{SEC:FAIL},
where the failure is at $i\!=$deg$_a\!-\!1(=\!4)$ too. 
The DAHA procedure in this case 
is $\{(\ga[2,3](P)\!\Downarrow)(Y)(\ga[2,3](P))\}$;
the superpolynomial is $\hat{\h}_{32,23}$=
{\footnotesize
\begin{align*}
&a^3 q^6+a^2 (q^3+q^4+q^5-q^3 t+q^4 t
+q^5 t+q^6 t-2 q^4 t^2+q^5 t^2+q^6 t^2-q^5 t^3+q^6 t^3)\\
&+a (q+q^2+q^3-q t+2 q^3 t+2 q^4 t+q^5 t-q^2 t^2-2 q^3 t^2
+2 q^4 t^2+2 q^5 t^2-q^3 t^3\\
&-2 q^4 t^3
+2 q^5 t^3+q^6 t^3-q^4 t^4+q^6 t^4-q^5 t^5+q^6 t^5)+
1-t+q t+q^2 t+q^3 t-q t^2\\
&+q^3 t^2+q^4 t^2-q^2 t^3-q^3 t^3+q^4 t^3+q^5 t^3
-q^3 t^4+q^5 t^4-q^4 t^5
+q^5 t^5-q^5 t^6+q^6 t^6.
\end{align*}
}  
 Then $\hat{H}_{32,23}^2\!=\!$ 
{\small $1\!-\!t\!+\!q t\!+\!q t^2\!+\!q^2 t^2\!-\!2 q t^3\!+\!
q^2 t^3\!+\!q^2 t^4
\!+\!q^3 t^4\!-\!q^2 t^5\!+\!q^3 t^5\!+\!q^3 t^6$}
has $1$ pair of real zeros 
approaching $\{1,\om\}$ for $\om>\!>0$, 
and $2$ conjugated pairs of 
complex zeros not satisfying $RH$ (though staying in the 
vicinity of $U_{\!\sqrt{\om}}$). 
One has:  $\varsigma_1 t^{\pi_2}S^2\!=\!-t(1\!+\!t^2)^2$ 
in this case; 
some ``irregular behavior" of the corresponding flagged Jacobian
factor can be expected.

We note that
unless for $(x^a-y^b)(x^b-y^a)$, the corresponding singularities
satisfy Weak $RH$ in the examples we calculated.
For instance,
$RH$ holds for all $i$ for  $\c_{34,23}\!=\!
\{(x^4\!-\!y^3)(x^2\!-\!y^3)\!=\!0\}$ with $Z=1+2q^6+2q^{12}
+q^{18}$
and the linking number $Z(1)=6$. The DAHA procedure here is
$\{(\ga[2,3](P)\!\Downarrow)(Y)(\ga[3,4](P))\}$.  


\vskip 0.2cm 
{\sf Non-algebraic knots.}
The first failures of $RH$ for $a=0$ in the family
$C\!ab(2m\!+\!1,2)T(3,2)$ are for $C\!ab(-1,2)T(3,2)$ and
 $C\!ab(-3,2)T(3,2)$. Let us provide
the corresponding $\hat{H}^0$ for the
latter:\, $1+2 q t^2+q t^3-q t^4+2 q^2 t^4+q^2
t^5+2 q^3 t^6+q^4 t^8$. Actually there are no zeros at all of
norm $\sqrt{\om}$ in this case for $\om>\!>0$. This remains equally
chaotic for all $C\!ab(-3\!-\!2m,2)T(3,2)$ as $m\ge 0$. We note
that $\hat{\h}$ become positive starting with
$C\!ab(-7,2)T(3,2)$; the corresponding $\hat{H}^0$ for $-7$ is
$1+2 q t^2+q t^3+2 q^2 t^4+q^2 t^5+q^2 t^6+2 q^3 t^6+q^3 t^7
+2 q^4 t^8+q^4 t^9+2 q^5 t^{10}+q^6 t^{12}.$ In this case,
$\varsigma_0 t^{\pi_0}S^0=t^6$, all zeros are non-$RH$ (and quite 
random).

\subsection{\bf Concluding remarks} \label{CONCL}
Let us begin with the computational aspects.
Superpolynomials have many symmetries: super-duality, evaluation
at $q=1$, color exchange, deg$_a$\~formula and more of these.
They are routinely checked by the programs that calculate
superpolynomials, including  extra evaluations
$a=-t^{n+1}$,  and it is very unlikely that there are 
mistakes with the formulas for $\hat{\h}$. 
{\em The attachment to this paper contains the formulas  
for quite a few (not all) superpolynomials used in the main 
table; 
the link is: 
{\tt \url{
http://intlpress.com/site/pub/files
/_supp/CNTP-2017-v12n3-cherednik-s1.zip
}}} .


Numerical finding the zeros of $\hat{H}^i$
is a relatively simple (and fast) task. We 
mostly rely here on the
standard software. The symmetry $\xi\mapsto \om/\xi$
provides a good independent test of the correctness of this
part of our programs. Then the program increases $\om$
to reach the $RH$\~range (if it exists) and then diminishes 
$\om$ to find the lower bonds $\varpi_i$. Then it checks
that they are (within the accuracy)
roots of the  reduced discriminants $D^i$.  
The (automated) comparison with the number of pairs of 
non-$RH$ zeros resulting from $\pi_i, S^i$ concludes
the analysis. 


\vskip 0.2cm
\subsubsection{\sf Toward Riemann's zeta}
The most optimistic expectations are that 
DAHA superpolynomials can be a move toward the Riemann 
zeta and Dirichlet $L$\~functions (and Grand $RH$). 
However,  quite a few steps are needed.

\comment{
To approach the spectral zeta and {\em Grand RH\,},
the passage to the {\em families\,}
alone seems insufficient without further efforts
of analytic nature. Given a family of plane curve singularities,
the distribution of zeros for the resulting series 
in terms of the DAHA superpolynomials is (numerically)
essentially no different from that for $\hat{H}^i$. It 
is essentially {\em uniform\,} in the corresponding
$U_\om$, though there are significant fluctuations within
some ``periods". We see here no traces of the famous
$T\log T$\~formula for the Riemann zeta.
}

\vskip 0.3cm
\vfil
{\sf Families.} First of all, $\hat{\h}$ must be extended
to the {\em families\,} of
iterated torus links; the {\em family superpolynomials\,}
$\hat{\mathbb H}(q,t,a,u)$
from Section \ref{SEC:FAM} are natural candidates (they are
actually rational functions). 
Algebraic links emerge in the DAHA theory as sequences of 
matrices $\ga\in PSL_2(\Z)$. The match of this interpretation
with the {\em splice diagrams\,} of \cite{EN} is 
a surprising outcome of \cite{ChD2}. The families are 
when we multiply one of these $\ga$ by $\tau_{\pm}^m$.
For instance,
$\{T(\rr+m\ss,\ss), m\in \Z_+\}$ and
$\{C\!ab(13+2m,2)T(3,2)\}$ for $m\in \Z_+$ are families.
 
\vskip 0.2cm
For algebraic {\em knots} and when 
$\ga_1\mapsto \tau_-^m\ga_1$, there are natural embeddings
of the corresponding rings $\r_m$. 
Geometrically, this means that we count submodules
$M$ from Section \ref{SEC:MOTZ} with some weights
in terms of its (full) ring of endomorphisms  
$\r$ from a given family. Algebraically, we sum the corresponding 
$\hat{\h}$\~polynomials over a given family with 
the weights $(u/t)^{\hbox{\tiny genus}}$ for 
uncolored algebraic knots, where $u$ is an additional parameter.

\vfil
We use that the {\em same\,} 
super-symmetry serves all rings $\r$. 
The corresponding $\hat{\mathbb H}(q,t,a,u)$ are 
generally algebraically simpler than individual 
$\hat{\h}_m(q,t,a)$ and the
$\varpi_i$ for their $a$\~coefficients are generally better 
(smaller) than those for individual  $\hat{\h}_m$ as $m>\!>0$.
Cf. formulas (\ref{3z+1}),(\ref{fam3-2}),(\ref{fam3-2}). 
 
\vskip 0.3cm
{\sf Analytic DAHA superpolynomials.}
The key step could be a passage
from algebraic superpolynomials to ``analytic" ones, 
parallel to 
Section ``Topological vertex" from \cite{ChD2}. 
As it was observed there, the analytic counterparts
of superpolynomials
for Hopf links extend (by adding $\,t$) the
Rogers-Ramanujan expansions.
The latter are interpreted in
\cite{ChF,GOW} as expansions of powers/products of 
theta-functions in terms of $q$\~Hermite and Hall-Littlewood
polynomials. They are (closely related)
limits $t\!\to\! 0,q\!\to\! 0$ of the Macdonald polynomials. 
Since the invariants $S^i,\pi_i$ we define are when $q\to 0$
(though this is not a direct substitution), one can expect
interesting connections here, which we will not discuss
in the present paper.

\vskip 0.2cm 
To define analytic DAHA superpolynomials,
we essentially replace the DAHA-Jones polynomials 
by some integrals 
of the products of the powers of the Gaussian $q^{x^2/2}$ and
their images  under the action of $\ga\in PSL(2,\Z)$.
The sums of such integrals with proper weights 
with respect to the {\em families\,} 
above generalize
the $q$\~analogs of the Riemann zeta
and Dirichlet $L$\~functions from \cite{ChZ}. For instance,
the $q$\~zeta there is the integral of $q^{x^2/2}
/(1-q^{x^2/2})$ with respect
to the Macdonald measure $\mu(X;q,t)$ in type $A_1$.
We will omit the exact definition here.
 
\vskip 0.2cm
The theory of {\em analytic\,}
DAHA superpolynomials is of clear 
independent interest regardless of zeta-functions.
Actually its main objective is in obtaining 
the invariants of Seifert and lens spaces; fruitful algebraic
applications are expected too. The details will be published
elsewhere. In contrast to knot invariants
(though these two theories are closely related), 
the invariants of Seifert
spaces are given in terms of modular functions, Maass
forms and Mock theta-functions. 
\vskip 0.2cm

An obvious problem with the passage to
the analytic superpolynomials is as follows. 
They are calculated in terms of a proper completion of
the polynomial representation $\v$ in contrast to the
algebraic theory (we present here) based on the
``adjoint representation", which is in End$_{\C}\v$
via the conjugation (actually in $\HH$).
However this sufficiently transparent relation  does not
guarantee any connection at the level of the
{\em zeros\,} of the corresponding superpolynomials.
Generally, approaches to the Riemann Hypothesis 
(Grand $RH$) via any theories of
``zeta-polynomials" satisfying $RH$, including 
the Hasse-Weil zeta functions, have little 
support in the classical and modern 
mathematics.


\subsubsection{\sf Further perspectives}
A connection of superpolynomials
with the zeta-functions
of Laplace/Dirac operators of Riemann surfaces
and $p$\~adic strings would be a fundamental development.

\smallskip
{\sf Spectral zeta-functions.}
The motivic zeta-functions are quite parallel to
the so-called {\em spectral zeta-functions\,}.
Namely, let us consider the Schottky uniformization
of Riemann surfaces and the corresponding
Dirac operators. The corresponding ``pure" zeta-functions 
then depend only on the genus in
the smooth case; see e.g. \cite{CM}. This fact (but not
the formula itself) matches Macdonald's formula, a starting point
of the Kapranov zeta-function. Then we switch to
plane curve singularities. The bound $q\le 1/2$ from
Conjecture \ref{CON:om2} resembles
the inequalities in the theory of spectral zeta.
Presumably we can arrive at the same superpolynomials
of plane curve singularities within this approach.
Importantly, the  $p$\~adic Schottky
uniformization is closely related to the curve 
singularities (as the closed fibers), which can be
potentially a tool for establishing the link with
superpolynomials.

\vskip 0.2cm
Let us mention here Witten's $p$\~adic strings,
which can be hopefully revisited and extended toward 
the superpolynomials of plain curve singularities.
At least, the {\em matrix models\,} 
can be used for this; see e.g. \cite{DMMSS}, which is
actually closely connected with the DAHA approach. 
The physics insight certainly can help here. 
\vskip 0.2cm

{\sf Adelic zeta-functions.}
Let us connect our considerations
with the classical theory
of zeta-functions of arithmetic varieties.
The compactified Jacobians and
flagged Jacobian factors can be naturally
defined over $\Z$. Accordingly, one can consider
their {\em adelic zeta-functions\,}, the
products of local zeta-functions. 
The latter are given in terms of the $q$\~coefficients 
of the motivic superpolynomials when $a=0, t=1$; 
see Section \ref{SEC:MOTZ}. If all $J_\De$ are
affine spaces, these coefficients simply give the
numbers of cells in each dimension and readily
result in the formula for the adelic zeta. It will be
the product of the corresponding
powers of the zeta-functions of affine spaces. 


Generally,  $J_\De$ are not always affine. However the 
flagged Jacobian 
factors are conjecturally {\em strongly polynomial-count\,}
due to the discussion at
the end of \cite{ChP1}. It is not impossible that they
are even paved by affine spaces (no counterexamples are known).
Thus their local zeta-functions (presumably) uniformly depend
on $|\mathbb F|$,  ignoring the points of
bad reduction (which are not a problem within a given
topological class of the singularity).  
Such adelic zeta-functions generalize
those of projective spaces, flag and Schubert varieties; 
flagged Jacobian factors can be 
naturally seen as the next level of  {\em  Schubert calculus}. 

\vskip 0.2cm

{\sf On Iwasawa polynomials.}
A similarity between the Iwasawa polynomials
and the Alexander polynomials observed by B.Mazur
is  basically as follows in our setting. 
We use that the DAHA superpolynomial
$\hat{\h}(q,t,a)$ (in the DAHA parameters) 
conjecturally coincide
with $\hat{\h}_{mot}(q,t,a)$. 
Due to  (\ref{czej}), the corresponding Alexander 
polynomial up to a normalization is 
\begin{align}\label{Iwasa}
&\hat{\h}\,(t,t,a=-1)=
\hat{\h}\,(q,t,a=-t/q)\big|_{q\mapsto t}
\\
=
L&(\Ga,q/t,t)\big|_{q\mapsto t}=L(\Ga,1,t)=
Z(\Ga,1,t)/(1-t).\notag
\end{align} 
Recall that 
$Z(\Ga,q,t)=\sum_{M\subset \r}^{\De(M)=\Ga} 
t^{dim_{\mathbb F}(\r /M)}$ (considered by 
Z\'u\oldt{n}iga-Galindo);  
i.e. the summation here is over {\em principal\,} ideals $M$.
Considering only 
the {\em group\,} of principal ideals (the generalized Jacobian)
matches the {\em group\,} of classes of ideals 
in the Iwasawa theory. Finding Iwasawa-type analogs of the whole
$\hat{\h}(q,t,a)$ (presumably coinciding with the 
reduced stable $KhR$\~polynomials of uncolored algebraic knots)
is a challenge. 

According to what we discussed above,
the passage to the {\em families\,} for 
$\ga_1\mapsto \tau_-^m\ga_1$ is natural here. 
The corresponding
Puiseux extensions play the role of Iwasawa towers.
A $u$\~counterpart of the Iwasawa polynomial is then a
weighted sum of the corresponding zeta-functions for 
principal ideals.
The limit of this construction at $q=1$ (the field
with one element) becomes the corresponding weighted sum
of Alexander polynomials, namely 
$\hat{\mathbb H}(t,t,a=-1,u)$ for 
$\hat{\mathbb H}$ in  (\ref{hsfam}). The techniques 
used to calculate the latter
allow to present $\hat{\mathbb H}(t,t,a=-1,u)$ as finite
sums of Alexander polynomials with sufficiently simple
denominators. 

\vfil
The deep connection of the 
Iwasawa polynomials with {\em $p$-adic analytic 
$L$\~functions\,}
is of obvious importance to us. See Section 7 of \cite{Mor},
especially formula (7.2) and its further discussion there.
Using {\em flags\,} and {\em families\,} (parameters $a,u$)
is beyond the approach there, and we have something else:
a connection with $q$\~zeta from \cite{ChZ}. We note that
motivic superpolynomials can be defined practically in the
same way for local $p$\~adic rings; we do not really need  
Jacobian factors to be algebraic varieties, but the
count of modules  becomes more involved. 
They may coincide with our ones 
(as in the Fundamental Lemma).

\vskip 0.2cm

{\sf What DAHA can provide.}
The coincidence of the DAHA 
superpolynomials with the motivic superpolynomials
and Galkin-St\"ohr zeta-functions can be checked as
follows. One uses the DAHA recurrence relations similar 
to those in Propositions  \ref{FAM1-1}, \ref{FAM1-m} and 
compare them with the transformations of geometric
superpolynomials under the blowups. This was checked
for some {\em families} and seems doable in general.  

For instance, this gives that the Galkin-St\"ohr zeta-functions 
depend only on the topological type of the singularity (i.e. on 
the corresponding link). For $L(\Ga,q,t)$ from (\ref{czej}) and
for any $L(\De,q,t)$ such that  
$J_\De$ is affine of the same dimension as over $\C$,
this follows from
St\"ohr's formula. However the affineness of {\em all\,} $J_{\De}$
holds  only for torus knots and some ``small" non-torus families. 
The DAHA 
superpolynomials are topological invariants, which is a 
relatively simple theorem.

Generally, the connection of the DAHA 
superpolynomials to the Khovanov-Rozansky stable polynomials 
requires the recurrence relations  for the latter
of Rosso-Jones type. They are not known, 
though the approach via Soergel modules seems quite relevant.
At $t\!=\!q$ in the DAHA parameters, the DAHA superpolynomials
were identified with the HOMFLY-PT polynomials in full generality.
The identification of the DAHA-Jones polynomials with the
corresponding {\em WRT\,} invariants was also done
in quite a few examples (including some cases of special root 
systems). The CFT approach, Rosso-Jones formulas, and the so-called
Skein are used here. 
\vskip 0.3cm

The DAHA theory and Macdonald polynomials are also connected 
with affine flag varieties, Hilbert schemes of $\C P^2$ 
(and some similar surfaces), Nekrasov's instanton
sums and the mixed Hodge polynomials 
of certain related character varieties. Linking 
these theories to (classical and motivic)
zeta-functions is quite a challenge. Also, the geometric 
superpolynomials can be expected to be connected with the
spectral zeta-functions of the plane curve singularities 
considered under the Schottky uniformization, but
this is only in the beginning and we do not see 
{\em apriori\,} reasons for DAHA to occur here.   
\vskip 0.3cm

A clear potential of the DAHA superpolynomials is their
connection with $q$\~analogs of Riemann's zeta
and $L$\~functions from  \cite{ChZ}.
Numerically, the zeros of these
$q$\~analogs are absent in the left/right half-spaces
in terms of $k\!=\!s\!-\!\frac{1}{2}$ (for $s$ from the zeta and
$t\!=\!q^k$); if true, this would 
give the Grand RH. We also suggested there a ``straight"
$q$\~Riemann hypothesis upon the symmetrization: 
Conjecture 6.3. The geometric
applications of DAHA outlined above may be not very surprising 
due to their origin: they  are deformations of Heisenberg 
and Weyl algebras.  However, their link to the classical zeta 
theory is certainly a surprising and promising development.

{\bf Acknowledgements.}
I thank Ian Philipp, David Kazhdan and the 
referees for very useful discussions and comments. Also, 
I would like to thank Kyoto University, MPI (Bonn) and
Hebrew University 
for kind invitations. I am very grateful to Yuri Manin; his works
on zeta-polynomials, iterated Shimura integrals and
our talks in Bonn inspired this project.


\vskip -2cm
\bibliographystyle{unsrt}

\medskip
\end{document}